\documentstyle[amssymb,12pt]{article}

\input amssym.def
\input amssym

\newtheorem{thm}{Theorem}[section]
\newtheorem{lem}[thm]{Lemma}
\newtheorem{prop}[thm]{Proposition}
\newtheorem{cor}[thm]{Corollary}
\newtheorem{fact}[thm]{Fact}

\newtheorem{defn}[thm]{Definition}

\newtheorem{nrmk}[thm]{Remark}
\newtheorem{expl}[thm]{Example}

\newcommand{\qed}{\hfill $\Box$ \vspace{.5cm}}
\newcommand{\pf}{{\bf Proof. }}

\title {O-minimal cohomology and definably compact definable groups}

\author {M\'{a}rio J. Edmundo \thanks{Supported by the EPSRC grant GR/M66332.}
\\The Mathematical Institute
\\24-29 St Giles
\\OX1 3LB  Oxford, U.K
\\edmundo@maths.ox.ac.uk
\\
\\}
\date{December 16, 2000}

\newcommand{\into}{\longrightarrow}

\renewcommand{\hat}{\widehat}
\renewcommand{\tilde}{\widetilde}
\renewcommand{\bar}{\overline}

\newcommand{\NN}{\mathbb{N}}
\newcommand{\ZZ}{\mathbb{Z}}
\newcommand{\QQ}{\mathbb{Q}}
\newcommand{\RR}{\mathbb{R}}
\newcommand{\CC}{\mathbb{C}}

\newcommand{\Ll}{\mbox{$\cal L$}}
\newcommand{\M}{\mbox{$\cal M$}}
\newcommand{\N}{\mbox{$\cal N$}}

\newcommand{\Pp}{\mbox{$\cal P$}}
\begin{document}

\maketitle
\begin{abstract}
Let $\N$ an
o-minimal expansion of a real closed field.
We develop the cohomology
theory for the category of $\N$-definable manifolds with
continuous $\N$-definable maps and use this to solve the
Peterzil-Steinhorn problem \cite{ps} on the
existence of torsion points on $\N$-definably compact $\N$-definable
abelian groups. Namely we prove the following result:
Let $G$ be an $\N$-definably compact $\N$-definably connected $\N$-definable
group of dimension $n$. Then the o-minimal Euler characteristic of $G$ is
zero. Moreover, if $G$ is abelian then $\pi _1(G)=\ZZ $$^n$ and 
for each $k>1$, the subgroup $G[k]$ of $k$-torsion points of $G$ is
isomorphic to $(\ZZ$$/k\ZZ$$)^n$.
We also compute the cohomology rings of $\N$-definably compact $\N$-definably
connected $\N$-definable groups.
\end{abstract}

\newpage

\begin{section}{Introduction}
\label{section introduction}

We work over an  o-minimal expansion $\N$$=(N,0,1,<,+,\cdot ,\dots)$
of a real closed field. Definable means $\N$-definable (possibly with 
parameters). The results of this paper are an extension to the work of
A. Woerheide in \cite{Wo} where o-minimal homology was
introduced. We develop the cohomology theory for the category of
definable manifolds with continuous definable maps and use this to
solve the problem of existence of torsion points on definably compact 
definable abelian groups (see \cite{ps}). Namely we prove the
following result (see theorem \ref{thm euler characteristic} and 
theorem \ref{thm main}).
  
\begin{fact}
Let $G$ be a definably compact definably connected definable
group of dimension $n$. Then the o-minimal Euler characteristic of $G$ is
zero. Moreover, if $G$ is abelian then $\pi _1(G)\simeq \ZZ $$^n$ and 
for each $k>1$, the subgroup $G[k]$ of $k$-torsion points of $G$ is
isomorphic to $(\ZZ$$/k\ZZ$$)^n$.
\end{fact}

We also compute the cohomology rings of definable groups (see theorem 
\ref{thm main}).
Other main results of this paper are: Poincar\'e duality (theorem
\ref{thm poincare duality}), Alexander duality (theorem \ref{thm
alexander duality}), Lefschetz duality (theorem \ref{thm lefschetz duality}),
the Lefschetz fixed point theorem (theorem
\ref{thm lefschetz}) in the category of Hausdorff 
definable manifolds with continuous definable maps. From this
we obtain the usual corollaries, in particular the generalised
Jordan-Brouwer separation theorem (corollary \ref{cor general
separation thm}) for Hausdorff definable manifolds.
Note that the Lefschetz fixed point theorem and  the generalised
Jordan-Brouwer separation theorem are also proved in \cite{bo1} for 
the category 
of affine $C^p$ definable manifolds with (definable) $C^p$  definable maps
with $p\geq 3$. The method used in \cite{bo1} is based on 
intersection theory for such definable manifolds adapted from classical 
differential topology. 

\medskip
The main motivation for the work presented here is the classification
problem for groups definable in o-minimal structures. The results from
\cite{e1} together with those from \cite{pps1} reduce this problem to
the classification problem for definably compact definable abelian
groups in o-minimal structures. In \cite{e2} we try to classify
definably compact definable abelian groups $G$ under the assumption that
$\pi _1(G)\simeq \ZZ$$^{dimG}$. As we mentioned above this assumption
is true when we are working in an o-minimal expansion of a real closed
field, we believe that this remains true in general, although to prove
it one would need to develop possibly an o-minimal analogue of
``Grothendieck theory of duality''.

\medskip
Most of the results of this paper are modifications of classical
results in algebraic topology (see for example \cite{d}, \cite{f},
\cite{g} or \cite{ro}). For the convenience of model theorists and
since the setting is different we will be careful enough to include
all the details. In particular, in order to introduce terminology and
notation we include a brief description of Woerheide work.
We will follow L. van den Dries book \cite{vdd}
for results on o-minimal structures and for the
non model theorists we recall in section \ref{section o-minimal structures}
some basic notions of o-minimality. 

\end{section}

\begin{section}{O-minimal structures}
\label{section o-minimal structures}

\begin{defn}\label{defn structure}
{\em
A {\it structure} $\N$ consists is a non empty set $N$ together with: $(i)$ a
set of constants $(c_k^N)_{k\in K}$, where $c_k^N\in N$; $(ii)$ a family
of maps $(f_j^N)_{j\in J}$, where $f_j^N$ is an $n_j$-ary map,
$f_j^N:N^{n_j}\into N$ and  $(iii)$ a family $(R_i^N)_{i\in I}$ of 
relations that is, for each $i$, $R_i^N$ is a subset of
$N^{n_i}$ for some $n_i\geq 1$. The {\it language} $\Ll$
associated to a structure $\N$ consists of:  $(i)$ for each constant
$c_k^N$, a constant symbol $c_k$; $(ii)$ for each map $f_j^N$ a
function symbol, $f_j$ of arity $n_j$ and  $(iii)$ for each relation
$R_i^N$, a relation symbol, $R_i$ of arity $n_i$.
}
\end{defn}

We often use the following notation
$\N$$=(N, (c_k^N)_{k\in K},(f_j^N)_{j\in J}, (R_i^N)_{i\in I})$, and 
sometimes we omit the superscripts.  Similarly we use the notation
$\Ll$$=\{ (c_k)_{k\in K}, (f_j)_{j\in J}, (R_i)_{i\in I}\}$. 
If $\Ll$ is the language associated with the structure $\N$ we say
that $\N$ is an
{\it $\Ll$-structure}. If $\Ll$$'\subseteq \Ll$ are two languages and
$\N '$ and $\N$ are respectively an $\Ll$$'$-structure and an
$\Ll$-structure, such that $N'=N$ then we say that $\N$ is an 
{\it expansion} of $\N '$ or that $\N '$ is a {\it reduct} of $\N$.

\medskip
Let $\Ll$ be a language and $\N$ an $\Ll$-structure. We are going to
define inductively the set of $\Ll$-formulas and satisfaction of an 
$\Ll$-formula $\phi $ in $\N$, in order to define the $\N$-definable sets. 

\begin{defn}\label{defn formulas}
{\em
The set of
{\it $\Ll$-terms} is generated inductively by the following rules: $(i)$
every variable $x_i$ from a countable set of variable $(x_q)_{q\in Q}$
is an $\Ll$-term, $(ii)$ every constant of $\Ll$ is an
$\Ll$-term and $(iii)$ if $f$ is in $\Ll$ is an $n$-ary function, and 
$t_1,\dots ,t_n$ are $\Ll$-terms, then $f(t_1,\dots ,t_n)$ is an
$\Ll$-term. An {\it atomic} $\Ll$-formula is an
expression of the form: $t_1=t_2$ or $R(t_1,\dots ,t_n)$ where $R$ is
an $n$-ary relation in $\Ll$ and $t_1,\dots ,t_n$ are $\Ll$-terms. We 
sometimes write $R(t_1(x_1,\dots ,x_k),\dots ,t_n(x_1,\dots ,x_k))$ if
we want to explicitly show the variables occurring in the atomic
$\Ll$-formula. Given a tuple $a\in N^k$, we say that $a$ {\it
satisfies} $R(t_1(x),\dots ,t_n(x))$ in $\N$ where $x=(x_1,\dots
,x_k)$ if $R^N(t_1(a),\dots ,t_n(a))$ holds. We denote this by 
$\N$$\models R(t_1(a),\dots ,t_n(a))$.

The set of $\Ll$-formulas is generated inductively by the following
rules: $(i)$ all atomic $\Ll$-formulas are $\Ll$-formulas; $(ii)$ if $\phi
_1(x)$ and $\phi _2(x)$ are $\Ll$-formulas, then $(\phi _1\wedge \phi
_2)(x)$ and $(\phi _1\vee \phi _2)(x)$  are $\Ll$-formulas, and for
$a\in N^k$, $\N$$ \models (\phi _1\wedge \phi _2)(a)$ iff $\N$$\models
\phi _1(a)$ and $\N$$\models \phi _2(a)$, we also have the obvious
clause for $\vee $; $(iii)$ if $\phi (x)$ is an $\Ll$-formula, $\neg
\phi (x)$ is an $\Ll$-formula, and for $a\in N^k$, $\N$$ \models  \neg
\phi (a)$ iff $\phi (a)$ does not hold in $\N$ (this is denote by $\N$
$\nvDash \phi (a)$); $(iv)$ if $\phi (x,x_{k+1})$ is an $\Ll$-formula,
then $\exists x_{k+1}\phi (x,x_{k+1})$ is an $\Ll$-formula, and for  
$a\in N^k$, $\N$$\models \exists x_{k+1} \phi (a,x_k)$ iff there
exists $b\in N$ such that $\N$$\models \phi (a,b)$; $(v)$ the obvious 
clauses  for $\forall $.
}
\end{defn}
 
We are now ready to define the important notion of $\N$-definable sets
and $\N$-definable maps.
 
\begin{defn}\label{defn definable sets}
{\em 
A subset $D\subseteq N^k$ is an {\it $\N$-definable subset (defined
over $A\subseteq N$)} if there is an $\Ll$-formula $\phi (x,y)$ with 
$x=(x_1,\dots ,x_k)$ and $y=(x_{k+1},\dots ,x_{k+m})$ and some $b\in
A^m$ such that $D=\{a\in N^k:\N$$\models  \phi (a,b) \}$. 
If $X\subseteq N^k$ and $Y\subseteq N^m$ are $\N$-definable sets (over 
$A\subseteq N$), a function $f:X\into Y$ is $\N$-definable (over $A$)
if its graph is an $\N$-definable set (over $A$). More generally, a 
structure $\M$$=(M, (c_k^M)_{k\in K}, (f_j^M)_{j\in J}, (R_i^M)_{i\in
I})$ is $\N$-definable (over $A$) if: $(i)$ $M\subseteq N^l$ is $\N$-
definable (over $A$); $(ii)$ for each $k\in K$ there is  a point
$m_k\in M$ corresponding to $c_k^M$; $(iii)$ for each $j\in J$ the
function $f_j^M:M^{n_j}\into M$ is $\N$-definable (over $A$) and
$(iv)$ for each $i\in I$ the relation $R_i^M \subseteq M^{n_i}$ is
$\N$-definable (over $A$). Note that, in this case every
$\M$-definable set is also an $\N$-definable set. 
}
\end{defn}

Finally we include here one more definition from basic model theory
which palys a crucial role in our paper.

\begin{defn}\label{defn elementary}
{\em
Given two $\Ll$-structures $\N$ and $\M$, we say that $\N$ is an 
{\it $\Ll$-substructure} of $\M$, denoted by $\N$$\subseteq \M$, if 
$N\subseteq M$ and: $(i)$  for every constant symbol $c$ in $\Ll$,
$c^N=c^M$; $(ii)$ for every $n$-ary function symbol $f$ in $\Ll$, for
every $a\in N^n$, $f^N(a)=f^M(a)$ and $(iii)$ for every $n$-ary 
relation symbol $R$ in $\Ll$ and for every $a\in N^n$, $R^N(a)$ iff
$R^M(a)$.  Let $\N$$\subseteq \M$. We say that
$\M$ is an {\it elementary extension} of $\N$ (or that $\N$ is an {\it 
elementary substructure} of $\M$), denoted by $\N$$\preceq \M$, if for
every $\Ll$-formula $\phi (x)$, for all $a\in N^k$, we have
$\N$$\models \phi (a)$ iff  $\M$$\models \phi (a)$. 
}
\end{defn}

Note that (by Tarski-Vaught test) $\N$$\preceq \M$ iff for every non empty
$\M$-definable set $E\subseteq M^l$, defined with parameters from $N$, 
$E(N):=E\cap N^l$ ("the set of $N$-points of $E$") is a non empty
$\N$-definable set. Clearly, if $S\subseteq N^l$ is an $\N$-definable
set defined with parameters from $N$ and $\N$$\preceq \M$, then the
$\Ll$-formula which determines $S$ determines an $\M$-definable set 
$S(M)\subseteq M^l$ ("the $M$-points of $S$"). 

\medskip
The {\it theory} $Th(\N$$)$ of an $\Ll$-structure $\N$ is the
collection of all $\Ll$-sentences (i.e., $\Ll$-formulas without free
variables) $\sigma $ such that $\N$$ \models  \sigma $. $\N$ is {\it 
elementarily equivalent} to $\M$, denoted $\N$$\equiv \M$ iff
$Th(\N$$)=Th(\M$$)$. Clearly, if $\N$$\preceq \M$ then $\N$$\equiv
\M$. The following two  facts (the L\"{o}wenheim-Skolem theorems) are
fundamental theorems of basic model theory: $(1)$ if 
$\N$ an $\Ll$-structure with $X\subseteq N$, 
then for every cardinal $\kappa $ such that $|X|+|\Ll$$|\leq \kappa
\leq |N|$, $\N$ has an elementary substructure $\M$ such that
$X\subseteq M$ and $|M|=\kappa $; and $(2)$ if $\N$ be an infinite
$\Ll$-structure, then for any cardinal $\kappa >|N|$, $\N$ has an
elementary extension of cardinality $\kappa $. 

\begin{defn}\label{defn o-minimal structure}
{\em
An {\it o-minimal structure} is an expansion $\N$$=(N,<,\dots )$ of a 
linearly ordered nonempty set $(N,<)$, such that every $\N$-definable
subset of $N$ is a finite union of points and intervals with endpoints in 
$N\cup \{-\infty, +\infty \}$. 
}
\end{defn}

Note the following important results:
let $\N$ be an o-minimal structures and $\Ll$ its language. 
Then: $(1)$ every $\N$-definable
structure $\M$ which is an expansion of a linearly ordered nonempty set
$(M,<_M)$ is also o-minimal; $(2)$ \cite{KPS} if $\M$ is a structure (in the
language of $\N$) such that $\N$$\equiv \M$ then $\M$ is also o-minimal; $(3)$
\cite{PiS1} for every $A\subseteq N$ there is a {\it prime model} of $Th(\N$$)$
over $A$ i.e., there is an o-minimal structure $\Pp$ such that
$A\subseteq P$ and for all $\M$ with $A\subseteq M$ and either $\M$$\preceq \N$
or $\N$$\preceq\M$ we have $\Pp$$\preceq \M$; and $(4)$ for every 
$\kappa >\max\{\aleph _0, |\Ll$$|\}$ there up to isomorphism
$2^{\kappa }$ o-minimal structures $\M$ such that $|M|=\kappa $ and 
$\M$$\equiv \N$ (see \cite{Sh}), and
if $\Ll$ is countable then up to isomorphism there are either
$2^{\aleph _0}$ or
$6^n3^m$ countable o-minimal structures $\M$ such that $\M$$\equiv \N$
(see \cite{M}).

\medskip
There are many geometric properties of $\N$-definable sets and
$\N$-definable maps in an o-minimal structures $\N$. For example, 
two of the most powerful results are the monotonicity theorem for 
definable one variable functions and the $C^p$-cell decomposition 
theorem for definable sets and definable maps. We will now explain the
$C^p$-cell decomposition theorem (here $p=0$ if $\N$ is not an
expansion of a (real closed) field) in order to introduce the notion of
o-minimal dimension and Euler characteristic.

\begin{defn}\label{defn cells and dimension}
{\em
$C^p$-cells and o-minimal dimension are defined inductively as
follows: $(i)$ the unique non empty $\N$-definable subset of $N^0$ is a 
$C^p$-cell of dimension zero, a point in $N^1$ is a $C^p$-cell of
dimension zero and an open interval in $N^1$ is a $C^p$-cell of
dimension one; $(ii)$ a $C^p$-cell in $N^{l+1}$ of dimension $k$
(resp., $k+1$) is an $\N$-definable set of the form $\Gamma
(f)$ (the graph of $f$) where $f:C\into N$ is a $C^p$-definable
function and $C$ is a $C^p$-cell in $N^l$ of dimension $k$ (resp., of
the form $(f,g)_C:=\{(x,y)\in C\times N:f(x)<y<g(x)\}$ where
$f,g:C\into N$ are $C^p$-definable function with $-\infty \leq f<g\leq
+\infty $ and $C$ is a $C^p$-cell in $N^l$ of dimension $k$.  The
Euler characteristic $E(C)$ of a $C^p$-cell $C$ of dimension $k$ is 
defined to be $(-1)^k$.
}
\end{defn}

\begin{defn}\label{defn cell decomposition}
{\em
A $C^p$-cell decomposition of $N^m$ is a special kind of partition of
$N^m$ into finitely many $C^p$-cells: a partition of $N^1$ into
finitely many disjoint $C^p$-cells of dimension zero and one is a
$C^p$-cell decomposition of $N^m$ and, a partition of $N^{k+1}$ into
finitely many disjoint $C^p$-cells $C_1, \dots , C_m$ is a $C^p$-cell 
decomposition of $N^{k+1}$ if $\pi (C_1), \dots , \pi (C_m)$ is a
$C^p$-cell decomposition of $N^k$ (where $\pi:N^{k+1}\into N^k$ is the
projection map onto the first $k$ coordinates). Let $A_1,\dots
A_k\subseteq A\subseteq N^m$ be $\N$-definable sets. A $C^p$-cell
decomposition of $A$ compatible with $A_1,\dots A_k$  is a finite
collection $C_1, \dots , C_l$ of $C^p$ partitioning $A$ obtained from
a $C^p$-cell decomposition of $N^m$ such that for every $(i, j)\in
\{1, \dots ,k\}\times \{1, \dots , l\}$ if $C_j\cap A_i\neq \emptyset $ 
then $C_j\subseteq A_i$. 
}
\end{defn}

\begin{fact}\label{fact cell decomposition}
{\it ($C^p$-cell decomposition theorem)}
Given $\N$-definable sets $A_1,\dots ,A_k$ $\subseteq A\subseteq N^m$ there
is a $C^p$-cell decomposition of $A$ compatible with $A_1,\dots A_k$
and, for every $\N$-definable function $f:A\into N$, $A\subseteq N^m$,
there is a $C^p$-cell decomposition of $A$, such that each restriction 
$f_{|C}:C\into N$ is $C^p$ for each cell $C\subseteq A$ of the
$C^p$-cell decomposition.
\end{fact}

The o-minimal dimension $dim(X)$ and Euler characteristic $E(X)$ of a
$\N$-definable set $X$ are defined by $dim(X)=\max \{dim (C):C\in
\mathcal{C}$$\}$ and $E(X)=\sum _{C\in \mathcal{C}}$$E(C)$ where
$\mathcal{C}$ is some (equivalently any) $C^p$-cell decomposition of
$X$. These notions are well behaved under the usual
set theoretic operations on $\N$-definable sets, are invariant under
$\N$-definable bijections and given an $\N$-definable family of $\N$-
definable sets, the set of parameters whose fibre in the family has a 
fixed dimension (resp., Euler characteristic) is also an $\N$-definable set.
The cell decomposition theorem is also used to
show that every $\N$-definable set has only finitely many $\N$-definably
connected components, and given an $\N$-definable family of $\N$-definable sets
there is a uniform bound on the number of $\N$-definably connected
components of the fibres in the family. 

\medskip
Since in this
paper we are concerned only with o-minimal expansions of ordered
fields (necessarily real closed fields) {\it we will from now on assume
that $\N$ is an o-minimal expansion of a real closed field and
definable means $\N$-definable}. Below we list the some of properties of
definable sets and definable maps that we will be using through this
paper.

\begin{defn}\label{defn triangulaton}
{\em
Let $S_1,\dots , S_k$ $\subseteq S\subseteq N^m$ be definable sets. A
{\it definable triangulation} in $N^m$ of $S$ compactible with $S_1,\dots ,
S_k$ is a pair $(\Phi ,K)$ consisting of a complex $K$ in $N^m$ and a
definable homeomorphism $\Phi :S\into |K|$ such that each $S_i$ is a
union of elements of $\Phi ^{-1}(K)$. We say that $(\Phi ,K)$ is a
{\it stratified definable triangulation} of $S$ compactible with 
$S_1,\dots , S_k$ if: $m=0$ or $m>0$ and there is a stratified
definable triangulation $(\Psi ,L)$ of $\pi (S)$ compactible with $\pi
(S_1),\dots , \pi (S_k)$ (where $\pi :N^m\into N^{m-1}$ is the
projection onto the first $m-1$ coordinates) such that $\pi
_{|Vert(K)}:K\into L$ is a simplicial map and the diagram 
\[
\begin{array}{clcr}
S\,\,\,\stackrel{\Phi }{\rightarrow}\,\,\,|K|\\
\,\,\,\,\,\,\,\,{\downarrow}^{\pi }\,\,\,\,\,\,\,\,\,\,\,\,\,\, 
{\downarrow}^{\pi }\,\,\,\,\,\,\,\,\,\\
\pi (S)\stackrel{\Psi }{\rightarrow}|L|.\\
\end{array}
\]
We say that $(\Phi ,K)$ is a {\it quasi-stratified definable
triangulation} of $S$ compatible with $S_1,\dots , S_k$ if there is a
linear bijection $\alpha :N^m\into N^m$ such that $(\alpha \Phi \alpha
^{-1}, \alpha K)$ is a stratified definable triangulation of $\alpha
(S)$ compatible with $\alpha (S_1),\dots , \alpha (S_k)$.
}
\end{defn}

\begin{fact}\label{fact definable triangulation theorem}
\cite{vdd} (Definable triangulation theorem).
Let $S_1,\dots , S_k$ $\subseteq S\subseteq N^m$ be definable
sets. Then, there is a definable triangulation of $S$ compatible with 
$S_1,\dots , S_k$. Moreover, if $S$ is bounded then, there is a
quasi-stratified definable triangulation of $S$ compatible with 
$S_1,\dots , S_k$.
\end{fact}

Several other geometric properties from semialgebraic and subanalytic 
geometry also hold for definable sets and definable maps: 
we have a definable curve selection theorem, a 
definable trivialization theorem and finally (see \cite{DM2}) a
definable analogue of the uniform bounds on growths 
theorem, the $C^p$-multiplier theorem, the generalised Lojasiewicz 
inequality, the $C^p$ zero set theorem, the $C^p$ Whitney
stratification theorem, etc.

\end{section}

\begin{section}{Definable manifolds}
\label{section definable manifolds}

\begin{subsection}{Definable manifolds}
\label{subsection definable manifolds}
              
\begin{defn}\label{defn definable n spaces}
{\em                  
A {\it definable manifold} (of dimension $m$ defined over $A$) 
is a triple $\mathbf{X}$
$:=(X,(X_i,\phi _i)_{i\in I})$ where $\{X_i:i\in I\}$ is a finite
cover of the set $X$ and for each $i\in I$, (1) we have injective maps 
$\phi _i:X_i\into N^m$ such that $\phi _i(X_i)$ is an open 
definably connected definable set (defined over $A$); (2) each 
$\phi _i(X_i\cap X_j)$ is an open definable subset of  $\phi_i(X_i)$
(defined over $A$)
and (3) the map $\phi _{ij}:\phi _i(X_i\cap X_j)\into \phi _j(X_i\cap
X_j)$ given by
$\phi _{ij}:=\phi _j\circ \phi _i^{-1}$ is a definable homeomorphism
(defined over $A$) for all $j\in I$ such that $X_i\cap X_j\neq \emptyset $. 
}  
\end{defn}

\begin{defn}
{\em
Let $\mathbf{X}$$=(X,(X_i,\phi _i)_{i\in I})$ and $\mathbf{Y}$
$=(Y,(Y_j,\psi _j)_{j\in J})$ be definable manifolds. 
A {\it definable subset} of $\mathbf{X}$
is a set $Z\subseteq X$ such that the sets $\phi _i(Z\cap X_i)$ are definable. 
If all of $\phi _{i_l}(Z\cap X_{i_l})$ are definable over $B$ we say
that $Z$ is definable over $B$. Let $Z$ be a definable subset of $\mathbf{X}$.
A map $f:Z\subseteq X\into Y$ is a {\it definable map} if
$\psi _j\circ f\circ \phi _i^{-1}:
\phi _i(X_i\cap Z)\into \psi _j(Y_j)$ is a  definable map
whenever it is defined. 
}
\end{defn}

A definable manifold $\mathbf{X}$ has a natural topology generated by the
definable open subsets of $X$ i.e., those definable subsets
$U$ such that for all $i\in I$, $\phi _i(U\cap X_i)$ is an open
definable subset of $\phi _i(X_i)$.
We say that a definable manifold $\mathbf{X}$ is {\it definably
connected} if there are no two disjoint definable open subsets of $X$
whose union is $X$. By \cite{e2}, $\mathbf{X}$ is definably connected 
iff it is definably path connected i.e., for all $x,y\in X$ there a definable 
continuous map $\alpha :[0,1]\into X$ such that $\alpha (0)=x$ and
$\alpha (1)=y$.

\begin{defn}\label{defn submanifold}
{\em
A {\it definable submanifold} of $\mathbf{X}$ is a definable manifold 
$\mathbf{Y}$ such that $Y\subseteq X$and each $\psi _j:Y_j\into \psi
_j(Y_j)$ are definable and whose topology is the
induced topology from $X$. The definable submanifolds of $N^m$ are
called {\it affine definable manifolds}. 
}
\end{defn}

\begin{nrmk}
{\em
Its easy to see that definable manifolds are $T1$, but one can construct
easily an example of a definable
manifold $\mathbf{X}$ which is not Hausdorff: Let $X_1:=N\times
\{0\}$, $X_2:=N\times \{1\}$ with
the subsets $X_1^{-}:=N^{<0}\times \{0\}$ and $X_2^{-}:=N^{<0}\times
\{1\}$ identified. But affine definable manifolds are definably normal
\cite{vdd} and by lemma \ref{lem regular} below
every Hausdorff definable manifold $\mathbf{X}$ is definably regular. 
By \cite{vdd} 
every definably regular definable manifold is definably homeomorphic to an
affine definable manifold. 
}
\end{nrmk}

\begin{lem}\label{lem regular}
Every Hausdorff definable manifold $\mathbf{X}$ is definably regular.
\end{lem}

\pf
This is contained in the proof of lemma 10.4 \cite{bo1}: For each
$i\in I$ and $x,y\in X_i$, let $d_i(x,y):= |\phi _i(x) - \phi _i(y)|$.  
Let $K$ be a closed  definable subset of $X$ and $a_0\in X\setminus K$. 
Let $I^K:=\{i\in I:K\cap X_i\neq \emptyset \}$, 
$I_{a_0}:=\{j\in I:a_0\in X_i\}$ and for each $i\in I^K$ let 
$K_i:=K\cap X_i$. 
For $\epsilon \in N$, $\epsilon >0$ and $i\in I^K$ define
$K_i^{\epsilon }$ to be the set of points $y\in \bigcup _{j\in
I}X_j$ such that there is a point $x$ in $K_i \cap \bigcup_{j\in
I}X_j$ with $d_j(x,y) < \epsilon $. Since $I$ is finite, 
$K_i^{\epsilon }$ is open (and contains $K_i$).  Let $K^{\epsilon
}:=\bigcup _{i\in I^K}K_i^{\epsilon }$. Then $K^{\epsilon }$ is an
open definable subset containing $K$. 
Similarly we define $L^{\epsilon }$ containing $a_0$ to be the open
definable subset of all $y\in \bigcup _{j\in I_{a_0}}X_j$ such that
there is $x\in \bigcup _{j\in I_{a_0}}X_j$ with $d_j(x,y)<\epsilon $. 

If for some $\epsilon > 0\;$ $K^{\epsilon }\cap L^{\epsilon
}=\emptyset $ we are done. Otherwise, there is a finite subset $J$ of
$I^K$ such that $K_J^{\epsilon }\cap L^{\epsilon }\neq \emptyset $ for
all sufficiently small $\epsilon >0$, where $K_J^{\epsilon }:=\bigcup
_{i\in J}K_i^{\epsilon }$.
Now by definable choice (chapter 6, proposition 1.2 \cite{vdd})  and
lemma 10.3 \cite{bo1}, there is a definable continuous map $a \colon
(0,\epsilon ) \to X$ such that $a(\epsilon )\in K_J^{\epsilon } \cap 
L^{\epsilon }$. Since ${\mathbf X}$ is Hausdorff $a_0 = \lim_{\epsilon
\to 0} 
a(\epsilon)$.  We reach a contradiction by showing that $a_0 \in K$.
Choose 
$i$ such that $a_0 \in X_i$. Then for all sufficiently small $\epsilon
> 0$ 
we have $a(\epsilon )\in X_i$ so $d_i(a(\epsilon ), K \cap X_i)$ is
well 
defined and must be less than $\epsilon $ since $a(\epsilon )$ belongs
to 
$K_J^{\epsilon }$. Therefore, $\lim_{\epsilon \to 0} d_i(a(\epsilon ), 
K \cap X_i) =0$ i.e., $d_i(a_0, K \cap X_i) =0$ and $a_0 \in K$.
\qed

Note that if $\mathbf{X}$ (resp., $\mathcal{X}$) is a definable
manifold over $A$ (resp., a collection of definable manifold
over $A$ and definable maps over $A$ between definable 
manifolds in $\mathcal{X}$) then $\mathbf{X}$ (resp., $\mathcal{X}$) 
can be identified with a definable structure over $A$. Moreover, 
if $\M$ is a structure in the language of $\N$ such that $A\subseteq M$
and either$\N$$\preceq \M$ or $\M$$\preceq \N$ then $\mathbf{X}$ (resp., 
$\mathcal{X}$) can also be seen as an $\M$-definable manifold over
$A$ (resp., a collection of $\M$-definable manifold over $A$ and 
$\M$-definable maps over $A$ between $\M$-definable manifolds in 
$\mathcal{X}$). 

Although many topological notions (such has Hausdorff, regular,
definably connected, etc.), of definable manifolds and
their definable maps are invariant under the process mentioned
above, other notions such as compact and locally compact are not good
notions in our context since, for example if $N$ is non standard then 
no nontrivial definable manifold will be locally compact or
compact. For this reason a weaker notion of compactness for 
definable manifold was introduced in \cite{ps}: 

\begin{defn}\label{defn definably compact}
{\em
A definable manifold $\mathbf{X}$ is {\em definably compact} if it is
Hausdorff and for every definable curve $\sigma\colon (a,b)\subseteq
N\to X$, where $-\infty\leq a<b\leq +\infty$ there are
$c,d\in X$ such that $\lim_{x\to a^+}\sigma(x)=c$ and
$\lim_{x\to b^-}\sigma(x)=d$. 
}
\end{defn}

By \cite{ps}, an affine definable manifold $\mathbf{X}$ is definably
compact iff $X$ is a closed and bounded definable subset of some
$N^m$. Therefore, this notion of definable compactness is invariant under
elementary extension/substructure.

\begin{expl}\label{expl definable groups}
{\em
Let $H:=(H,\cdot )$ be a definable a group. 
Results from \cite{p} and \cite{pst}
show that: (1) $H$ has a unique structure of a 
definable manifold such that the group operations are definable
continuous maps. 
Note that,  since definable groups are Hausdorff (as a definable manifolds), 
every definable group is isomorphic to an affine definable group.
Moreover,
there is in a uniformly definable family of subsets of $H$ containing the
identity element $e$, $\{V_a:a\in S\}$
such that $\{V_a: a\in S\}$ is a basis for the open
neighbourhoods of $e$; (2) 
the topology of a definable subgroup  $G$
of $H$ agrees with the topology induced on $G$ by $H$, $G$ is closed in $H$
and if $dim G=dim H$ then $G$ is open in $H$; and (3) 
a definable homomorphism of definable groups $\alpha :H\into K$ is a
continuous (in fact $C^p$ for every $p$) definable map.
}
\end{expl}

\end{subsection}

\begin{subsection}{Definable manifolds with boundary}
\label{subsection definable manifolds with boundary}

\begin{defn}\label{defn nspaces with boundary}
{\em                  
A {\it definable manifold with boundary} of dimension $m$ is a
triple $\mathbf{X}$
$:=(X,(X_i,\phi _i)_{i\in I})$ where
$\{X_i:i\in I\}$ is a finite cover of the set $X$ and for each $i\in I$, (1) we
have maps $\phi _i:X_i\into \{(x_1,\dots ,x_m)\in N^m:x_m\geq 0\}$ such that
$\phi _i(X_i)$ is a definable open 
definably connected set; (2) each $\phi _i(X_i\cap X_j)$
is definable and open in $\phi_i(X_i)$;
(3) the map $\phi _{ij}:\phi _i(X_i\cap X_j)\into \phi _j(X_i\cap
X_j)$ given by
$\phi _{ij}:=\phi _j\circ \phi _i^{-1}$ is a definable homeomorphism
for all $j\in I$ such that $X_i\cap X_j\neq \emptyset $ and (4), there
is $i\in I$ such that $\{x\in X_i:\phi _i(x)\in 
N^{m-1}\times \{0\}\}\neq \emptyset $. 
}
\end{defn}

Just like for definable manifolds, we have the notions of
definable subset, open definable subset, definable submanifold and
definable maps etc., for
definable manifolds with boundary. 

\begin{defn}\label{defn the boundary}
{\em
We have the following definable submanifolds of $\mathbf{X}$ of
dimension $m$ and $m-1$ respectively:
$\mbox{\textbf{\. X}}$$:=(\mbox{\. X},(\mbox{\. X}_i,\phi
_{i|})_{i\in I})$ where $\mbox{\. X}:=\bigcup _{i\in
I}\mbox{\. X}_i$ and $\mbox{\. X}_i:=\{x\in X_i:\phi _i(x)\in 
\{(x_1,\dots ,x_m)\in N^m:x_m>0\} \}$; and $\partial
\mathbf{X}$$:=(\partial X, (\partial X_i, \phi _{i|})_{i\in I})$ where 
$\partial X:=\bigcup _{i\in I}\partial X_i$ and 
$\partial X_i:=\{x\in X_i:\phi _i(x)\in N^{m-1}\times \{0\} \}$. 
$\partial \mathbf{X}$ is called the {\it boundary of} $\mathbf{X}$.
}
\end{defn}

We include here the following remark which will be use in subsection 
\ref{subsection lefschetz duality theorem}.

\begin{nrmk}\label{nrmk with boundary}
{\em
Then ${\mathbf X}$ can be definably embedded in two different ways
in a definable manifold $2{\mathbf X}$ of dimension $m$ such that 
$X$ is a closed definable submanifold, and there are open definable 
submanifolds $Y_1$ and $Y_2$ of $Y$ 
each containing a copies $X_1$ and $X_2$ of $X$ and such that $X_i$ is
a definable deformation retract of
$Y_i$ and $Y_1\cap Y_2$ is definably homotopically equivalent to 
$\partial X$.

In fact let $2{\mathbf X}$$=(Y,(Y_i,\psi _i)_{i\in I})$ be given by: 
$Y_i=\mbox{\. X}^1_i\cup \partial X_i\times [0,-1]
\cup \mbox{\. X}^2_i$ where $\mbox{\. X}^1_i$ and $\mbox{\. X}^2_i$
are 
copies of $\mbox{\. X}_i$ and $\psi _i$ is given by
$\psi _i=\phi _i$ on $\mbox{\. X}^1_i$, $\psi _i=\phi _i\times
1_{[0,-1]}$ on  
$\partial X_i\times [0,-1]$ and
for $x\in \mbox{\. X}^2_i$ define $\psi _i(x)=(-1-\phi _i^1(x),\dots ,
-1-\phi _i^m(x))$ where $\phi _i=(\phi _i^1,\dots ,\phi _i^m)$. Let 
$Y_1=\mbox{\. X}^1\cup \partial X\times [0,-1)$ and $Y_2=\partial
X\times (0,1]\cup \mbox{\. X}^2$ where for $l=1,2$, $\mbox{\. X}^l=
\bigcup _{i\in I}\mbox{\. X}^l_i$.
}
\end{nrmk}

{\it We will from now on assume that, all definable manifolds are
Hausdorff definable manifolds, hence affine.}
\end{subsection}

\begin{subsection}{The fundamental group}
\label{subsection fundamental group}

Let $\mathbf{X}$$=(X, (X_i,\phi _i)_{i\in I})$ and $\mathbf{Y}$
$=(Y, (Y_j,\psi _j)_{j\in J})$ be definably connected definable
manifolds. Given $A_1,\dots ,A_l\subseteq X'$ and $B_1,\dots ,B_l
\subseteq Y'$ definable subsets of $X$ and $Y$ respectively and a
continuous definable map $f:X'\into Y'$ such that $f(A_i)\subseteq
B_i$, we write $f:(X',A_1,\dots ,A_l)\into (Y',B_1,\dots ,B_l)$ (if
$l=1$, $A_1=\{x\}$ and $B_1=\{y\}$ we write $f:(X',x)\into (Y',y)$).

\begin{defn}\label{defn homotopy type}
{\em
Two continuous definable maps $f,g:(X',A_1,\dots ,A_l)\into
(Y',B_1,\dots ,B_l)$ are {\it definably homotopic} if there is a
continuous definable map $H:(X\times [0,1],A_1\times [0,1],\dots
,A_l\times [0,1])\into (Y,B_1,\dots ,B_l)$ such that $H(x,0)=f(x)$ and 
$H(x,1)=g(x)$ for all $x\in X'$. $H$ is called a {\it definable
homotopy}. This is an equivalence relation compatible with
composition in the set of all such continuous definable maps. We
denote by $[f]$ the equivalence class of $f$ and by $[(X',A_1,\dots
,A_l),(Y',B_1,\dots B_l)]$ the set of all such classes. $(X',A_1,\dots
,A_l)$ and $(Y',B_1,\dots ,B_l)$ are said to  have the same {\it
definable homotopy type} if there are continuous definable maps
$f:(X',A_1,\dots ,A_l)\into (Y',B_1,\dots ,B_l)$ and $g:(Y',B_1,\dots
,B_l)\into (X',A_1,\dots ,A_l)$ such that $[g\circ f]=[1_{X'}]$ and 
$[f\circ g]=[1_{Y'}]$. This is an equivalence relation and we denote
by $[(X',A_1,\dots ,A_l)]$ the equivalence class of $(X',A_1,\dots ,A_l)$.
}
\end{defn}

The following result is proved just like in the classical case, with
the word definable add where necessary (for details see chapter 12 \cite{f}):

\begin{fact}\label{fact fundamental group}
There is a covariant functor $\pi _1 $ from the category of pointed
definable manifolds into the category of groups. $\pi
_1(X,x)$$:=[([0,1],\{0,1\}), (X,x)]$ is a group called the definable 
fundamental group of ${\mathbf X}$ at $x$ with the product defined by 
$[\sigma ][\gamma ]:=[\sigma \cdot \gamma ]$ and given definable maps 
$f,g:(X,x)\into (Y,y)$, $\pi _1$$(f)$ (which is denoted $f_*$) is
defined by $f_*([\sigma ])=[f\circ \sigma ]$. Moreover, if $[f]=[g]$
then $f_*=g_*$;  if there is a definable path in $X$ from $x_0$ to
$x_1$ then $\pi _1(X,x_0)$ $\simeq \pi _1(X, x_1)$; and finally 
$\pi _1(X\times Y,(x,y))$ $\simeq \pi _1(X,x) \times \pi _1(Y,y)$.
\end{fact}

\begin{defn}
{\em
A definably connected definable manifold ${\mathbf X}$ is called 
{\it definably simply connected} if $\pi _1(X,x)=0$ for some
(equivalently for all) $x\in X$.
}
\end{defn}

We finish this subsection, with the following result from \cite{bo2}
(for a different proof see \cite{e2})
on how to compute the definable fundamental group of definable
manifolds. Note however that the main results of this paper on
definable groups were proved (in an old version of this paper) without 
using fact \ref{fact tietze} below. In fact, although we need to know
that the fundamental group of a definable group is a finitely
generated abelian group, this will be proved directly in fact
\ref{fact hurewicz} and lemma \ref{lem pi of an h space}. 

\begin{defn}\label{defn tietze}
{\em
Let $K$ be a simplicial complex, $v$ a vertice of $K$ and $T$ a
maximal tree in $K$ ($T$ contains all the vertices of $K$). An edge
path in $K$ is a sequence $\sigma =u_1, u_2, \dots , u_l$ of vertices
of $K$ such that for all $i=1, \dots , l-1$, $u_i, u_{i+1}$ are
vertices of an edge in $K$. If $\gamma =w_1, w_2,\dots ,w_k$ is
another edge path in $K$ and $u_l, w_1$ are vertices of an edge in $K$
then the concatenation $\sigma \cdot \gamma $ is also an edge path in
$K$. $\sigma $ is an edge loop in $K$ at $v$ if $v=u_1=u_l$. $E(K,v)$
is the group under the operation $[\sigma ][\gamma ]:=[\sigma \cdot
\gamma ]$ of classes $[\sigma ]$ of edge loops $\sigma $ in $K$ at $v$
under the equivalence relation: $v, u, a, b, c,\dots , w, v ~ v, u, a,
c, \dots , w, v$ iff $abc$ spans a $k$-simplex of $K$ where
$k=0,1,2$. $E(K,v)$ is isomorphic to the group $G(K,T)$ generated by
the $1$-simplexes of $K$, denoted $g_{a,b}$ for each edge $a,b$ and
with relations: $(1)$ $g_{a,b}=1$ if $a,b$ spans a simplex of $L$ and
$(2)$ $g_{a,b}g_{b,c}=g_{a,c}$ if $a,b,c$ spans a simplex of $K$.
}
\end{defn}

\begin{fact}\label{fact tietze}
(Tietze theorem, \cite{bo2}).
Suppose that ${\mathbf X}$ is a definably compact, definable
manifold. If $(\Phi ,K)$ is a definable triangulation of $X$ and $T$
is a maximal tree in $K$ then, $\pi _1(X,x)\simeq G(K,T)$. In
particular, $\pi _1(X,x)$ is invariant under taking elementary
extensions, elementary substructures of $\N$ (containing the
parameters over which ${\mathbf X}$ is defined) and under taking
expansions of $\N$ and reducts of $\N$ on which ${\mathbf X}$ is definable. 
\end{fact}

\end{subsection}
\end{section}

\begin{section}{Basic homological algebra}
\label{section basic homological algebra}

The category of chain complexes, denoted $\mathbf{Comp}$ is the
category whose objects are chain complexes $(E_*,\partial _*)$ - which
are sequences $(E_n)_{n\in \ZZ}$ of abelian groups with morphisms
$\partial _n:E_n\into E_{n-1}$ such that $\partial _{n-1}\circ
\partial _n=0$ for all $n\in \ZZ$ - and whose morphisms are chain maps
$f:E_*\into F_*$ between chain complexes $(E_*, \partial _*)$ and
$(F_*,\delta _*)$, i.e., $f$ is a sequence $(f_n)_{n\in \ZZ}$ of group
homomorphisms $f_n:E_n\into F_n$ such that $f_{n-1}\circ \partial
_n=\delta _n \circ f_n$ for all $n\in \ZZ$. A chain complex $E_*$ is
nonnegative if $E_n=0$ for $n<0$; 

The category of augmented chain complexes $\mathbf{\tilde{C}omp}$, is
the subcategory whose objects are those chain complexes $(E_*,\partial
_*)$ with $E_{-1}=\ZZ$ and $E_n=0$ for all $n<-1$, and whose morphisms
are all those chain maps $f$ between such chain complexes with
$f_{-1}=1_{\ZZ}$. An augmented chain complex $(\tilde{E}_*,
\tilde{\partial }_*)$ is an augmentation of the nonnegative chain
complex $(E_*,\partial _*)$ if $E_n=\tilde{E}_n$ for $n\geq 0$ and 
$\partial _n=\tilde{\partial }_n$ for $n>0$. 

\begin{defn}
{\em
The covariant functors $H_n:\mathbf{Comp}$$\into Ab$ for $n\in
\ZZ$ are defined by: $H_n(E_*)=ker\partial _n/im\partial _{n+1}$ and for
$f:E_*\into F_*$, $H_n(f):H_n(E_*)\into H_n(F_*)$ is given by
$cls_{E_*}z\into cls_{F_*}f_n(z)$ where $cls_{E_*}z$ denotes the
equivalence class of $z\in ker\partial _n$ in $H_n(E_*)$. $H_n(E_*)$
is the {\it $n$-th homology group of $E_*$} and $f_*$ is usually used to
denote the sequence $(H_n(f))_{n\in \ZZ}$.
}
\end{defn}

The homology groups $H_n(E_*)$ measure how far from being exact the
chain complex $(E_*,\partial _*)$ is. We say that $E_*$ is exact (or
acyclic) if $H_n(E_*)=0$ for all $n$. Under some conditions we may
have $H_n(E_*)\simeq H_n(F_*)$ for all $n$, for example if $F_*$ is an 
adequate subcomplex of $E_*$ - the theorem on removing cells, tells us
that under certain conditions, we can remove an $n$-cell and an
$n+1$-cell (i.e., elements from the set of generators of $E_n$ and
$E_{n+1}$ respectively) in order to obtain an adequate subcomplex
$F_*$ of $E_*$.  Also, if  a chain map $f:E_*\into F_*$ is a chain
equivalence i.e., there is a chain map $g:F_*\into E_*$ such that
$g\circ f\simeq 1_{E_*}$ and $f\circ g\simeq 1_{F_*}$, then
$H_n(f):H_n(E_*)\into H_n(F_*)$ is an isomorphism for every $n$. Here,
$f\simeq h$ means that the chain maps $f:E_*\into F_*$ and $h:
F_*\into E_*$ are chain homotopic i.e., there  is a sequence of maps 
$s_n:E_n\into F_{n+1}$ such that $f_n-h_n=\delta _{n+1}\circ s_n
+s_{n-1}\circ \partial _n$ for all $n$.

\begin{fact}\label{fact long exact sequence}
(Long exact sequence). Let $0{\rightarrow }E'_*\stackrel{i}{\rightarrow }
E_*\stackrel{p}{\rightarrow }E''_*{\rightarrow }0$ 
be a short exact sequence of chain
maps. There exists a sequence $d=(d_n:H_n(E"_*)\into H_{n-1}(E'_*))$
of (connecting) homomorphisms such that the following sequence is exact.
$$\cdots \into H_n(E'_*)\stackrel{i_*}{\rightarrow
}H_n(E_*)\stackrel{p_*}{\rightarrow }
H_n(E''_*)\stackrel{d_n}{\rightarrow }H_{n-1}(E'_*)\into \cdots $$
\end{fact}

\begin{fact}\label{fact naturality of the connecting homomorphim}
(Naturality of the connecting homomorphism). Assume that there is a 
commutative diagram of chain  complexes with exact rows:
\[
\begin{array}{clcr}
0 {\rightarrow }E'_*\stackrel{i}{\rightarrow}E_*\stackrel{p}{\rightarrow}E''_*
{\rightarrow } 0\\
\,\,\,\,\,\,\,\,\,\,\,\,\,\,\,\,
{\downarrow}^{f'}\,\,\,\,\,\,\,\, 
{\downarrow}^{f}\,\,\,\,\,\,\,\,\,
{\downarrow}^{f''}\,\,\,\,\,\,\,\,\,\,\,\,\\
0 {\rightarrow }F'_*\stackrel{i}{\rightarrow}F_*\stackrel{p}{\rightarrow}F''_*
{\rightarrow } 0\\
\end{array}
\]

Then there is a commutative diagram of abelian groups with exact rows:
\[
\begin{array}{clcr}
\cdots \into H_n(E'_*)\stackrel{i_*}{\rightarrow
}H_n(E_*)\stackrel{p_*}{\rightarrow }
H_n(E''_*)\stackrel{d_n}{\rightarrow }H_{n-1}(E'_*)\into \cdots\\
\,\,\,\,\,\,\,\,\,\,\,\, 
{\downarrow}^{f'_*}\,\,\,\,\,\,\,\,\,\,\,\,\,\,\,\,\,\,\,\,\,\,\, 
{\downarrow}^{f_*}\,\,\,\,\,\,\,\,\,\,\,\,\,\,\,\,\,\,\,\,\,\,\,
{\downarrow}^{f''_*}\,\,\,\,\,\,\,\,\,\,\,\,\,\,\,\,\,\,\,\,\,\,\,
{\downarrow}^{f'_*}\,\,\,\,\,\,\,\, \\
\cdots \into H_n(F'_*)\stackrel{i_*}{\rightarrow
}H_n(F_*)\stackrel{p_*}{\rightarrow }
H_n(F''_*)\stackrel{d'_n}{\rightarrow }H_{n-1}(F'_*)\into \cdots
\end{array}
\]
\end{fact}

\begin{defn}
{\em
Let $F:\mathbf{C}$$\into Ab$ be a functor, and let $\mathbf{M}$ be  a
set of models for $\mathbf{C}$ ( i.e., a subset of
$Obj\mathbf{C}$). We say that $F$ is free with base in $\mathcal{M}$
if: $(1)$ $FC$ is a free abelian group for every $C\in Obj\mathbf{C}$; 
$(2)$ there is an indexed family $(M_j){j\in J}$ of models in
$\mathbf{M}$ and an element $x=(x_j)\in \Pi _{j \in J}FM_j$ such that,
for every $C\in Obj\mathbf{C}$, the family $((F\sigma )(x_j))_{j\in
J,\sigma :M_j\into C}$ is a basis of $FC$. In this situation, we call
$x$ a basis of $F$ in $\mathbf{M}$.
}
\end{defn}

Let $E:\mathbf{C}$$\into \mathbf{Comp}$ be a functor. We say that $E$
is nonnegative if $EC$ is a nonnegative chain complex for all $C\in
Obj\mathbf{C}$; $E$ is acyclic on a set of models $\mathbf{M}$ for
$\mathbf{C}$ if $EM$ is an acyclic chain complex for every $M\in
\mathbf{M}$. A natural transformation between functors
$E,F:\mathbf{C}$$\into \mathbf{Comp}$ is called a  natural chain map.

\begin{fact}\label{fact acyclic models}
(Acyclic Models). Let $\mathbf{C}$ be a category with set of models 
$\mathbf{M}$ and let $E,F:\mathbf{C}$$\into \mathbf{Comp}$ be
nonnegative functors. Assume that $F_n$ is free with basis in
$\mathbf{M}$ for all $n>0$ and that $E$ is acyclic on
$\mathbf{M}$. Then: $(1)$ every natural transformation $\tau
_0:F_0\into E_0$ lifts to a natural chain map $\tau :F\into E$; $(2)$
any two liftings  $\tau , \tau ':F\into E$ are naturally chain
homotopic i.e., there are natural transformations $s_n:F_n\into
E_{n+1}$ such that $\tau _n-\tau '_n=\delta _{n+1}\circ
s_n+s_{n-1}\circ \partial _n$ for all $n$. Furthermore, we can choose $s_0=0$.
\end{fact}

An augmented natural chain map $\tau :E\into F$ is by definition a
natural transformation $\tau :E\into F$ between two functors
$E,F:\mathbf{C}$$\into $$\mathbf{\tilde{C}omp}$. Note that in this case,
for every $C\in Obj\mathbf{C}$, the map $(\tau
_{\mathbf{C}}$$)_{-1}:E_{-1}C=\ZZ$$\into F_{-1}C=\ZZ$ is $1_{\ZZ}$.

\begin{fact}\label{fact acyclic models in augmented}
(Acyclic Models in $\mathbf{\tilde{C}omp}$).  Let $\mathbf{C}$ be a
category with set of models $\mathbf{M}$ and let
$E,F:\mathbf{C}$$\into $$\mathbf{\tilde{C}omp}$ be functors. Then: $(1)$
if $F_n$ is free with basis in $\mathbf{M}$ for all $n>0$ and if $E$
is acyclic on $\mathbf{M}$ then there exist an augmented natural chain
map $\tau :F\into E$ and any two are naturally chain homotopic; $(2)$
if both  $E$ and $F_n$ are free with basis in $\mathbf{M}$ for all
$n>0$ and if both $E$ and $F$ are acyclic on $\mathbf{M}$ then every 
augmented natural chain map $\tau :F\into E$ is a natural chain equivalence.
\end{fact}

\end{section} 

\begin{section}{Homology}\label{section homology}

We denote by $DTOP(\N)$ (resp., $DCTOP(\N)$) the category of definable
manifolds (resp., the subcategory of definably compact manifolds) with 
continuous definable maps.

\begin{subsection}{Homology}
\label{subsection homology}

\begin{defn}
{\em
The {\it category of pairs} over $\N$, $DTOPP(\N)$ is the category whose
objects are the pairs $(X,A)$ where $X$ is a definable manifold and
$A$ a definable subset of $X$, and whose morphisms 
$f:(X,A)\into (Y,B)\in MorDTOPP(\N)$,  are the 
continuous definable maps $f:X\into Y$, such that $f(A)\subseteq B$.

The lattice of $(X,A)\in ObjDTOPP(\N)$ consists of the pairs
$(\emptyset ,\emptyset)$, $(A ,\emptyset)$, $(A,A)$, $(X,\emptyset)$, 
$(X,A)$ and $(X,X)$
all their identity maps, the inclusion maps and all their compositions.
$f:(X,A)\into (Y,B)\in MorDTOPP(\N)$,  induces a map of every
pair of the lattice of $(X,A)$ into the corresponding pair of the
lattice of $(Y,B)$.
}
\end{defn}

\begin{defn}
{\em
A subcategory $C$ of $DTOPP(\N)$ is called {\it admissible} if the
following conditions are satisfied.
\begin{enumerate}
\item
If $(X,A)\in ObjC$, then the lattice of $(X,A)$ is in $C$.
\item
If $f:(X,A)\into (Y,B)\in MorC$, then the induced  maps on every
pair of the lattice of $(X,A)$ into the corresponding pair of the
lattice of $(Y,B)$ is in $MorC$.
\item 
If $(X,A)\in ObjC$, then $(X\times [0,1],A\times [0,1])\in ObjC$,
and for $i\in \{0,1\}$, $g_i:(X,A)\into (X\times [0,1],A\times [0,1]):x\into (x,i)\in MorC$.
\item
Let $P=\{p\}$ be a one point set. Then $(P,\emptyset)\in ObjC$,
and for each $(X,\emptyset)\in ObjC$ and each point $x\in X$, then 
$(P,\emptyset)\into (X,\emptyset):p\into x \in MorC$. 
\item
If $(X,A), (Y,B)\in ObjC$, $X\subseteq Y$ and $A\subseteq B$ then 
the inclusion $(X,A)\into(Y,B)\in MorC$. 
\end{enumerate}
}
\end{defn}

The subcategory of $DTOPP(\N)$ of all definably
compact pairs over $\N$ is denoted by $DCTOPP(\N)$. For example $DTOPP(\N)$
and $DCTOPP(\N)$ are admissible subcategories of $DTOPP(\N)$.

\begin{defn}
{\em
Let $C$ be an admissible subcategory of $DTOPP(\N)$. A {\it homology}
$(H_*,d_*)=(H_n,d_n)_{n\geq 0}$ on $C$ is a sequence of covariant
functors 
$H_n:C\into Ab$ for $n\geq 0$
such that the following axioms hold.

{\it Homotopy Axiom}. If $f_0,f_1:(X,A)\into (Y,B)\in MorC$ and there
is a definable homotopy in $C$ from $f_0$ to $f_1$ (i.e., there is  
$F:(X\times [0,1],A\times [0,1])\into (Y,B)\in MorC$, with $F\circ
g_0=f_0$ 
and  $F\circ g_1=f_1$ where for 
$i\in \{0,1\}$,  $g_i:(X,A)\into (X\times [0,1],A\times [0,1]):x\into
(x,i)
\in MorC$), then  
\[ H_n(f_0)=H_n(f_1):H_n(X,A)\into H_n(Y,B) \] for all $n\geq 0$. 

{\it Exactness Axiom}. There is a sequence of natural transformations  
$d_n:H_n\into H_{n-1}\circ G$, where 
$G:C\into C$ is the functor $(X,A)\into (A,\emptyset)$, such that for
each 
pair $(X,A)\in C$ with 
inclusions $i:(A,\emptyset )\into (X,\emptyset )$ and $j:(X,\emptyset
)\into 
(X,A)$, the following sequence is 
exact. 
\[ \cdots \into H_n(A,\emptyset)\stackrel{H_n(i)}{\rightarrow} 
H_n(X,\emptyset)\stackrel{H_n(j)}{\rightarrow}H_n(X,A) 
\stackrel{d_n}{\rightarrow}H_{n-1}(A,\emptyset)\into \cdots \]

{\it Excision Axiom}. For every $(X,A)\in ObjC$ and every definable
open 
subset $U$ of $X$ such that 
$cl_X(U)\subseteq int_X(A)$ and $(X-U,A-U)\in ObjC$, the inclusion 
$(X-U,A-U)\into (X,A)$ induces 
isomorphisms 
\[ H_n(X-U,A-U) \into H_n(X,A) \] for all $n\geq 0$.

{\it Dimension Axiom}. If $X$ is a one point set, then
$H_n(X,\emptyset)=0$ for all $n>0$. (The group $H_0(X,\emptyset)$, is
called 
the {\it coefficient group}). 
}
\end{defn}

\textbf{Notation:}
We will write $X\in ObjC$ for $(X,\emptyset )\in ObjC$, $f:X\into Y\in
MorC$ for $f(X,\emptyset )\into (Y,\emptyset )\in MorC$ and $H_n(X)$
for $H_n(X,\emptyset )$.

\begin{defn}
{\em
Two homologies $(H_*,d_*)$ and $(H'_*,d'_*)$ on $C$ are {\it
isomorphic} if 
there is a sequence of natural
equivalences $\tau _n:H_n\into H'_n$ for all $n\geq 0$ such that 

\[
\begin{array}{clcr}
\,H_{n+1}(X,A)\,\,\,\,\stackrel{d_{n+1}}{\longrightarrow} \,\,\,\,\,\,\,\, H_n(A)  \\
{\downarrow}^{\tau _{n+1}}\,\,\,\,\,\,\,\,\,\,\,\,\,\,\,\,\,\,\,\,\,\,\,\,\,\,\,\,\,\,\,\,\,\,\,\,\,\, 
{\downarrow}^{\tau _{n}}\\
\,H'_{n+1}(X,A)\,\,\,\,\stackrel{d'_{n+1}}{\longrightarrow}\,\,\,\,\,\,\,\,H'_n(A)
\end{array}
\]
for all $(X,A)\in ObjC$ and for all $n\geq 0$.
}
\end{defn}

\begin{fact}\label{fact eilenberg steenrod}
(Eilenberg-Steenrod Theorem)
Homology functors on $DTOPP(\N)$ with isomorphic coefficient group are 
isomorphic. 
\end{fact}

\pf
The proof in \cite{h} can be adapted to show that the
Eilenberg-Steenrod axioms characterise homology on $DTOPP(\N)$.
\qed

\end{subsection}

\begin{subsection}{O-minimal simplicial homology}
\label{subsection o-minimal simplicial homology}

A. Woerheide \cite{Wo} defines {\it simplicial homology} on $DCTOPP(\N)$. Like 
in the case for simplicial homology over $\RR$, the main complication
is defining the induced homomorphisms between the homology groups. The
standard method, using repeated barycentric subdivision and the 
Lebesgue number property fails because $(N,0,1,<,+,\cdot )$ may be 
non-archemedian. But the use of the method of acyclic models and the
o-minimal triangulation theorem makes the construction possible.

We give here a brief description of Woerheide's construction. First we
need some basic terminology. 

\begin{defn}
{\em
Let $K$ be a closed simplicial complex in $N^n$. The {\it simplicial chain
complex} $(C_*(K),\partial _*)$ is defined by: $C_l(K):=A_l(K)/A'_l(K)$
where $A_l(K)$ is the free abelian group generated by the set of
$(l+1)$-tuples $(v_0,\dots ,v_l)$  and $A'_l(K)$ is the subgroup
generated by the
elements of the form $(v_0,\dots ,v_l)-sgn(\alpha )(v_{\alpha
(0)},\dots ,$$v_{\alpha (l)})$ where $\alpha $ is a permutation of
$\{0,\dots ,l\}$ and $v_0,\dots ,v_l$ spans an
$l$-simplex in $K$; and for $l>0$,
$\partial _l$ is induced by the homomorphism $(v_0,\dots ,v_l)\into
\sum _{i=0}^l(-1)^i(v_0,\dots ,\hat{v}_i,\dots ,v_l)$. We define
$\tilde{\partial }_0:C_0(K)\into \ZZ$ by $\sum \alpha _iv_i\into \sum
\alpha _i$, in this way we obtain an augmentation $\tilde{C}_*(K)$ of
$C_*(K)$. Note that each $\tilde{C}_l(K)$ is a free abelian group, and
given any total order on the vertices of $K$, the set of all classes
$<v_0,\dots ,v_l>\in \tilde{C}_l(K)$ such that $v_0<\cdots <v_l$
provides a basis for $\tilde{C}_l(K)$.
If $(K,K')$ is a closed simplicial pair, then $\tilde{C}_*(K')$ is a
subcomplex of $\tilde{C}_*(K)$ and we define the relative simplicial
chain complex $C_*(K,K'):=\tilde{C}_*(K)/\tilde{C}_*(K')$. Note that
$C_*(K)$$\simeq C_*(K,\emptyset)$. We define
$H_*(K,K'):=H_*(C_*(K,K'))$. 
}
\end{defn}

For a closed simplicial complex $K$, let $\mathbf{A}$$(K)$ denote the
category whose objects are closed simplicial subcomplexes of $K$ and
whose morphism are the inclusion maps. A map $f:|K|\into |L|$ between
the geometric realizations of closed simplicial complexes $K$ and $L$
is said to compatible if, for each simplex $t\in K$, there is a
simplex $s\in L$ such that $f(t)\subseteq s$. For each compatible map
$f:|K|\into |L|$, we define the functor
$\mathbf{A}$$(f):\mathbf{A}$$(L)\into \mathbf{A}$$(K)$ by
$\mathbf{A}$$(f)(L')=\{t\in K:f(t)\subseteq |L'|\}$. Let
$\tilde{C}_*^L:$$\mathbf{A}$$(L)\into \mathbf{\tilde{C}omp}$ be the
functor which sends a closed subcomplex $L'$ of $L$ to
$\tilde{C}_*(L')$ and sends an inclusion $L''\into L'$ to the inclusion 
$\tilde{C}_*(L'')\into \tilde{C}_*(L')$.  

An easy application of the
theorem on acyclic models for $\mathbf{\tilde{C}omp}$ (see \cite{Wo})
shows that: $(1)$ there is an augmented natural chain map between the
functors $\tilde{C}_*^K\circ \mathbf{A}$$(f)$ and $\tilde{C}_*^L$ and
any two are naturally chain homotopic; $(2)$ if $K$ is a closed
simplicial complex and $(\Psi ,L)$ is a quasi-stratified triangulation
of $|K|$ (i.e., after a linear change of coordinates its a stratified
triangulation) such that $\Psi ^{-1}:|L|\into |K|$ is compactible then
every augmented natural chain map $\tilde{C}_*^L\circ
\mathbf{A}$$(\Psi ^{-1})\into \tilde{C}_*^K$ is a natural chain
equivalence. The proof of $(2)$ also involves the method of acyclic
models to establish that if $|K|$ is convex then $\tilde{C}_*(L)$ is
an acyclic chain complex (see \cite{Wo}).

We are now recall Woerheide's definition of the o-minimal simplicial 
homology.

\begin{defn}\label{defn h of f}
{\em 
\cite{Wo}
Let $f:(K,K')\into (L,L')$. If $f$ is compatible then there is an
augmented natural chain map $\tau :\tilde{C}_*^K\circ
\mathbf{A}$$(f)\into \tilde{C}_*^L$, then chain map $\tau
_L:\tilde{C}_*(K)\into \tilde{C}_*(L)$ induces a chain map $\bar{\tau
}_L:C_*(K,K')\into C_*(L,L')$.
Define $H_*(f)$ to be $H_*(\bar{\tau }_L)$. If $f:(K,K')\into (L,L')$
is not necessarily compactible, by the definable traingulation
theorem, there is a quasi-stratified triangulation $(\Phi ,M)$ of $|K|$ 
such that $\Phi ^{-1}:|M|\into |K|$ and $f\circ \Phi ^{-1}:|M|\into
|L|$ are compactible, so that $H_*(\Phi ^{-1})$ and $H_*(f\circ \Phi
^{-1})$ are defined. On the other hand, by the above, every augmented
natural chain map $\tilde{C}_*^M\circ \mathbf{A}$$(\Phi ^{-1})\into 
\tilde{C}_*^K$ is a natural chain equivalence. So $H_*(\Phi ^{-1})$ is
an isomorphism. Define $H_*(f)$ to be the composition $H_*(f\circ \Phi
^{-1})\circ (H_*(\Phi ^{-1}))^{-1}$. 
}
\end{defn}

Its easy to see that $H_*$ is a well defined functor.

\begin{defn}\label{defn o-minimal simplicial homology}
{\em
\cite{Wo}
For $(X,A)\in ObjDCTOPP(\N)$, let $T(X,A)$ denote the set of all definable 
triangulations $(\Phi ,K)$ of $X$ which respects $A$, and for each 
$(\Phi ,K)\in T(X,A)$, $K'$ denotes the closed subcomplex of $K$ such
that $\Phi(A)=|K'|$. The {\it $n$-th homology group} $H_n(X,A)$ is
defined to be the subgroup of $\Pi _{(\Phi, K)\in T(X,A)}H_n(K,K')$ 
consisting of all elements $\alpha $ such that, for all $(\Phi
_1,K_1),(\Phi _2,K_2)\in T(X,A)$, we have $\alpha _{(\Phi _2,K_2)}=
H_n(\Phi _2\circ \Phi _1^{-1})(\alpha _{(\Phi _1,K_1)}$. 

Given $f:(X,A)\into (Y,B)\in MorDCTOPP(\N)$ we define {\it $n$-th induced 
homomorphism} $H_n(f):H_n(X,A)\into H_n(Y,B): \alpha \into \beta $
where, for all $(\Phi ,K)\in T(X,A)$ and $(\Psi ,L)\in T(Y,B)$, we
have $\beta _{(\Psi ,L)}=H_n(\Psi \circ f\circ \Phi ^{-1})(\alpha
_{(\Phi ,K)})$. 
}
\end{defn} 

The verification of the Eilenberg-Steenrod axioms is now easy and shows that:

\begin{fact}\label{fact simplicial homology}
\cite{Wo}
There is a homology $(H_*,d_*)$ on $DCTOPP(\N)$ such that if $(X,A)\in
ObjDCTOPP(\N)$, $(\Phi, K)$ is a 
definable triangulation
of $X$ which respects $A$ and $K'$ is the subcomplex of $K$ such that
$|K'|=\Phi (A)$, then we have isomorphisms
$\pi _n^{(\Phi ,K)}:H_n(X,A) \into H_n(K,K')$ for all $n\geq 0$ and
such that if 
$f:(X,A)\into (Y,B)\in MorDCTOPP(\N)$ and $(\Psi ,L)$ is a definable 
triangulation of $Y$ respecting $B$ and such that $|L'|=\Psi (B)$, then

\[
\begin{array}{clcr}
\,H_{n}(X,A)\,\,\,\,\,\,\,\stackrel{H_n(f)}{\longrightarrow} \,\,\,\,\,\,\,\,\,\,\, H_n(Y,B)  \\
\,\,\,\,{\downarrow}^{\pi _{n}^{(\Phi ,K)}}\,\,\,\,\,\,\,\,\,\,\,\,\,\,\,\,\,\,\,\,\,\,\,\,\,\,\,\,\,\,\,\,\,\,\,\,\,\,
{\downarrow}^{\pi _n^{(\Psi ,L)}}\\
\,\,\,\,H_{n}(K,K')\stackrel{H_n(\Psi \circ f\circ \Phi
^{-1})}{\longrightarrow} H_n(L,L')
\end{array}
\]
\end{fact}

\end{subsection}

\begin{subsection}{O-minimal singular homology}
\label{subsection o-minimal singular homology}

A. Woerheide \cite{Wo} also defines {\it singular homology} on $DTOPP(\N)$.
In this case the construction is essentially the same as for the
standard singular
homology, only with the word ``definable'' added here and there. But,
for the same reasons as above the standard proof of the excision axiom
fails, and the difficulty is avoided by the use of the o-minimal 
triangulation 
theorem and the results obtained while constructing the simplicial homology. 

\begin{defn}\label{defn singular complex}
{\em
The {\it standard $n$-simplex} over $N$, $\Delta ^n$ is the convex
hull of 
the standard basis vectors $e_0,
\dots ,e_{n}$ in $N^{n+1}$. Let the standard $(-1)$-simplex $\Delta
^{-1}$ be 
the empty set. Let $X\in ObjDTOP(\N)$ 
be a definable set. For $n\geq -1$, we define $\tilde{S}_n(X)$ to be
the free 
abelian group on the set of definable
continuous maps $\sigma :\Delta ^n\into X$ (the {\it (singular) 
$n$-simplexes in $X$})
and for $n<-1$, we define $\tilde{S}_n(X)=0$. Note that $\tilde{S}_{-1}
(X)=\ZZ$. The elements of $\tilde{S}_n(X)$ are called the {\it $n$-chains}.

For $n>0$ and $0\leq i\leq n$ let $\epsilon ^n_i:\Delta ^{n-1}\into
 \Delta ^n$
 be the continuous definable map
given by $\epsilon ^n_i(\sum _{i=0}^{n-1}a_je_j):=\sum _{j\neq
 i}^na_je_j$. 
Let $\epsilon ^0_0:\Delta ^{-1}
\into \Delta ^0$ be the unique map. We define the {\it boundary}
 homomorphism 
$\partial _n:\tilde{S}_n(X)
\into \tilde{S}_{n-1}(X)$ to be the trivial homomorphism for $n<0$ and
 for 
$n\geq 0$, $\partial _n$ is given on basis
elements by 
$$\partial _n(\sigma ):=\sum _{i=0}^n(-1)^i\sigma \circ \epsilon ^n_i.$$
One verifies that $\partial ^2=0$ and so
$(\tilde{S}_*(X),\partial )$ 
is a chain complex, the {\it
augmented singular chain complex}.

Given $(X,A)\in ObjDTOPP(\N)$, then the {\it relative singular chain complex} 
$(S_*(X,A),\partial )$ is the quotient
chain complex $(\tilde{S}_*(X)/\tilde{S}_*(A),
\bar{\partial})$.
 We define the {\it singular chain
complex} $S_*(X)$ to be $S_*(X,\emptyset )$.

For $f:(X,A)\into (Y,B)\in MorDTOPP(\N)$, we have an induced chain map 
$f_{\sharp}:S_*(X,A)\into S_*(Y,B)$ given on the basis
elements of $\tilde{S}_*(X)$ by $f_{\sharp}(\sigma )=f\circ
\sigma $. We get a sequence
of functors $H_n:DTOPP(\N)\into Ab$ by setting
$H_n(X,A):=H_n(S_*(X,A))$ and $H_n(f):=H_n(f_{\sharp})$.
We set $\tilde{H}_n(X):=H_n(\tilde{S}_*(X))$ and 
$H_n(X):=H_n(X,\emptyset )$.
}
\end{defn}

\begin{defn}\label{defn cone construction}
{\em
Let $X\subseteq N^m$ be a convex definable set and let $p\in X$. The 
{\it cone construction over $p$ in $X$} is a sequence of homomorphisms 
$\tilde{S}_*(X)\into \tilde{S}_{*+1}(X)$$:z\into p.z$ defined as
follows: For $n<-1$, $p.$ is defined as the trivial homomorphism and
for $n\geq -1$ and a  basis element $\sigma $, 
$p.\sigma (\sum _{i=0}^{n+1}t_ie_i)$$=p$ if $t_0=1$ or 
$t_0p+(1-t_0)\sigma (\sum
_{i=1}^{n+1}\frac{t_i}{1-t_0}e_i)$ if $t_0\neq 1.$

Given definable sets $X,Y$ and $Z$, let  $F(X,Y)$  denote the free
abelian group on the set of all definable continuous maps from $X$
into $Y$. Given $\alpha =\sum _is_i\alpha _i\in F(X,Y)$ and $\beta
=\sum _jt_j\beta _j\in F(Y,Z)$, we define the {\it sharp operator} by
$$\beta _{\sharp }\alpha :=\sum _{i,j}s_it_j(\beta _j\circ \alpha _i) 
\in F(X,Z).$$   
}
\end{defn}

Note that $\sharp $ is associative, and since $\tilde{S}_*(X)=
F(\Delta ^n,X)$, every $z\in \tilde{S}_*(X)$ yields a chain map
$z_{\sharp }:\tilde{S}_*(\Delta ^n)\into \tilde{S}_*(X)$.

\begin{defn}\label{defn barycentric subdivision}
{\em
Let $X$ be a definable set. The barycentric subdivision
$Sd_n:\tilde{S}_n(X)\into \tilde{S}_n$ is defined as follows: for
$n\leq -1$, $Sd_n:\tilde{S}_*(X)\into \tilde{S}_*(X)$ is the trivial 
homomorphism, $Sd_{-1}:\tilde{S}_{-1}(X)\into \tilde{S}_{-1}(X)$ is
the identity and for $n\geq 0$, $Sd_n(z):=z_{\sharp }(b_n.Sd_{n-1}
\partial 1_{\Delta ^n})$, where $b_n$ is the barycenter of $\Delta ^n$.
}
\end{defn}

Note that $Sd_0$ is the identity and $Sd_n$ is natural i.e., it
commutes with $f_{\sharp }:\tilde{S}_n(X)\into \tilde{S}_n(Y)$. Some
lemmas on the cone construction and the sharp operator show that $Sd
=(Sd_n)_{n\in \ZZ}$ is a chain map $\tilde{S}_*(X)\into \tilde{S}_*(X)$.

\begin{defn}\label{defn tau k}
{\em
Let $K$ be a closed simplicial complex in $N^m$. We define the chain
map $\tau _{K}:\tilde{C}_*(K)\into \tilde{S}_*(|K|)$ by setting $\tau
_{K}=1_{\ZZ}:\tilde{C}_{-1}(K)\into \tilde{S}_{-1}(|K|)$ and for $n\geq 0$, 
$\tau _{K}(<v_0,\dots ,v_n>):=\sigma $ where $\sigma :\Delta ^n\into
|K|$ is such that $\sigma (\sum _{i=0}^nt_ie_i)=\sum _{i=0}^nt_iv_i$.
}
\end{defn}

\begin{defn}\label{defn the subdivision operator}
{\em
Let $(\Phi ,K)$ be a definable triangulation of $\Delta ^n$
compatible with the standard simplicial complex of $E^n$. Let $X$ be
a definable set. We define the subdivision operator
$Sd_i^K:\tilde{S}_i(X)\into \tilde{S}_i(X)$ where $i\leq n$ by:
$Sd_{-1}^K=1_{\ZZ}$ and for $0\leq i\leq n$, 
$$Sd_i^K(z):=(Sdz)_{\sharp }(\gamma _{i}^{n})_{\sharp }(\Phi ^{-1})_{\sharp
}\tau _{K}F_n<e_{n-i},\dots ,e_n>$$
 where $F_n:\tilde{C}_*(E^n)\into \tilde{C}_*(K)$ is the chain map
induced the unique (up to natural chain homotopy) augmented natural
chain map $\tilde{C}_*^{E^n}\into \tilde{C}_*^K\circ \mathbf{A}$
$(\Phi ^{-1})$, and $\gamma ^n_i:\Delta ^n\into \Delta ^i$ is defined by 
$$\gamma ^n_i(\sum _{j=0}^na_je_j)=\sum _{j=0}^i(a_{n-i+j}+\frac{\sum
_{k=0}^{n-i-1}a_k}{i+1})e_j.$$
}
\end{defn}

Its clear from the definition that $Sd^K_i$ is natural. A simple
computation shows that $Sd^K:=(Sd^K_i)_{i\leq n}$ is a natural partial
chain map $\tilde{S}_*\into \tilde{S}_*$ of order $n$. The proof of
the theorem on acyclic models shows that there is a natural partial
chain homotopy between the subdivision operator $Sd^K$ and the
restriction $1_{*\leq n}$ of the identity natural chain map
$1_*:\tilde{S}_*\into \tilde{S}_*$. This latter fact together with the 
result that says that given a definable set $X$ with definable open
sets $U$ and $V$ such that $X=U\cup V$, if $z\in \tilde{S}_n(X)$ then
there is a definable triangulation $(\Phi ,K)$ of $\Delta ^n$
compatible with $E^n$ such that $Sd^K_n(z)\in
\tilde{S}_n(U)+\tilde{S}_n(V)$ implies the excision axiom for the
singular o-minimal homology, and therefore:

\begin{fact}
\cite{Wo}
The sequence $H=(H_n)_{n\in \ZZ}$ of functors defined above is a
homology for $DTOPP(\N)$, called singular homology. 
\end{fact}

\end{subsection}

\begin{subsection}{Some properties of homology}
\label{subsection some properties of homology}

Like in the classical case, we have the following results which are
consequence of the axioms for homology. 

\begin{nrmk}\label{nrmk deformation retract}
{\em
Consider the diagram in $DTOPP(\N)$:
\[
\begin{array}{clcr}  
A\,\,\,\,\,\,\stackrel{i}{\rightarrow } \,\,X\\
\,\,\,\,\,\,\downarrow ^{1_A}\,\,\swarrow ^{r}\,\,\uparrow ^{1_X}\\
A\,\,\,\,\,\stackrel{i}{\rightarrow }\,\,X.
\end{array}
\]
We have a {\it definable weak retract} if $[r\circ i]=[1_A]$ (i.e., the map 
$r\circ i$ is definably homotopic to the map $1_A$); a {\it definable 
retract} if $r\circ i=1_A$; a {\it definable weak deformation retract}
if $[r\circ i]=[1_A]$ and $[i\circ r]=[1_X]$; and a {\it definable
deformation retract} if $r\circ i=1_A$ and $[i\circ r]=[1_X].$

If $i:A\into X$ is a definable (weak) retract, then the exactness axiom 
implies for all $n$ that $H_n(X)\simeq H_n(A)\oplus H_n(X,A)$ (we have 
$ker r_*=H_n(X,A).$) If $i:A\into X$ is a definable weak deformation
retract then for all $n\geq 0$, $H_n(X,A)=0$. In particular, $H_n(X,X)=0$ for
all $n\geq 0$.
}
\end{nrmk}

\begin{thm}\label{thm fg}
Let $X\in ObjDTOP(\N)$, then $H_q(X)$ is finitely generated for
every $q\ge 0$, $H_q(X)=0$ for all $q>dim X=:m$, $H_m(X)$ is a
free abelian group (possibly trivial) and $H_0(X)=\ZZ$$^k$ where $k$
is the number of definable connected components of $X$. Moreover,
$H_*(X)$ is invariant under taking elementary
extensions or substructures of $\N$ that contain the parameters over
which $X$ is defined and under taking expansions and reducts of $\N$ on
which $X$ is definable.
\end{thm}

\pf
Note that by remark \ref{nrmk deformation retract} and since by
proposition 3.3 in chapter 8 of \cite{vdd} $X$ is a definable
deformation retract of some $Y\in DCTOP(\N)$, we may assume that $X\in
DCTOP(\N)$. Let $(\Phi ,K)$ be a definable triangulation of $X$. Then
$H_q(X)=H_q(K)$ which by theorem 7.14 \cite{ro} is finitely generated for
every $q\ge 0$, $H_q(X)=0$ for all $q>dim X=:m$, $H_m(X)$ is a
free abelian group (possibly trivial). The fact that $H_0(X)=\ZZ$$^k$ where 
$k$ is the number of definable connected components of $X$ is proved
in theorem 4.14 \cite{ro}. 
\qed

\begin{fact}\label{fact hurewicz}
(Hurewicz Theorem). Suppose that ${\mathbf X}$ is a definable
manifold. Then $\pi _1(X,x)/\pi _1$$(X,x)^{(1)}\simeq H_1(X)$. In
particular, $\pi _1(X,x)/\pi _1(X,x)^{(1)}$ is finitely generated with 
finitely many relations. 
\end{fact}

\pf
See theorem 4.29 \cite{ro} (on lemma 4.26 \cite{ro} instead of the
function $u:[0,1]\into S^1$ given by $u(t)=e^{2\pi it}$, use any
definable continuous $v:[0,1]\into S^1$ such that $v(0)=v(1)$ and such
that $v _{|[0,1)}$ is a bijection of $[0,1)$ and $S^1$).
\qed

Note that the excision axiom is equivalent to the following: for all 
$X\in DTOP(\N)$ and all $X_1,X_2$ subsets of $X$ such that 
$X=int_X(X_1)\cup int_X(X_2)$, the inclusion $(X_1,X_1\cap X_2)\into
(X,X_2)$ 
induces isomorphisms $$H_n(X_1,X_1\cap X_2)\into H_n(X,X_2)$$ for all 
$n\geq 0.$

\begin{fact}\label{fact mayer vietoris}
(Mayer-Vietoris)
Consider the commutative diagram of inclusions in $DTOPP(\N)$
\[
\begin{array}{clcr}
(X_1\cap X_2, Z)\stackrel{i_1}{\longrightarrow}(X_1, Z)\stackrel{p}{\longrightarrow}(X_1,X_1\cap X_2)\\
{\downarrow}^{i_2}\,\,\,\,\,\,\,\,\,\,\,\,\,\,\,\,\,\,\,\,\,\,\, {\downarrow}^{g}\,\,\,\,\,\,\,\,\,\,\,\,\,\,\,\,\,\,\,\,\,\,\,{\downarrow}^{h}\\
(X_2,Z)\,\,\,\,\stackrel{j}{\longrightarrow}(X,Z)\stackrel{q}{\longrightarrow }\,\,\,\,(X,X_2),
\end{array}
\]
where $X_1, X_2,Z\subseteq X$ with
$X=int_X(X_1)\cup int_X(X_2)$ and $Z\subseteq X_1\cap X_2$. Then there
is an exact sequence for all $n\in \NN$
\[
\begin{array}{clcr}
H_n(X_1\cap X_2,Z)\stackrel{(i_{1*},i_{2*})}{\rightarrow}
H_n(X_1,Z)\oplus H_n(X_2,Z)\stackrel{g_*-j_*}
{\rightarrow}H_n(X,Z)\stackrel{D}{\rightarrow}\\
\stackrel{D}{\rightarrow}H_{n-1}(X_1\cap X_2,Z) 
\end{array}
\]
with $D=l_*dh^{-1}_*q_*$, where $d$ is the connecting homomorphism of the
pair $(X_1,X_1\cap X_2)$ and $l:(X_1\cap X_2,\emptyset )\into (X_1\cap
X_2,Z)$ is the inclusion.
\end{fact}

\begin{fact}\label{fact exactness for triples}
(Exactness Axiom for triples). If we have inclusions
$$(A,\emptyset )\stackrel{c}{\rightarrow }(A,B)\stackrel{a}{\rightarrow }
(X,B)\stackrel{b}{\rightarrow }(X,A)$$ in $DTOPP(\N)$, then
there is an exact sequence for all $n\in \NN$
$$\cdots \into H_n(A,B)\stackrel{a_*}{\rightarrow
}H_n(X,B)\stackrel{b_*}{\rightarrow }
H_n(X,A)\stackrel{c_*\circ d_n}{\rightarrow }H_{n-1}(A,B)\into .
$$
\end{fact}

\begin{defn}\label{defn hom with coefficients}
{\em
Let $(X,A)\in DTOPP(\N)$ and let $G$ be an abelian group. The 
{\it singular chain complex $(S_*(X,A;G),\partial ^G)$
with coefficients in $G$} is defined as 
$S_n(X,A;G):=S_n(X,A)\otimes G$ and $\partial ^G_{n}:=\partial _n
\otimes 1$. 
For $f:(X,A)\into (Y,B)\in MorDTOPP(\N)$, we have an induced chain map 
$f_{\sharp}:S_*(X,A;G)\into S_*(Y,B;G)$ given on the basis
elements by $f_{\sharp}(\sigma \otimes g)=(f\circ
\sigma )\otimes g$. 
We have in this way a sequence
of functors $H_n(\,\, ;G):DTOPP(\N)\into Ab$ by setting $H_n(X,A;G):=
H_n(S_*(X,A;G))$ 
and $H_n(f):=H_n(f_{\sharp})$.
}
\end{defn}

Below, $Tor=Tor^{\ZZ}$ where for a ring $R$, $Tor ^R$ is the {\it
torsion 
functor}, and for $R$-modules $A$ and $B$, $Tor ^R(\,\,,B)$ measures
the 
failure of functor $\,\,\otimes _RB$ to be exact. Thus, if $0\into
C\into F
\into A\into 0$ is a free resolution of $A$, then $$ Tor^R(A,B)=ker
(C\otimes _RB\into F\otimes _RB).$$

We have the following result. For a proof see theorem 9.32 \cite{ro}.

\begin{fact}\label{fact coefficient for hom}
(Universal Coefficients Theorem for Homology).
For every $X\in ObjDTOP({\N})$ and for all 
$n\ge 0$, there are canonical exact sequences
$$0\rightarrow H_n(X)\otimes G\stackrel{\alpha}{\rightarrow}
H_n(X;G)\rightarrow Tor(H_{n-1}(X),G)\rightarrow 0 .$$
These sequences split, that is there are canonical isomorphisms 
$H_n(X;G)\simeq H_n(X)\otimes G\oplus Tor(H_{n-1}(X),G)$. In
particular, if $K_p$ is a field of characteristic $p$ (a prime or
zero) then we have an isomorphism $H_n(X;K_p)\simeq H_n(X)\otimes K_p$
iff $H_*(X)$ is $p$-torsion free.
\end{fact}

\begin{thm}\label{thm kunneth for hom}
(K\"unneth Formula for Homology). Let $R$ be a principal ideal
domain. For every $X,Y\in ObjDTOP(\N)$ and 
for all $n\ge 0$, there are canonical exact sequences of $R$-modules
$$0\rightarrow \sum _{i+j=n}H_i(X;R)\otimes _RH_j(Y;R)\stackrel{\alpha ''}
{\rightarrow}H_n(X\times Y;R)\rightarrow $$ $$\rightarrow \sum
_{p+q=n-1}Tor^R(H_p(X;R),H_q(Y;R))\rightarrow 0 .$$
These sequences split, that is there are canonical isomorphisms  
$$H_n(X\times Y;R)\simeq \sum _{i+j=n} H_i(X;R)\otimes _RH_j(Y;R)\oplus $$ 
$$\oplus \sum _{p+q=n-1}Tor^R(H_p(X;R),H_q(Y;R)).$$
In particular, if $H_*(X;R)$ or $H_*(Y;R)$ is a finitely generated
free $R$-module, then the homology (external) cross product $\alpha
''$ is an isomorphism of graded $R$-modules.
\end{thm}

\pf
We include here a proof for the case $R=\ZZ$. For the general case see
proposition 2.6 \cite{d}.
By the Eilenberg-Zilber theorem (see theorem 9.33 \cite{ro}) there is
a 
natural chain equivalence
$\zeta :S_*(X\times Y)\into S_*(X)\otimes S_*(Y)$
which 
induces isomorphisms
$H_n(X\times Y)\simeq H_n(S_*(X)\otimes S_*(Y))$ for all 
$n\geq 0$. Here,  
$S_*(X)\otimes S_*(Y)$ is the chain complex with 
$$(S_*(X)\otimes S_*(Y))_n=\sum_{i+j=n}S_i(X)\otimes S_j(Y)$$
and whose differentiation $D_n: (S_*(X)\otimes S_*(Y))_n\into $
$(S_*(X)\otimes S_*(Y))_{n-1}$ is defined on the generators by
$$D_n(a_i\otimes b_j)=\partial a_i\otimes b_j+(-1)^ia_i\otimes 
\partial b_j,\,\,\,\,\,\,\,\, i+j=n.$$

The K\"unneth theorem (see theorem 9.36 \cite{ro}), gives a split
exact 
sequence with middle term
$H_n( S_*(X)\otimes S_*(Y))$ and the Eilenberg-Zilber
theorem identifies this term with
$H_n(X\times Y)$. Here $\alpha ''$ is given by
$$\alpha ''((cls a_i)\otimes (clsb_j)):=cls(\zeta '(a_i\otimes
b_j))=:clsa_i\times clsb_j$$
where $\zeta '$ is the inverse of an Eilenberg-
Zilber chain equivalence $\zeta $.

By theorem 12.25 \cite{ro}, the Eilenberg-
Zilber chain equivalence $\zeta :S_*(X\times Y)\into
S_*(X)\otimes S_*(Y)$ can be defined as the {\it
Alexander-Whitney
map}
$$\zeta _n(\sigma ):=\sum _{i+j=n}\sigma '\lambda _i^n\otimes \sigma
''\mu _j^n,$$
where $\sigma :\Delta ^n\into X\times Y$ and $\sigma '=\pi '\sigma $,
$\sigma ''=\pi ''\sigma $ (and where $\pi '$, $\pi ''$ are the
projections of $X\times Y$ onto $X$, $Y$, respectively). And for
$0\leq i\leq n$, the (affine) maps 
$\lambda _i^n,\mu _i^n:\Delta ^i\into \Delta ^n$ are given by
$$\lambda _i^n(\sum _{j=0}^ia_je_j):=\sum _{j=0}^ia_je_j+\sum_{j=i+1}^n0e_j$$ 
and
$$\mu _i^n(\sum _{j=0}^ia_je_j):=\sum _{j=0}^{n-i}0e_j+
\sum _{j=0}^{n-(i+1)}a_{j}e_{j+(i+1)}.$$
\qed

For the proof of the next fact see proposition 2.6 \cite{d}.

\begin{fact}\label{fact homology external cross product}
The homology cross product satisfies the following properties:
(1) $(f\times g)_*(\alpha \times \beta )=(f_*\alpha )\times (g_*\beta )
$ (naturality);
(2) $t_*(\alpha \times \beta )=(-1)^{deg\alpha deg\beta }\beta \times
\alpha $ (skew-commutativity) where $t:X\times Y\into Y\times X$
commutes factors;
(3) $(\alpha \times \beta )\times \gamma =\alpha \times (\beta \times \gamma
)$ (associativity);
(4) $1\times \alpha =\alpha \times 1=\alpha $ (unit element).
\end{fact}

\begin{nrmk}
{\em
Note that all the results presented in this subsection can be
generalised to relative homology. More precisely, we have a relative
version of the universal coefficient
theorem, the K\"unneth formula with corresponding cross product.
For details see \cite{d}.  
}
\end{nrmk}
\end{subsection}
\end{section}

\begin{section}{Cohomology}\label{section cohomology}
\begin{subsection}{Cohomology}
\label{subsection cohomology}

\begin{defn}
{\em
Let $G$ be an abelian group. Suppose that $C$ is an admissible
subcategory of $DTOPP(\N)$. A {\it cohomology $(H^*,d)=(H^n,d^n)_{n\geq
0}$ on $C$
with coefficients in $G$}
is a sequence of contravariant functors $H^n(\,\, ;G):C\into Ab$ for $n\geq 0$
satisfying the axioms dual to those satisfied by a homology functor. 
We denote by $H^n(X,A;G)$ the image of $(X,A)\in ObjC$ under $H^n$
and by $H^n(X;G)$ the image of $(X,\emptyset )\in ObjC$ under
$H^n(\,\, ;G)$.
}
\end{defn}

A. Woerheide constructions together with classical arguments (see
chapter 12 of \cite{ro}) show that there exist simplicial cohomology
functor for the category $DCTOPP(\N)$ and singular
cohomology functors for the category $DTOPP(\N)$. 

\begin{defn}\label{defn cochains}
{\em
For a fixed abelian group $G$, recall that $Hom(\,\, ,G):Ab\into Ab$ is
a contravariant functor.
For $(X,A)\in ObjDTOPP(\N)$, we have the {\it augmented singular cochain
complex $(\tilde{S}^*(X;G),\delta )$} with coefficients in $G$,
where $\tilde{S}^n(X;G):=Hom(\tilde{S}_n(X),G)$ and $\delta
^n:=Hom(\partial _{n+1},G)$. We also have the {\it relative singular
cochain complex $(S^*(X,A;G),\delta )$} with
coefficients in $G$ given by
$S^n(X,A;G):=Hom(S_n(X,A),G)$ 
and $\delta ^n:=Hom(\partial _{n+1},G)$.
Given 
$f:(X,A)\into (Y,B)\in MorDTOPP(\N)$ we have a chain homomorphism 
$f^{\sharp}:S^*(Y,B;G)\into S^*(X,A;G)$ given
by $f^{\sharp}:=Hom(f_{\sharp},G)$.

The {\it singular cohomology functors $H^n(\,\,;G):DTOPP(\N)\into Ab$}
with coefficients in $G$ are given by $H^n(X,A;G):=H^n(S^n(X,A;G))$
and $H^n(f):=H^n(f^{\sharp})$. We set
$\tilde{H}^n(X;G):=H^n(\tilde{S}^*(X;G))$ and
$H^n(X;G):=H^n(X,\emptyset ;G)$. As usual, $H^n(X):=H^n(X,\ZZ)$.

}
\end{defn}

We now list some properties of the (singular) cohomology groups. There
are analogue results for homology groups, and the proofs are obtained
by taking their ``dual''. For fact \ref{fact dual universal coh} see
theorem 12.11 \cite{ro} and for fact \ref{fact universal coefficient
coh} see theorem 12.15 \cite{ro}.

Below, $Ext(\,\,,\,\,)$ is the {\it extension functor} defined as
follows. 
Let $A$ and $B$ be abelian groups and consider a free resolution
$0\into R
\into F\into A\into 0$ of $A$ (i.e., $F$ is a free abelian group). Then 
$$Ext(A,B)=coker(Hom(F,B)\into Hom(R,B)),$$
so that $Ext(\,\,,B)$ measures the failure of $Hom(\,\,,B)$ to be exact.

\begin{fact}\label{fact dual universal coh}
(Dual Universal Coefficients).
For every $X\in ObjDTOP({\N})$ and for all 
$n\ge 0$, there are canonical exact sequences
$$0\rightarrow Ext(H_{n-1}(X),G)\rightarrow H^n(X;G)\stackrel{\beta}
{\rightarrow} Hom(H_n(X),G)\rightarrow 0 .$$
These sequences split, that is there are canonical isomorphisms 
$H^n(X;G)\simeq Hom(H_n(X),G) \oplus Ext(H_{n-1}(X),G)$. In
particular, if $K_p$ is a field of characteristic $p$, then for all
$n\ge 0$, $H^n(X;K_p)\simeq Hom(H_n(X),K_p)$ iff $H_*(X)$ is
$p$-torsion free.
\end{fact}

\begin{fact}\label{fact universal coefficient coh}
(Universal Coefficients Theorem for Cohomology).
For every $X\in ObjDTOP({\N})$ and for all 
$n\ge 0$, there are canonical exact sequences
$$0\rightarrow H^n(X)\otimes G\stackrel{\alpha}{\rightarrow}
H^n(X;G)\rightarrow Tor(H^{n+1}(X),G)\rightarrow 0 .$$
These sequences split, that is there are canonical isomorphisms 
$H^n(X;G)\simeq H^n(X)\otimes G\oplus Tor(H^{n+1}(X),G)$. In
particular, if $K_p$ is a field of characteristic $p$, then for all
$n\ge 0$, $H^n(X;K_p)\simeq H^n(X)\otimes K_p$ iff $H_*(X)$ is
$p$-torsion free (equivalently iff $H^*(X)$ is $p$-torsion free).
\end{fact}

Recall that a ring $R$ (resp., an algebra over a ring $S$) is a graded
ring (resp., graded $S$-algebra) if there are additive subgroups
(resp., $S$-submodules) $R^n$, $n\ge 0$ such that $R=\sum _{n\ge
0}R^n$ (direct sum of additive subgroups (resp., $S$-submodules)) with
a homomorphism $R\otimes R\into R$ (called multiplication)
 of graded additive subgroup
(resp., $S$-modules) such that $R^nR^m\subseteq R^{n+m}$. A unit is an
element $1\in R^0$ such that $1r=r1=r$ for all $r\in R$. The
multiplication is associative if $a(bc)=(ab)c$ for all $a,b,c\in R$. 
If $x\in R^n$ we say that $x$ is homogeneous
of degree $n$, we denote this by $deg(x)=n$. If $x=x_1+\cdots +x_k$,
with $x_i\in R^{n_i}$ we define the degree of $x$ to be $deg(x):=\sum
_{i=1}^kdeg(x_i)$. The multiplication is skew-commutative if
$ab=(-1)^{degadegb}ba$ for all $a,b\in R$.

\begin{defn}\label{defn kronecker product}
{\em
The product of $c\in S_n(X,A;G)$ and $\sigma \in S^n(X,A;G)$
defined by
$(\sigma ,c):=\sigma (c)$ satisfies $(\delta \sigma ,c)=(\sigma
,\partial 
c)$ and hence induces the {\it Kronecker product}
$$ (\,\, ,\,\,):H^n(X,A;G)\otimes H_n(X,A;G)\into G.$$
}
\end{defn}

\begin{thm}\label{thm cup product}
Let $R$ be a commutative ring. For every $X\in ObjDTOP({\N})$,
$H^*(X;R)=\sum _{n\ge 0}H^n(X;R)$ is a graded $R$-algebra with a
canonical multiplication (called cup product) $\cup :H^*(X;R)\otimes
H^*(X;R)\into H^*(X;R)$ satisfying:
$$ \theta \cup \psi =(-1)^{deg(\theta )deg(\psi )} \psi \cup \theta .$$ 
Moreover, $H^*(\,\,;R)$ is a functor from $DTOP(\N)$ into the category of
graded skew-commutative (associative) $R$-algebra with
unit element.
\end{thm} 

\pf
The cup product $\cup :S^n(X,R)\times S^m(X,R)\into S^{n+m}(X,R)$ is
defined 
by
$$(c ,\phi \cup \theta ):=(c\lambda _n^{n+m},\phi )(c\mu _m^{n+m},\theta ).$$
Here and through $(c,\phi ):=\phi (c)$.
Lemma 12.19 \cite{ro} shows that $S^*(X,R)$ is a graded
$R$-algebra under cup product, corollary 12.21 
\cite{ro} shows that $S^*(\,\, ,R)$ is a contravariant functor
from $DTOP(\N)$ into the category of graded 
$R$-algebras. Lemma 12.22 \cite{ro} and theorem 12.23 \cite{ro} shows
that $\cup $ can be defined on
$H^n(X,R)\times H^m(X,R)\into H^{n+m}(X,R)$ by $cls\phi \cup cls\theta 
:=cls(\phi \cup \theta)$ and 
$H^*(\,\, ,R)$ is a contravariant functor from $DTOP(\N)$ into the
category of graded $R$-algebras.

Theorem 12.26 \cite{ro} shows that the cup product is the composite
$\Delta ^{\sharp}\zeta ^{\sharp}\pi $, where $\Delta :X\into X\times X$
is the diagonal and 
$$\pi:S^*(X,R)\otimes S^*(Y,R)\into
Hom(S_*(X)\otimes
S_*(Y),R)$$ is defined as follows: If $\phi \in S^n(X,R)$ and
$\theta \in S^m(Y,R)$, then there is a function 
$S^{n}(X,R)\otimes S^{m }(Y,R)\into
Hom(S_{n}(X)\otimes S_{m}(Y),R)$ by $(\phi \otimes\theta )\into
(\,\, ,\phi\otimes \theta  )$, where $(\sigma \otimes \tau ,\phi \otimes
\theta)=(\sigma ,\phi)(\tau ,\theta)$.
Since this  function is bilinear, it extends to the homomorphism $\pi $. 
This is used in theorem 12.29 \cite{ro}
to show the skew-commutativity.
\qed

\begin{nrmk}\label{nrmk relative cup products}
{\em
The same proof also shows that we have the following {\it relative cup 
products}:
$$ \cup :H^*(X,A;R)\otimes _RH^*(X;R)\into H^*(X,A;R);$$
$$ \cup :H^*(X;R)\otimes _RH^*(X,A;R)\into H^*(X,A;R);$$
$$ \cup :H^*(X,A;R)\otimes _RH^*(X,A;R)\into H^*(X,A;R).$$
satisfying the above properties. They are interrelated via the
homomorphism 
$j^*:H(X,A;R)\into H^*(X;R).$ The relation with $d^*:H^*(A;R)\into 
H^{*+1}(X,A;R)$ is given by 
$d^*(\alpha \cup i^*\beta )=d^*\alpha \cup \beta ,$ for $\alpha \in 
H^p(A;R)$, $\beta \in H^q(X;R)$, where $i^*:H^*(X;R)\into H^*(A;R).$
}
\end{nrmk}

\begin{thm}\label{thm kunneth for cohomology}
(K\"unneth Formula for Cohomology).
Let $R$ be any principal ideal domain.
For every $X,Y\in ObjDTOP({\N})$ and for all 
$n\ge 0$, there are canonical exact sequences of skew-commutative
(associative) $R$-algebras with unit
$$0\rightarrow \sum _{i+j=n}H^i(X;R)\otimes _RH^j(Y;R)\stackrel{\alpha '}
{\rightarrow}H^n(X\times Y;R)\rightarrow $$
$$\rightarrow \sum
_{p+q=n+1}Tor^R(H^p(X;R),H^q(Y;R))\rightarrow 0 .$$
These sequences split, that is there are canonical isomorphisms 
$$H^n(X\times Y;R)\simeq \sum _{i+j=n} H^i(X;R)\otimes _RH^j(Y;R)\oplus $$
$$\oplus \sum _{p+q=n+1}Tor^R(H^p(X;R),H^q(Y;R)).$$
In particular, if $H^*(X;R)$ or $H^*(Y;R)$ is a finitely generated
free $R$-module, then the cohomology (external) cross product $\alpha
'$ is an isomorphism of graded skew-commutative (associative) $R$-algebra with
unity element.
\end{thm}

\pf
See theorem 12.16 \cite{ro} and theorem 12.31 \cite{ro}. The cross
product homomorphism $\alpha '$ is determined by the {\it (external)
cross product homomorphism}
$$\zeta ^{\sharp}\pi :S^*(X,R)\otimes S^*(Y,R)\into
S^*(X\times Y,R)$$ (it is common to write $\phi \times \theta$
for $\zeta ^{\sharp}\pi (\phi \otimes\theta )$ and also
$\phi \times \theta$ for $\alpha '(\phi \otimes\theta )$).
\qed

\begin{nrmk}\label{nrmk properties of cohomology cross product}
{\em The cohomology cross product satisfies the dual properties of
those satisfied by the homology cross product. Moreover, it also
satisfies the following : $(\alpha \times \beta ,\sigma \times \tau )=
(-1)^{deg\tau deg\alpha }(\alpha ,\sigma )(\beta ,\tau )$ {\it (duality)}.
For details see chapter VII section 7 in \cite{d}.
}
\end{nrmk}

\begin{nrmk}\label{nrmk relative kunneth formula}
{\em
Note that, the cohomology cross product is related to the cup product
by 
$\alpha \times \beta =p_X^*\alpha \cup p_Y^*\beta $ where $p_X:X\times 
Y\into X$ and $p_Y:X\times Y\into Y$. This can be use to prove
relative versions of the K\"unneth formula for cohomology. 
}
\end{nrmk}

The following two results are proved just like in the classical case.

\begin{thm}\label{thm cap product}
Let $R$ be a commutative ring. For all $X\in ObjDTOP(\N)$, $H_*(X;R)$
is graded (right) $H^*(X;R)$-module under the bilinear pairing 
$$\cap :H_{p+q}(X;R)\times H^p(X;R)\into H_{q}(X;R)$$
called the cap product. Furthermore, the cap product satisfies the
following properties: (1) $(\alpha \cap \beta ,\gamma )=(\alpha ,\beta
\cup \gamma )$ (duality); (2) for $f:X\into Y \in MorDTOP(\N)$, we have
$f_*(\alpha \cap f^*\beta )=f_*\alpha \cap \beta $ (naturality). 
\end{thm}

\pf
The cap product is induced by $\cap :S_{p+q}(X,R)\times S^p(X,R)\into
S_q(X,R)$
give by 
$$(\sigma \otimes r)\cap c:=\sigma \lambda ^{p+q}_q\otimes r(\sigma
\mu ^{p+q}_q,c).$$
For details see theorem 66.2 \cite{m} and also corollary 24.22
\cite{g}. 
\qed

\begin{thm}\label{thm slant product}
There is a bilinear pairing $$/:H^{p+q}(X\times Y ,R)\times H_p(X,R)\into
H^q(Y,R)$$ called the slant product which satisfies: (1)
$(\beta ,\gamma /\alpha )=(\alpha \times \beta ,\gamma )$ (duality); (2)
$1/\alpha =(\alpha ,1)1$ (units); (3) for  $f:X\into X' \in MorDTOP(\N)$
 $f:Y\into Y' \in MorDTOP(\N)$ we have $(f\times g)^*(\gamma )/\alpha
=g^*(\gamma /f_*\alpha )$ (naturality) and (4)
$[(\zeta \times \beta )\cup \gamma ]/\alpha =(-1)^{q(p+q+r-s)}\beta
\cup [\gamma /\alpha \cap \zeta ]$.
\end{thm}

\pf
For details see (29.19), (29.20), (29.21), (29.22) and (29.23) in \cite{g}.
\qed

\begin{nrmk}
{\em
Note that as for the homology case,
all the results presented in this subsection can be
generalised to relative cohomology. More precisely, we have a relative
version of the Mayer-Vietoris sequence, the universal coefficient
theorem, the K\"unneth formula with corresponding cross product. We
also have cup products, cap products and slant products for the
relative case.
For details see \cite{d} or \cite{g}.
In subsequent sections we will need the following versions of cap products:
 $$\cap :H_{p+q}(X,A;R)\times H^p(X,A;R)\into H_q(X;R)$$ and
$$\cap :H_{p+q}(X,A;R)\times H^p(X;R)\into H_q(X,A;R).$$ 
}
\end{nrmk}
\end{subsection}
\end{section}

\begin{section}{The Euler characteristic}
\label{section the euler characteristic}
\begin{subsection}{Hopf algebras}
\label{subsection hopf algebras}

A graded skew-commutative associative $R$-algebra $H=\sum _{k\ge 0}H^k$
with unity element is called a {\it quasi Hopf $R$-algebra}
if each $H^k$ is a finite dimensional $R$-module and there is a degree
preserving $R$-algebra homomorphism $\mu :H\into H\otimes _RH$ called
{\it comultiplication}, for which $\mu (H^k)\subseteq \sum
_{i+j=k}H^i\otimes _RH^j$. A quasi Hopf $R$-algebra $H=\sum _{k\ge 0}H^k$
is {\it connected} if $H^0$ is an $R$-algebra of dimension one with
generator $e$ and the map $\epsilon :H\into R$, defined by $\epsilon
(e)=1$ and $\epsilon (h)=0$ for all $h\in H^k$ with $k\geq 1$ is a {\it
co-unit} i.e., for all $h\in H$,
$(\epsilon \otimes _R1)\mu (h)=h\otimes _R1$ and 
$(1\otimes _R\epsilon )\mu (h)=1\otimes _Rh$. A quasi-Hopf algebra $H$
is 
called an {\it Hopf algebra} if $\mu $ is associative i.e., 
$(\mu \otimes _R1)\mu =1\otimes _R\mu)\mu .$ $\mu $ is called
commutative if 
$T\circ \mu =\mu $ where $T(x\otimes _Ry)=(-1)^{deg(x)deg(y)}y\otimes _Ry.$ 

An example of a connected Hopf $R$-algebra, 
is the {\it free, skew-commutative graded Hopf $R$-algebra}
$$R[x_1,\dots ,x_k,\dots]\otimes _R\bigwedge [y_1,\dots ,y_l,\dots
]_{R},$$
where the $x_i$'s are of even degrees and the $y_j$'s are of odd
degrees, and we have the relations: $$y_j^2=-y_j^2=0,\,\,
y_iy_j=-y_jy_i,\,\, y_jx_i=x_iy_j,\,\, x_ix_j=x_jx_i$$ and the number of
$x_i$ and $y_j$ of each degree is finite. In view of the freeness of
this algebra, comultiplication is determined by its values on the
generators $x_i$ and $y_j$:
$$\mu (x_i)=x_i\otimes _R 1+1\otimes _Rx_i+\sum
_{deg(u_p)+deg(u_q)=deg(x_i)}u_p\otimes _Ru_q $$ and
 
$$\mu (y_j)=y_j\otimes _R 1+1\otimes _Ry_j+\sum
_{deg(u_p)+deg(u_q)=deg(y_j)}u_p\otimes _Ru_q .$$

\begin{fact}\label{fact hopf algebra over a field}
(Theorem 7.6 \cite{dfn}) If $R$ is a field of characteristic
zero, then a connected Hopf $R$-algebra is a free, skew-commutative
graded Hopf $R$-algebra. 
\end{fact}

We also have the following result:

\begin{fact}\label{fact hopf algebra}
(\cite{mt} chapter VII, corollary 1.4). Let $H$ be a quasi Hopf
algebra over 
a perfect field $K_p$ of characteristic $p$. Then we have the
following ring 
isomorphisms:\\
(1) For $p=0$; $H\simeq (\bigotimes _{\alpha }\bigwedge [x_{\alpha
}]_{K_0})
\otimes (\bigotimes _{\beta }K_0[x_{\beta }])$, where $degx_{\alpha }$
is odd 
and $degx_{\beta }$ is even.\\
(2) For $p=2$; $H\simeq (\bigotimes _{\alpha }K_2[x_{\alpha }]/
(x_{\alpha }^{h_{\alpha }}))\otimes (\bigotimes _{\beta
}K_2[x_{\beta }])$,
 where $h_{\alpha }$ is a power of $2$.\\
(3) For $p\neq 0,2$; $H\simeq (\bigotimes _{\alpha }\bigwedge
[x_{\alpha }]_{K_p})
\otimes (\bigwedge _{\beta }K_p[x_{\beta }])\otimes (\bigotimes
_{\gamma }
K_p[x_{\gamma }]/(x_{\gamma }^{h_{\gamma }}),$ where $degx_{\alpha }$
is 
odd, $degx_{\beta }$ and $degx_{\gamma }$ are even, and $h_{\gamma }$ is a 
power of $p$. 

Here, if $dimH<\infty $, then there is no term of $K_p[x_{\beta }]$. 
\end{fact}

\begin{defn}\label{defn more on hopf spaces}
{\em
Let $(X,e,m)$ be a {\it definable $H$-manifold} i.e., $X$ is a
definable manifold with a continuous 
{\it definable $H$-multiplication} $m:(X\times X,(e,e))\into (X,e)$
and a {\it $H$-unit} $e\in X$ such that $[m\circ i_1]=[1_X]=[m\circ i_2]$
where
$i_1,i_2:X\into X\times X$ are the inclusions $i_1(x):=(x,e)$ and
$i_2(x):=(e,x)$. We say that $m$ is {\it definably
$H$-commutative} if $[m]=[m\circ t]$ (where $t:X\times X\into X\times X$
commutes factors); $m$ is {\it definably $H$-associative} if 
$[m\circ (m\times 1_X)]=[m\circ (1_X\times m)]$. A definable $H$-manifold
$(X,e,m)$ is a {\it definable $H$-group} if $m$ is $H$-associative and
has an {\it $H$-inverse} i.e., a definable continuous map 
$\iota :X\into X$ such
that $[m\circ (\iota \times 1_X)\circ \Delta _X]=[e]=$
$[m\circ (1_X\times \iota )\circ \Delta _X]$ where 
$\Delta _X:X\into X\times X$ is the diagonal map.

A definable continuous map $f:(X,e)\into (X',e')$ between definable 
$H$-manifolds (resp., definable $H$-groups) 
$(X,e,m)$ and $(X',e',m')$ is called a definable $H$-map (resp., 
$H$-homomorphism) if $[h\circ m]=[m'\circ (h\times h)]$ (resp., also 
$[h\circ \iota ]=[\iota '\circ h]$).
}
\end{defn}

Fact \ref{fact hopf algebra over a field} together with the K\"unneth
formula gives the following result.

\begin{thm}\label{thm hopf algebra of a group}
Let $X\in ObjDTOP({\N})$ and let $R$ be either a field or
$\ZZ$ in which case we assume that
$H_n(X)$ are all free abelian groups. If $X$ is a definably connected
definable $H$-manifold
(in particular, if $X$ is a definable group) then $H^*(X;R)$ is a 
connected quasi Hopf $R$-algebra. Moreover, if $R$ is a field of characteristic
zero, then $H^*(X;R)$ is isomorphic to the finitely generated free
exterior $R$-algebra $\bigwedge [y_1,\dots ,y_r]_{R}$ and therefore
comultiplication is given by $\mu (x)=x\otimes _R1+1\otimes_Rx$.
If $K_p$ is a field of characteristic $p$ then graded $K_p$-algebra
$H^*(X;K_p)$ is given in fact \ref{fact hopf algebra}. 
\end{thm}

\pf
For the fact that $H^*(X;R)$ is a connected quasi Hopf $R$-algebra see 
the proof of theorem 12.42 \cite{ro} and the remark that follows it.
Note that comultiplication $\mu $ is given by $\mu :=(\alpha
')^{-1}\circ m^*$.

By fact \ref{fact hopf algebra over a field}, if $R$ is a field of 
characteristic
zero then $H^*(X;R)$ is isomorphic to a  free, skew-commutative 
graded Hopf $R$-algebra. 
If there were any free
generators $x_i$ of positive even degree, then there would be elements
of arbitrary high degree and this is impossible 
since for $n>dimX$ we have $H^n(X;R)=0$. Therefore,
$H^*(X;R)$ is isomorphic to the finitely generated free
exterior $R$-algebra $\bigwedge [y_1,\dots ,y_r]_{R}$ and 
comultiplication is given by $\mu (x)=x\otimes _R1+1\otimes_Rx$.
\qed

Let $X$ be as above. We call $r$ the {\it rank of $X$} and denote this
by $r=rank(X)$. Also for each $i=1,\dots ,r$, $g_i:=degy_i$ and we
write $y_{i_j}$ instead of $y_i$ if $g_i=j$.

\medskip
We will now prove the dual of theorem \ref{thm hopf algebra of a
group}.

\begin{thm}\label{thm pontrjagin product}
Let $R$ be a
commutative ring. Then $H_*(\,\,;R)$ is a covariant functor from the
category of definable (resp., definable
$H$-commutative, $H$-associative) definable $H$-manifold into the
category of graded (resp., skew-
commutative, associative) $R$-algebras with unit element.

If $(X,e,m)$ be a definably connected definable
$H$-commutative and $H$-associative definable $H$-manifold of finite
type, and $R$ is either a field or  $\ZZ$ in which case we assume that
$H_n(X)$ are all free abelian groups. Then $H_*(X;R)$ is a 
connected quasi Hopf $R$-algebra. Moreover, if $R$ is a field of characteristic
zero, then $H_*(X;R)$ is isomorphic to the finitely generated free
exterior $R$-algebra $\bigwedge [y_1,\dots ,y_r]_{R}$ and therefore
comultiplication is given by $\mu (x)=x\otimes _R1+1\otimes_Rx$.
\end{thm}

\pf
Let $(X,e,m)$ be a definable $H$-manifold. The multiplication in 
$H_*(X;R)$ is called the {\it Pontrijagin product} and is defined by 
$\alpha ''\circ m_*$.
For details see chapter VII section 3 in \cite{d}. Co-multiplication
in the other case is given by $(\alpha '')^{-1}\circ \Delta _{X*}$ (for
details see chapter VII section 10 in \cite{d}).
\qed

\end{subsection}

\begin{subsection}{The Euler characteristic}
\label{subsection the euler characteristic}

\begin{defn}\label{defn Lefschetz number}
{\em
Let $X\in ObjDTOP({\N})$, 
$f:X\into X\in MorDTOP(\N)$ and let $R$ be a
commutative ring. Suppose that $dimX=m$.
The {\it Lefschetz number of $f$ over $R$}
is defined by $\lambda (f;R):=\sum _{i=0}^m(-1)^itr \, f^{*}_i$, where
$tr f^{*}_i$ is the trace of the $R$-module homomorphism
$f^{*}:H^{*}(X;R)\into H^{*}(X;R)$.
We denote $\lambda (f;\QQ)$ by $\lambda (f)$.  
}
\end{defn}

\begin{fact}\label{fact hopf trace theorem}
(Hopf Trace Theorem) $$\sum _{i=0}^m(-1)^itr \, f_{i*}=
\sum _{i=0}^m(-1)^itr\, f_{i\sharp }=
\sum _{i=0}^m(-1)^itr\, f^{\sharp}_i
=\sum _{i=0}^m(-1)^itr \, f^{*}_i.$$
\end{fact}

\pf
For the first equality see lemma 9.18 \cite{ro}, for the last
equality see theorem I.D.2 \cite{b} and for the second equality
see corollary I.D.3 \cite{b}. This shows that $\lambda (f;R)$ is an
integer and is independent of $R$.
\qed

Recall that the {\it o-minimal Euler characteristic} $E(X)$ of a
definable set $X$ of dimension $m$ is defined by 
$E(X):=\sum _{i=0}^m(-1)^ia_i$ where for each $i$, $a_i$ is the number
of cells of dimension $i$ in a cell decomposition of $X$. This number
does not depend on the chosen cell decomposition.

\begin{thm}\label{thm euler characteristic}
If $X\in ObjDTOP({\N})$ then $E(X)=\lambda (1_X)=\sum
_{p=0}^m(-1)^pb_p$ where, $b_p:=dim_{\QQ}$$H^p(X;\QQ)$ are the Betti numbers.
\end{thm}

\pf
By the triangulation theorem, we can assume that $X=|K|$ for  a simplicial
complex $K$. Its clear that
$E(X)=\sum _{p=0}^m(-1)^pa_p$ where for each $p$, $a_p$ is the number
$p$-simplexes in $X$. There is a canonical chain equivalence
$j_{\sharp}:C_*(K)\into S_*(|K|)$ 
and for each $p$ the rank of the $\ZZ$-module $C_{p}(K)$ is exactly
$a_p$ and so $\sum _{p=0}^m(-1)^pa_p=\lambda (1_X)$. The other
equality is clear.
\qed

We now combine theorem \ref{thm euler characteristic} with theorem 
\ref{thm hopf algebra of a group} to prove the following result.

\begin{thm}\label{thm euler of a group}
Let $G$ be a definably compact definable group and for each $k\ge 1$
let $p_{k}:G\into G$ be given by $p_k(x):=x^k$, $p_0(x)=1$ and
$p_{-k}(x)=(x^{-1})^k$. Then there is $r>0$
such that for all $l\in \ZZ$, 
$\lambda (p_l)=(1-l)^r$. In particular, $E(G)=0$. 
\end{thm}

\pf
By theorem \ref{thm hopf algebra of a group},
$H^*(G;\QQ)$ is isomorphic to the finitely generated free
exterior $\QQ$-algebra $\bigwedge [y_1,\dots ,y_r]_{\QQ}$ and therefore
comultiplication is given by 
$\mu (x)=x\otimes _{\QQ}$$1+1\otimes _{\QQ}$$x$. An
induction on $k$ shows that $(p_k)^{*}(x)=kx$ for $x\in \{y_1,\dots ,y_r\}$
(use the fact that if 
$\Delta :G\into G\times G$ is the
diagonal map, then $(p_{k+1})^*(x)=\Delta ^*((p_k)^*\otimes
_{\QQ}1_G^*)\mu (x)$, for details see theorem F.1
chapter III in \cite{b})
 and therefore, if
$x_i=y_{i_1}\cdots y_{i_u}$ is a generating monomial for $H^*(G;\QQ)$
then $(p_k)^*(x_i)=k^{len(x_i)}x_i$ where $len(x_i)=u$. This implies
that
$\lambda (p_k)=\sum (-1)^{deg(x_i)}k^{len(x_i)}+1$ where the sum is taken
over all monomials $x_i$ generating $H^*(G;\QQ)$. Since $y_j$'s have odd
degrees, we have $(-1)^{deg(x_i)}=(-1)^{len(x_i)}.$ Using this, a simple
calculation shows that $\lambda (p_k)=(1-k)^n$. For negative $k$ use
the same argument together with the fact that $\iota ^*(z)=(-1)^pz$
for all $z\in H^p(G;\QQ)$
where $\iota $ is the inverse in $G$. The case $k=0$ is trivial.

It will follow from remark \ref{nrmk orientation of G} that $r>0$. 
\qed

Combining theorem \ref{thm euler of a group}, the main theorem of
\cite{e1} together with results from \cite{s} we get.

\begin{cor}\label{cor main}
If $G$ is a definable abelian group such that $E(G)=0$ then for
each $k>1$ the subgroup $\{x\in G: kx=0\}$ of $k$-torsion points is finite
and non trivial.
\end{cor}

\pf
The fact that the subgroup $\{x\in G: kx=0\}$ of
$k$-torsion points is finite and non trivial follows from proposition
6.1 \cite{s} and lemma 5.11 \cite{s}.
\qed

\begin{cor}\label{cor euler zero and s}
(1) A definable solvable group has Euler characteristic zero if and
only if it has a definably compact part.
(2) If $G$ is a definable group such that $E(G)=0$ then it has a
definable abelian subgroup $H\neq 1$ such that $E(H)=0$. 
(3) Maximal definably compact definably connected definable abelian
subgroups of a definable group are conjugate.
\end{cor}

\pf
(1) follows from the main theorem of \cite{e1} and 
theorem \ref{thm euler of a group}. (2) is lemma 2.9 \cite{s} and  (3) is
corollary 5.19 \cite{s} since by (1) and corollary 5.17 \cite{s}
``maximal tori'' (in the sense of \cite{s}) are exactly 
maximal definably compact definably connected definable abelian
subgroups.
\qed

\end{subsection}
\end{section}
\begin{section}{Orientable definable manifolds}
\label{section orientable definable manifolds}

Through this section, ${\mathbf X}=(X,(X_i,\phi _i)_{i\in I})$ is a
definable manifold of dimension $n$.

\begin{subsection}{Lemma on definable triangulation}
\label{subsection lemma on definable triangulation}

Since $\mathbf{X}$ if affine, $X\subseteq N^k$ for some $k$ and so we
can definably triangulate the definable set $X$. But below we will be
interested in a special definable triangulation of $X$ compatible
with the definable charts $(X_i,\phi _i)_{i\in I}$. We make this
notion precise in the following definition:

\begin{defn}\label{defn good triangulation}
{\em
Let $A_1, \dots , A_n, B, Z$ be definable subsets of $X$.  Let $I=\{1,
2, \dots , k\}$ be a numbering of $I$. Define inductively 
$(X'_i, (\Psi _i, M_i), (\Phi _i, N_i))$ for $i\in I$ by:  $X'_1=X_1$, 
$(\Psi _1, M_1)$ is a definable triangulation of $\phi _1(X_1)$
compatible with the definable subsets $\phi _1(X_1\cap X_j)$, 
$\phi _1(X_1\cap A_l)$ and $\phi _1(X_1\cap X_j\cap A_l)$ for all 
$l\in \{1, \dots , n\}$ and all $j\in I_1$, and $(\Phi _1, N_1)=(\Psi
_1, M_1)$; let $X'_{i+1}:=\setminus \cup \{C:\phi _r(C)\in N_i, \,\,\,
r=1, \dots , i\}$ and be $(\Psi _{i+1}, M_{i+1})$ be a definable 
triangulation of $\phi _{i+1}(X'_{i+1})$ compatible with the definable 
sets $\phi _{i+1}(X'_{i+1}\cap X_j)$, $\phi _{i+1}(X'_{i+1}\cap A_l)$
and $\phi _{i+1}(X'_{i+1}\cap X_j\cap A_l)$ for all $l\in \{1, \dots ,
n\}$ and all $j\in I_{i+1}$ and such that $(\Phi _{i+1}, N_{i+1})$
which is equal to $(\Psi _{i+1}, M_{i+1})$ together with all the 
$(\Phi _{j|}, N_j)$ for $j\in I_{i+1}\cap \{1, \dots , i\}$ is a
definable triangulation of $\phi _{i+1}(X_{i+1})$. 

By a definable triangulation of $X$ compatible compatible with
the definable charts and with $A_1, \dots , A_n$ we mean a sequence 
$(\Phi , K)=\{(\Phi _i, N_i): i\in I\}$ some $\{(\Phi , N_i): i\in
I\}$ like above. By a definable triangulation of $B$ compatible with
the definable charts and with $A_1, \dots ,A_n$ we mean a sequence 
$(\Phi _B, K_B)=\{(\Lambda , L)\in (\Phi , K):|L|\subseteq B\}$ for
some definable triangulation $(\Phi , K)$ of $X$ compatible with the
definable charts and with $A_1, \dots , A_n, B$.  Note that
for each $C\in K$ such that $\phi _i(C)\in K_i$ for some $i\in I$, $C$
is definably homeomorphic to a $k$-simplex. If $C\subseteq Z$ we say
that $C$ is a $k$-simplex of $(\Phi , K)$ in $Z$. 
}
\end{defn}

\begin{lem}\label{lem good triangulation}
There is a finite cover of $X$ by open definable subsets definably 
homeomorphic to open balls in $N^n$. Moreover, if $A\subseteq X$ is a
definable closed subset then, there is a finite family $A_1,\dots A_l$ 
closed under intersection of closed definable subsets $A_i$ 
and a finite cover of $X$ by open definable subsets $U_j$ definably 
homeomorphic to open balls in $N^n$, such that $A=\cup A_i$ and
for each $i$ there is a $j_i$ such that $A_i\subseteq U_{j_i}$
\end{lem}

\pf
Let $(\Phi ,K)$ be a definable triangulation of $X$ compatible with
the definable charts. Clearly each $\phi _i(X_i)$ has a finite cover by
open definable subsets $V$ definably homeomorphic to open balls in
$N^n$: for each $0$-simplex $v$ of the induced definable triangulation of
$\phi _i(X_i)$ let $V$ be the star centred at $v$ (i.e, the union of
$\{v\}$ together with all open simplexes of the triangulation which
contain $v$ as a point in the closure. Now its also clear that the
definable open sets $\phi _i^{-1}(V)$ satisfy the lemma.

Now suppose that $A\subseteq X$ is a definable closed subset. Take a
definable triangulation of $X$ compatible with the definable charts
and with $A$.  For each $0$-simplex $v$ of such triangulation
let $U$ be the star centred at $v$ and let $A_i$ be the closed
simplexes obtained from the barycentric subdivision of the closed
simplexes contained in $A$. Then clearly each $A_i$ is contained in
the a star of a $0$-simplex of the triangulation contained in $A$.     
\qed

\end{subsection}

\begin{subsection}{The orientation sheaf of a definable manifold}
\label{subsection the orientation sheaf of a  definable manifold}

Note that the standard proof of the fact that $\tilde{H}_q(S^n)=\ZZ$ iff $q=n$
and zero otherwise, remains valid in our case.

\begin{lem}\label{lem local orientation}
For any point $x\in X$, $H_n(X,X-x;R)=R$.
\end{lem}

\pf
Suppose that $x\in X_i$.
Let $U$ be an open definable subset of $X_i$ such that $x\in U$
and 
$\phi _i(U)$ is an open ball in 
$\phi _i(X_i)$. Then 
$$H_n(X,X-x;R)\simeq H_n(U,U-x;R)\simeq \tilde{H}_{n-1}(U-x;R)\simeq 
\tilde{H}_{n-1}(S^{n-1})$$
and
$\tilde{H}_{n-1}(S^{n-1})\simeq R.$

The first equality is obtained by excising the closed definable subset $X-U$
of the definable open set $X-x$,
the second follows from the exact sequence of the pair $(U,U-x)$ since
$U$ is definably contractible, the third equality
is obtained from the fact that $U-x$ is definably homotopically equivalent
to $S^{n-1}$ and the last equality from remark above.
\qed

\begin{lem}\label{lem continuation}
Given an element $\alpha _x\in H_n(X,X-x;R)$ there is an open
definable neighbourhood $U$ of $x$ and $\alpha
\in H_n(X,X-U;R)$ such that $\alpha _x=j^U_x(\alpha )$, where 
$$j^U_x:H_n(X,X-U;R)\into H_n(X,X-x;R)$$ is the 
canonical homomorphism induced by inclusion.
\end{lem}

\pf
Let $a$ be the relative cycle representing $\alpha _x$. Then the
support 
$|\partial a|$ of $\partial a$ is a definably
compact subset of $X$ contained in $X-x$, so that $U:=X-|\partial
a|
$ is an open definable neighbourhood of $x$.
Let $\alpha \in H_n(X,X-U;R)$ be the homology class of $a$ relative to $X-U$.
\qed

We call $\alpha $ of lemma \ref{lem continuation} the {\it
continuation of 
$\alpha _x$ in $U$}.

\begin{lem}\label{lem locally constant}
There is a finite cover of $X$ by open definable subsets
$U_i$ definably homeomorphic to open balls in $N^n$
such that if $x\in U_i$ then for every $y\in U_i$, $j^{U_i}_y$
is an isomorphism (hence $\alpha _x$ has a unique continuation in $U_i$).
\end{lem}

\pf
The existence of a finite cover of $X$ by open definable subsets $U_i$  
definably homeomorphic to open balls in $N^n$ is proved in lemma
\ref{lem good triangulation}.

Let $V$ be one such open definably contractible definable subset  of $X$,
let $x\in V$ and let $H:V\times [0,1]\into V$ be a definable contraction
of $V$ to a point in $V$. Let $U:=V_t:=H(V,t)$ be such that $x\in V_t$
and $0<t<1$.
Then we have the following commutative diagram for any 
$y\in U$:
\[
\begin{array}{clcr}
H_n(X,X-U;R)\,\,\,\simeq \,\,\, H_n(V,V-U;R)\,\,\,\simeq \,\,\, \tilde{H}_{n-1}(V-U;R)\\
\downarrow
^{j^U_y}\,\,\,\,\,\,\,\,\,\,\,\,\,\,\,\,\,\,\,\,\,\,\,\,\,\,\,\,\,\,\,\,\,\,\,\,\,\,\,\,\,\,\,\,\,\,\,
\downarrow \,\,\,\,\,\,\,\,\,\,\,\,\,\,\,\,\,\,\,\,\,\,\,\,\,\,\,\,\,\,\,\,\,\,\,\,\,\,\,\,\,\,\,\,\,\,\,\,\,\,\,\downarrow \\
H_n(X,X-y;R)\,\,\,\simeq \,\,\, H_n(V,V-y;R)\,\,\,\simeq \,\,\, \tilde{H}_{n-1}(V-y;R)
\end{array}
\]
in which the left horizontal isomorphism are excisions and the right
ones are connecting homomorphisms 
($V$ is definably contractible). The right vertical arrow is an isomorphism 
because the inclusion $V-U\into V-y$ is a definable homotopy
equivalence (move out radially from $y$). Therefore,
$j^U_y$ is an isomorphism. Since $V=\cup _{0<t<1}V_t$ the lemma follows.
\qed

\begin{defn}\label{defn orientation presheaf}
{\em 
The {\it $R$-orientation presheaf} $O^X$ on $X$ is the presheaf of
$R$-modules  given by $O^X(U):=H_n(X,X-U;R)$ i.e., $O^X$ is a
contravariant functor from the category of open definable subsets of $X$
with the inclusion maps into the category of $R$-modules.
If $V\subseteq U$ are open definable subsets of $X$,
$$j^U_V:H_n(X,X-U;R)\into H_n(X,X-V;R)$$
denotes the homomorphism induced by inclusion. $j^U_V$ is called a
{\it restriction map}.
}
\end{defn}

Note that by lemma \ref{lem continuation} the {\it stalk} $O^X_x$ of
$O^X$ at $x\in X$, which is by definition, the direct limit
$O^X_x:=\lim _{x\in U}O^X(U)$ with respect to the restriction maps is 
$H_n(X,X-x;R)$. And  $j^U_x:H_n(X,X-U;R)\into H_n(X,X-x;R)$ is exactly
the natural map induced by the direct limit.

\begin{defn}\label{defn the etale space}
{\em
The {\it \'etal\'e space} $\tilde{O}^X$ associated to the
$R$-orientation presheaf $O^X$ is the topological space 
$$\tilde{O}^X:=\{(x,\alpha _x): x\in X,\, \alpha _x\in H_n(X,X-x;R)\}$$
where the basis for the topology on $\tilde{O}^X$ is given by 
$(U,\alpha _U):=\{(x,\alpha _x):x\in U,\,\, \alpha _x=j^U_x(\alpha _U)\}$
for  open definable subsets $U$ of $X$ together with the \'etal\'e map 
$p:\tilde{O}^X\into X$ given by $p(x,\alpha _x)=x$ (i.e., $p$ is
locally a homeomorphism).   
}
\end{defn}

For any definable subset $A\subseteq X$, a continuous map
$s:A\into \tilde{O}^X$ such that $p\circ s=1_A$ is called a 
{\it section over $A$}. A section over $X$ is called a {\it global section}.
The set $\Gamma (A;R)$ of all sections over $A$ is in a natural way an
$R$-module. We denote by $\Gamma _c(A;R)$ the $R$-submodule of $\Gamma (A;R)$
of all sections $s\in \Gamma (A;R)$ which agree with the zero section
outside some definably compact definable subset of $A$.   

We have a canonical homomorphism
$$j_A:H_n(X,X-A;R)\into \Gamma (A;R)$$ defined by $j_A(\alpha
)(x):=(x,j_x^A(\alpha ))$  for $x\in A$ (see remark 22.23 \cite{g}). 
If $B\subseteq A$, we have
the commutative diagram
\[
\begin{array}{clcr}
H_n(X,X-A;R)\,\,\,\stackrel{j_A}{\rightarrow}\,\,\, \Gamma (A;R)\\
\downarrow
^{j^A_B}\,\,\,\,\,\,\,\,\,\,\,\,\,\,\,\,\,\,\,\,\,\,\,\,\,\,\,\,\,\,\,\,\,\,\,\,\,\,\,\,\downarrow ^r\\
H_n(X,X-B;R)\,\,\,\stackrel{j_B}{\rightarrow} \,\,\, \Gamma (B;R)
\end{array}
\]
where $r$ is the restriction map.

\begin{defn}\label{defn orientation sheaf}
{\em
$\tilde{O}^X$ determines the presheaf of sections
$U\into \tilde{O}^X(U):=\Gamma (U;R)$.  And we have a morphism of
presheafs $j:O^X\into \tilde{O}^X$ given by $j_U:O^X(U)\into
\tilde{O}^X(U)$. The presheaf $\tilde{O}^X$ is actually a {\it sheaf}
i.e., for every collection of open definable subsets $U_i$ of $X$ with
$U=\cup _{i}U_i$ then $\tilde{O}^X$ satisfies (1) if $\alpha ,\beta
\in \tilde{O}^X(U)$ and $j^U_{U_i}(\alpha )=j^U_{U_i}(\beta )$ for all
$i$, then $\alpha =\beta $ and (2) if $\alpha _i\in \tilde{O}^X(U_i)$
and if for $U_i\cap U_j\neq \emptyset $ we have $j^{U_i}_{U_i\cap
U_j}(\alpha _i)=j^{U_j}_{U_i\cap U_j}(\alpha _j)$ for all $i$, then
there exists an $\alpha \in \tilde{O}^X(U)$ such that
$j^U_{U_i}(\alpha )=\alpha _i$ for all $i$.

The sheaf $\tilde{O}^X$ is called the {\it $R$-orientation sheaf} of $X$.
} 
\end{defn}

\end{subsection}

\begin{subsection}{Orientable definable manifolds}
\label{subsection orientable definable manifolds}
 
\begin{defn}\label{defn orientation along U}
{\em
If $x\in X$, an {\it $R$-orientation of $X$ at $x$} is a
generator $\alpha _x$ of the $R$-module $H_n(X,X-x;R)$.
Given a definable subset $A\subseteq X$, an {\it $R$-orientation of
$X$ along $A$} is a section $s\in \Gamma (A,R)$ such that for each
$a\in A$, $s(a)$ is an $R$-orientation of $X$ at $a$. An 
{\it $R$-orientation of $X$} is an $R$-orientation of
$X$ along $X$. $X$ is
{\it $R$-orientable along $A$} if such $s$ exists, $X$ is {\it
$R$-orientable} if $X$ is $R$-orientable along $X$. We also say that
$X$ is {\it orientable} if it is $\ZZ$-orientable.
}
\end{defn}

Note that, if there is an element $\alpha \in H_n(X,X-A;R)$ such that
$j^A_y(\alpha )$ generates $H_n(X,X-y;R)$ for each $y\in A$ 
then $s(a):=(a,j_a^A(\alpha ))$ is an $R$-orientation of
$X$ along $A$. We call such $\alpha $ an {\it $R$-orientation of
$X$ along $A$}.
If there is a family $(U_i,\alpha _i)$ with $U_i$'s open definable 
neighbourhoods
which cover $X$ and each $\alpha _i$'s an $R$-orientation of
$X$ along $U_i$ which satisfies the compatibility
conditions: for any
$x\in X$, if $x\in U_i\cap U_j$, then $j^{U_i}_x(\alpha _i)
=j^{U_j}_x(\alpha _j)$ then, $s\in \Gamma (X,R)$ given by $s(x):=
(x,j_x^{U_i}(\alpha _i))$ if $x\in U_i$ is an $R$-orientation of $X$.
We call such family $(U_i,\alpha _i)$ is called 
an {\it $R$-orientation system on $X$}.

By the universal coefficient theorem, if $X$ is orientable then it is 
$R$-orientable for all coefficient rings $R$.

\begin{nrmk}\label{nrmk orientation of G}
{\em
Let $G$ be an definably connected, definably compact definable group of
dimension $n$. 
Then $G$ is $R$-orientable: 
let $\alpha _1$ be an $R$-orientation in an open definable neighbourhood
$U_1$ of 
$1$ (the identity element of $G$).
For $a\in G$ let $l_a:G\into G$ be $l_a(x):=ax$, $U_a:=l_a(U_1)$ and 
$\alpha _a:=l_{a*}(\alpha _1)$. Then 
$(U_a,\alpha _a)_{a\in G}$ is an $R$-orientation system for $G$. This
follows from the fact that for $x\in U_1$ and $a\in G$ we have 
$l_{a*}\circ j^{U_1}_x=j^{U_a}_{ax}\circ l_{a*}$.
}
\end{nrmk}

Below, we talk of {\it definable covering maps}, for details see \cite{e2}. 

\begin{thm}\label{thm non orientable}
(1) If $X$ is definably connected and non-orientable then there
is a 2-fold definable covering map
$p:E\into X$ such that $E$ is a definable orientable definably connected 
manifold. In particular, every definably connected 
manifold whose fundamental group contains no
subgroup of index 2 is orientable.
(2) An open definable submanifold of an $R$-orientable $X$ is
$R$-orientable. In particular, $X$ 
is $R$-orientable if and only if all
its definably connected components are $R$-orientable.
(3) If $X$ is definably connected then any two $R$-orientations on $X$ that
agree at one point are equal. In particular, if
$X$ is orientable then it has exactly two distinct orientations.
(4) Every definable manifold has a unique $\ZZ$$/2\ZZ$-orientation.
\end{thm}

\pf
$(1)$ Let $E$ be the set of pairs $(x,\alpha _x)$, where $x\in X$ and 
$\alpha _x$ is one of the generators of
$H_n(X,X-x;\ZZ)$. $E$ is a definable manifold with a basis of open 
definable neighbourhoods given by
$(U,\alpha _U)$, where $U$ is an open definable neighbourhood in $X$ and 
$\alpha _U$ is a local orientation 
of $X$ along $U$. To see this use lemma \ref{lem continuation}. The
definable map 
$p:E\into X$ given by $p(x,\alpha _x)
=x$ is a $2$-fold definable covering map, and $E$ is oriented by the system
whose opens are the $(U,\alpha _U)$ and 
whose orientations are the elements $H_n(p_{|(U,\alpha
_U)})^{-1}(\alpha _U)$ in 
$H_n(E,E-(U,\alpha _U);\ZZ)$. Moreover, $E$ is definably connected for
otherwise, $p$ induces a homeomorphism  between
$X$ and a component of $E$ contradicting the orientability of $E$ and
its components. The rest of (1) 
follows from the fact that $p_*\pi _1(E,e_0)$ has index $2$ in 
$\pi _1(X,x_0)$.

(2) Follows from the fact that if $V\subseteq X$ is an open
definable subset and $x\in V$ then $H_n(V,V-x;R)\simeq H_n(X,X-x;R)$ and
definably connected components of $X$ are open definable subsets.

(3) Follows from lemma \ref{lem locally constant} since the set on
which two $R$-orientations agree is a clopen definable subset of $X$.

(4) Follows also from lemma \ref{lem locally constant}.
\qed

\begin{nrmk}\label{nrmk orientation of the product}
{\em
Let $X$ and $Y$ be $R$-orientable definable manifolds. Then an 
$R$-orientation of $X$ and of $Y$ determine in a 
canonical way an $R$-orientation of $X\times Y$. If fact we have by
the 
relative K\"unneth formula for homology
a continuous map $\mu :\tilde{O}^X\times \tilde{O}^Y\into \tilde{O}^
{X\times Y}$ commuting with the 
natural maps $p^{X\times Y}:\tilde{O}^{X\times Y}\into X\times Y$ and 
$p^X\times p^Y:\tilde{O}^X\times \tilde{O}^Y
\into X\times Y$. This map determines a map $\Gamma (A;R)\times 
\Gamma (B;R)\into \Gamma (A\times B;R)$ for every
definable subsets $A\subseteq X$ and $B\subseteq Y$. For details
see 
VII.2.13 \cite{d}.
}
\end{nrmk}

\begin{nrmk}\label{nrmk orientation of the boundary}
{\em
Suppose now that $X$ is a definable manifold with boundary. Then
just like 
in 28.7\cite{g} and 28.8\cite{g}, 
we have that if $V$ is an open definable subset of $X$ and
$\mbox{\. V}=
V\cap \mbox{\. X}$ and 
$\partial V=V\cap \partial X$ then there are unique homomorphisms
$\partial _V
:\Gamma (\mbox{\. V})\into \Gamma
(\partial V)$ which are compatible with restriction to smaller $V$ and
which 
take local $R$-orienattion 
of $\mbox{\. X}$ along $\mbox{\. V}$ into a local $R$-orientation of 
$\partial X$ along $\partial V$. In particular,
if $\mbox{\. X}$ is $R$-orientable then so is $\partial X$.

Note that by 28.12 \cite{g}, if $\mbox{\. X}$ is $R$-orientable, then
there 
is a unique $R$-orientation of
$2X$ (see remark \ref{nrmk with boundary}) inducing the given
$R$-orientation 
on 
$\mbox{\. X}^1$ and $\mbox{\. X}^2$.
}
\end{nrmk}

\end{subsection}

\begin{subsection}{The fundamental class}
\label{subsection the fundamental class}

We start this subsection with the following easy remark.

\begin{nrmk}\label{nrmk orientations and sections}
For any definable subset $A\subseteq X$, $X$ is $R$-orientation along $A$ iff 
the covering map $p:p^{-1}(A)\into A$ is trivial. In which case,
$\Gamma (A;R)$ is isomorphic to the
$R$-module of all continuous maps $A\into R$. In particular, if $A$
has $k$ connected components then $\Gamma (A;R)\simeq R^k$.
\end{nrmk}

\pf
Suppose that $s\in \Gamma (A;R)$ is an $R$-orientation of $X$ along
$A$. Then for each 
$a\in A$, $s(a)=(a,s'(a))$ and $s'(a)$
is an $R$-orientation at $a$. If $(x,\alpha _x)\in p^{-1}(A)$ then
there is a unique $\lambda _x\in R$ such that $\alpha _x=\lambda
_xs'(x)$. 
The map
$\phi :p^{-1}(A)\into A\times R$ given by $\phi (x,\alpha
_x)=(x,\lambda _x)$ 
is a 
homeomorphism by lemma \ref{lem continuation}. Conversely, given $\phi
$ we 
can recover 
$s$ by $s(x)=\phi ^{-1}(x,1)$ for $x\in A.$  
\qed

The main result of this subsection is the following theorem.

\begin{thm}\label{thm fundamental class}
Suppose that $A\subseteq X$ is a closed definable subset. Then, for all $q>n$
$H_q(X,X-A;R)=0$ and
$$j_A:H_n(X,X-A;R)\into \Gamma _c(A;R)$$
is an isomorphism.
\end{thm}

\pf
By lemma \ref{lem good triangulation}, there is a finite family
$A_1,\dots A_l$ closed under intersection of closed definable subsets $A_i$ 
and finitely many open definable subsets $U_j$ definably homeomorphic
to open balls such that $A=\cup A_i$ and for each $i$ there is a $j_i$ 
such that $A_i\subseteq U_{j_i}$.

\medskip
{\it Claim (1):} If the result holds for closed definable subsets $A$, $B$
and 
$A\cap B$ then it holds for $C:=A\cup B.$

\medskip
{\it Proof of Claim (1):} 
Using the relative Mayer-Vietoris sequence for the triad $(X, X-A, X-B)$ we get
$H_q(X,X-C;R)=0$ for $q>n$ and we have the commutative diagram
\[
\begin{array}{clcr}
0\,\,\,\,\into \,\,\,G_n(C)\,\,\,\,\,\,\into \,\,\,\,\,\,\,G_n(A)\oplus G_n(B)\,\,\,\,\,\,\,\into
\,\,\,\,G_n(A\cap B)\\
\,\,\,\,\,\,\,\,\,\,\,\,\,\,\,\,\,\downarrow ^{j_C}\,\,\,\,\,\,\,\,\,\,\,\,\,\,\,\,\,\,\,\,\,\,\,\,\,\,\,\,\,\,\,\,\,\,\,\,\,\,\,\,\,\,\,\downarrow ^
{j_A\oplus j_B}\,\,\,\,\,\,\,\,\,\,\,\,\,\,\,\,\,\,\,\,\,\,\,\,\,\,\,\,\,\,\,\,\,\,\,\,\,\,\,\,\,\,
\downarrow ^{j_{A\cap B}}\\
0\into \Gamma _c(C;R)\stackrel{(r_A,r_B)}{\rightarrow }\Gamma
_c(A;R)\oplus 
\Gamma _c(B;R)\stackrel{r_A+r_B}{\rightarrow }\Gamma _c(A\cap B;R)
\end{array}
\] 
where $G_n(Z):=H_n(X,X-Z;R)$,
and chasing the diagram shows that $j_C$ is an isomorphism.
\qed

\medskip
{\it Claim (2):} The result holds for $A$ definably compact.

\medskip
{\it Proof of Claim (2):} Note that since $A$ is definably compact, it is a
closed 
complex i.e.,
it contains all its faces, and also $A_i$'s are closed simplexes (all
its 
faces are in $X$).
By claim (1), the theorem follows by induction on the number of
$A_i$'s such 
that $A=\cup A_i$. So suppose that $A=A_1$ and let $A\subseteq
U_{j_1}=:U$. 
Since
$U$ is definably homeomorphic to an open ball
of dimension $n$ the result follows from $H_q(U, U-A;R)\simeq
\tilde{H}_{q-1}
(U-A;R)$$
\simeq \tilde{H}_{q-1}(S^{n-1};R)$, $\tilde{H}_{n-1}(S^{n-1};R)=R$ and
remark \ref{nrmk orientations and sections}. 
\qed

\medskip
{\it Claim (3):} If $A\subseteq U$, where $U$ is an open definable subset
of $X$ 
with definably compact closure $\bar{U}$ then the result holds for $U$ and $A$.

\medskip
{\it Proof of Claim (3):} 
We use the exact homology sequence for the triple $(X,U\cup
(X-\bar{U}),$$
(U-A)\cup (X-\bar{U}))$. Note that by excision, $H_q(U,U-A;R)\simeq 
H_q(U\cup (X-\bar{U}),(U-A)\cup (X-\bar{U});R)$.
For $q>n$ we have
$H_{q+1}(X,U\cup (X-\bar{U});R)\into H_q(U,U-A;R)\into H_q(X,(U-A)\cup 
(X-\bar{U});R).$

For $q=n$, we have
\[
\begin{array}{clcr}
0\into H_n(U,U-A;R)\into H_n(X,W;R)\into
H_n(X,V;R)\\
\,\,\,\,\,\,\,\,\,\,\,\,\,\,\,\,\downarrow
^{j^U_A}\,\,\,\,\,\,\,\,\,\,\,\,\,\,\,\,\,\,\,\,\,\,\,\,\,\,\,\,\,\,\,\,\,\,\,
\,\,\,\,\,\,\,\,\,\,\,\,\downarrow
\,\,\,\,\,\,\,\,\,\,\,\,\,\,\,\,\,\,\,\,\,\,\,\,\,\,\,\,\,\,\,\,\,\,\,\,\,\,\,
\,\,\,\,\,\,\,\,\,\,\,\,\downarrow \\
0\,\,\,\into \,\,\,\Gamma _c(A;R)\,\,\stackrel{i}{\rightarrow
}\,\,\Gamma (\bar{A}\cup (\bar{U}-U);R)\,\,\stackrel{r}{\rightarrow
}\,\,\Gamma (\bar{U}-U;R)
\end{array}
\]
where $W:=(U-A)\cup (X-\bar{U})$, $V:=U\cup 
(X-\bar{U})$, $\Gamma _c(A;R)$ and $j_A$ are computed in the definable
manifold $U$, and
the monomorphism $i$ is defined as follows: Let $s\in \Gamma _c(A;R)$
be zero 
outside
an definably compact $K\subseteq A$. Then $i(s)_{|A}=s$ and $i(s)=0$
outside  $K$. 

Applying the result for the definably compacts $\bar{U}-U$ and $\bar{A}\cup 
(\bar{U}-U)$ 
we see that $H_q(U,U-A;R)=0$ for $q>n$ and $j^U_A$ is an isomorphism. 
\qed

\medskip
{\it Claim (4):} The result holds for any closed definable subset $A$.

\medskip
{\it Proof of Claim (4):} Given $s\in \Gamma _c(A;R)$ zero outside the 
definably compact $K\subseteq A$. There is an open definable subset
$K\subseteq U$ 
such that $\bar{U}$ is definably compact: As we have seen, by the
triangulation 
theorem, $X$ is covered by finitely many open definable subsets $U_i$'s 
definably homeomorphic to open balls. We get $U$ by taking it to be
$\cup _iV_i$ where $V_i$ is obtained from $U_i$ after a small
definable contraction.

Consider $A'=A\cap U$, $s'=s_{|A}.$ By Claim (3) applied to $U$ and
$A'$, and 
the commutative diagram
\[
\begin{array}{clcr}
H_n(U,U-A';R)\into H_n(X,X-A;R)\\
\,\,\,\,\,\,\,\,\,\downarrow ^{j_{A'}}\,\,\,\,\,\,\,\,\,\,\,\,\,\,\,\,\,\,\,\,\,\,
\,\,\,\,\,\,\,\,\,\,\,\,\downarrow ^{j_A}\\
0\into s'\in \Gamma _c(A';R)\,\,\,\,\,\,\stackrel{i}{\rightarrow
}\,\,\,\,\,\,s\in \Gamma _c(A;R)
\end{array}
\]
we see that $j_A$ is surjective. 

Now let $\alpha \in H_q(X,X-A;R).$ If $q=n $ suppose that $j_A(\alpha
)=0.$ 
Let $z$ be the relative cycle representing $\alpha $. Applying the
above 
argument to $|z|$, there is an open definable subset $|z|\subseteq U$
such that 
$\bar{U}$ is definably compact. Let $A'=A\cap U$. By the same commutative
diagram, 
we have $\alpha =0$. For $q>n$, we know that the class of $z$ in 
$H_q(U,U-A';R)$ is zero by Claim (3), so $\alpha =0$. 
\qed

\begin{defn}\label{defn fundamental class}
{\em
We see that, an $R$-orientation of $X$ along a definably compact 
$A$ determines a
generator $\zeta _{X,A}$ of $H_n(X,X-A;R)$ called the {\it
fundamental class} of the $R$-orientation of $X$ along $A$. If $A=X$ we let
$\zeta _X:=\zeta _{X,X}$. 

The {\it orientation class of $X$} is the element $\omega _X\in
H^n(X;R)$ such that $(\zeta _X, \omega _X)=1$.
}
\end{defn}

\begin{cor}\label{cor vanishing}
If $A$ is a closed definably connected and not definably compact
definable subset of $X$ then
$H_n(X,X-A;R)=0$. In particular, if $X$ is definably connected
and not definably compact then $H_n(X;R)=0$ . 
\end{cor}

\pf
$j_A$ is locally constant and zero outside an definably compact subset of
$A$. By definably connectedness of $A$, if $\alpha \in H_n(X,X-A;R)$,
then $j_A(\alpha )$ is zero and so $\alpha =0.$
\qed
      
\begin{cor}\label{cor fundamental class}
Suppose that $X$ is definably compact and definably connected. 
Assume that for any $a\neq 0$, $
a\in R$ and any 
unity $u\in R$, $ua=a$ implies $u=1$. Then $H_n(X;R)=R$ if $X$ is 
$R$-orientable and $H_n(X;R)=0$ otherwise.
\end{cor}

\pf
If $X$ is $R$-orientable, apply remark \ref{nrmk orientations and sections}.
Supppose that there is a global section $s\in \Gamma (X;R)$, $s\neq
0$. Then 
by lemma
\ref{lem continuation} there is $a\in R$, $a\neq 0$, such that $s'(x)$
is $a$ 
times a
generator of $H_n(X,X-x;R)$ for all $x\in X$, where $s(x)=(x,s'(x)).$
The 
hypothesis on
$R$ implies that $s'(x)/a$ is a well defined generator, and so $s/a\in 
\Gamma (X;R)$ given
by $(s/a)(x):=(x,s'(x)/a)$ is an $R$-orientation of $X.$
\qed

\end{subsection}

\begin{subsection}{Degrees}
\label{subsection degrees}

\begin{defn}\label{defn degrees}
{\em
Let $f:X\into Y$ be a continuous definable map between orientable 
definable manifolds of dimension $n$,
and let $K\subseteq Y$ be a definably compact definably connected
non-empty definable subset such that 
$f^{-1}(K)$ is definably compact. The {\it degree of $f$ over
$K$} is the integer $deg_Kf$
such that $f_*(\zeta _{X,f^{-1}(K)})=(deg_Kf)\zeta _{Y,K}$.
}
\end{defn}

We say that $f$ is a {\it definably proper} map if for all definably
compact 
$K\subseteq Y$, $f^{-1}(K)$ is 
definably compact. For example $f$ is a definably proper map if $X$ is 
definably compact.

\begin{lem}\label{lem degrees}
Let $f:X\into Y$ and $K\subseteq Y$ be as in definition \ref{defn
degrees}. 
Then we have:
(1) if $K'\subseteq f^{-1}(K)$ is 
definably compact then $f_*(\zeta _{X,K'})=(deg_Kf)\zeta _{Y,K}$;
(2) if 
$L$ is a
definably compact definable subset of $K$ then we have
$f_*(\zeta _{X,f^{-1}(L)})=(deg_Kf)\zeta
_{Y,K}$ in 
particular $deg_Lf=deg_Kf$;
(3) if $X$ is a finite union of open definable subsets $X_1,\dots
, X_r$ 
such that the sets
$K_i=f^{-1}(K)\cap X_i$ are mutually disjoint then $deg_Kf=\sum
_{i=1}^r
deg_{K_i}f_{|K_i}$;
(4) if $g:Z\into X$ is a continuous definable map between
orientable 
definable manifolds of dimension $n$,
and $g^{-1}(f^{-1}(K))$ is definably compact then $deg_K(f\circ
g)=(deg_{f^{-1}(K)}g)deg_Kf$;
(5) if $f$ is a definably proper map and if $Y$ is definably connected then
$deg_kf=
degf$ and is independent from $K$. 
\end{lem}

\pf
(1) The inclusion homomorphism $$H_n(X,X-K';R)\into
H_n(X,X-f^{-1}(K);R)$$ takes 
$\zeta _{X,K'}$ into $\zeta _{X,f^{-1}(K)}$. Therefore, the composition with 
$$f_*:H_n(X,X-f^{-1}(K);R)\into H_n(X,X-K;R)$$ takes $\zeta _{X,K'}$ into 
$f_*(\zeta _{X,f^{-1}(K)})$$=(deg_K)\zeta _{X,K}.$

(2) Consider the commutative diagram
\[
\begin{array}{clcr}
H_n(X,X-f^{-1}(K);R)\stackrel{f_*}{\rightarrow }H_n(X,X-K;R)\\
\,\,\,\downarrow
^{i_*}\,\,\,\,\,\,\,\,\,\,\,\,\,\,\,\,\,\,\,\,\,\,\,\,\,\,\,
\,\,\,\,\,\,\,\,\,\,\, \downarrow ^{i_*}\\
H_n(X,X-f^{-1}(L);R)\stackrel{f_*}{\rightarrow }H_n(X,X-L;R)
\end{array}
\]
Chasing $\zeta _{X,f^{-1}(K)}$ through the diagram and using (1) gives
the 
result.

(3) Consider the maps
$$\oplus _{j=1}^rH_n(X_j,X_j-K_j;R)\stackrel{\oplus _{j=1}^ri^j_*}
{\rightarrow }
H_n(X,X-f^{-1}(K);R)\stackrel{i^Q_*}{\rightarrow }H_n(X,X-Q;R),$$
where the $i^j$ are inclusions and $Q\in f^{-1}(K).$  We  have
$$i^Q_*(\oplus _{j=1}^ri^j_*(\oplus _{j=1}^r\zeta _{X,K_j}))=\zeta
_{X,Q}$$ for 
every
$Q\in f^{-1}(K)$, hence $\oplus _{j=1}^ri^j_*(\oplus _{j=1}^r\zeta _{X,K_j})=
\zeta _{X,f^{-1}(K)}$ by (1). Now, we have
\[
\begin{array}{clcr}
(deg_Kf)\zeta _{X,K}=f_*(\zeta _{X,f^{-1}(K)})=f_*(\oplus
_{j=1}^ri^j_*
(\oplus _{j=1}
\zeta _{X,K_j}))=\\
\oplus _{j=1}^rf_{|K_j*}(\oplus _{j=1}^r\zeta _{X,K_j})=
\sum _{j=1}^rf_{|K_j*}(\zeta _{X,K_j})=(\sum _{j=1}^rdeg_Kf_{|K_j})
\zeta _{X,K}. 
\end{array}
\]

(4) Is obvious. (5) Let $K,L\subseteq Y$ be definably compact
and let 
$k\in K$ and $l\in L$ then by (2)
$deg_Kf=deg_kf$ and $deg_Lf=deg_lf$. Since $Y$ is definably connected there
is a 
definable path 
$\alpha :[0,1]\into Y$ from $k$ to $l$ (see \cite{e2}). $\alpha
([0,1])$ is 
definably compact and by (2) again
$deg_kf=deg_{\alpha ([0,1])}f=deg_lf$. 
\qed

\end{subsection}
\end{section}

\begin{section}{Duality on definable manifolds}
\label{section dualiy on definable manifolds}

Through this section, ${\mathbf X}=(X,(X_i,\phi _i)_{i\in I})$ is an
$R$-orientable  definable manifold of dimension $n$ with
$R$-orientation $s\in \Gamma (X,R)$.
Otherwise we take $R=\ZZ$$/2\ZZ$ and take the unique $\ZZ$$/2\ZZ$-orientation.

\begin{subsection}{Poincar\'e duality}
\label{subsection poincare duality}

\begin{defn}\label{defn singular cohomology with compact supports}
{\em
The {\it singular cohomology with compact supports} is defined to be
the direct limit $$H^q_c(X;R):=\lim _{K}H^q(X,X-K;R)$$ 
where $K$ are definably compact. Note that if $X$ is definably
compact then $H^q_c(X;R)=H^q(X;R)$. Definably proper definable maps 
$f:X\into Y$ induce homomorphisms $H^q_c(f):H^q_c(Y;R)\into
H^q_c(X;R)$. 

Now for any definably compact $K\subseteq X$,
consider the homomorphism $\zeta _{X,K}\cap :H^q(X,X-K;R)\into
H_{n-q}(X;R)$ given by cap product. Passing to the limit, these
homomorphisms give a homomorphism
$$D:H^q_c(X;R)\into H_{n-q}(X;R).$$
}
\end{defn}

\begin{thm}\label{thm poincare duality}
(Poincar\'e Duality Theorem).
$D:H^q_c(X;R)\into H_{n-q}(X;R)$ is an isomorphism (for all q).
\end{thm}

\pf
First note that, by lemma \ref{lem good triangulation}, 
$X$ is the union of a finite family $U_1,\dots ,U_k$ closed under
intersection of open definable subsets
definably homeomorphic to open balls. We prove the result by
induction on $k$.

Case $k=1$: Let $U_1=U$. Since $U$ is definably homeomorphic to an open
ball $B$ 
of dimension $n$
centred 
at the origin and of radius $1$, in computing the inductive limit 
$\lim _{K}H^q(B,B-K;R)$ it suffices to let $K$ run through the final
system 
of closed balls of radius $<1$ centred at the origin. But for such
$K$, the 
modules in question are zero unless $q=n$, and
$\zeta _{K}\cap :H^n(B,B-K;R)\into H_0(B;R)\simeq R$ is an isomorphism
and so 
the limiting homomorphism is also an isomorphism.

\medskip
{\it Claim:} If the result holds for open definable subsets $U$, $V$ and 
$W=U\cap V$, then the result holds for $Y=U\cup V$.

\medskip
{\it Proof of Claim:} Let $K$ (resp., $L$) be an definably compact subset
of $U$ 
(resp., $V$). We use the Mayer-Vietoris sequence for the triad
$(Y,Y-K,Y-L)$. 
The diagram (where we put $W':=W-K\cap L$, $Y':=Y-K\cup L$,$U':=U-K$,
$V':=V-L$; and for $Z=W,Y,U,V$ we set $T^m(Z):=H^m(Z,Z';R)$,
and $k:=n-q+1$ and $l:=n-q$),
\[
\begin{array}{clcr}
T^{q-1}(W)\,\,\,\,\,\,\into \,\,\,\,\,\,\,\,
T^q(Y)\,\,\,\,\,\,\into \,\,\,\,T^q(U)\oplus
T^q(V)\,\,\,\,\into \,\,\,\,\,T^q(W) \\
\,\,\,\,\,\,\downarrow ^{\zeta _{K\cap L}\cap }
\,\,\,\,\,\,\,\,\,\,\,\,\,\,\,\,\,\,\,\,\,\,\,\,\,\,\,\,\,\downarrow ^
{\zeta _{K\cup L}\cap }\,\,\,\,\,\,\,\,\,\,\,\,\,\,\,\,\,\,\,\,\,\,\,\,\,
\downarrow ^{\zeta _K\cap \oplus \zeta _L\cap }
\,\,\,\,\,\,\,\,\,\,\,\,\,\,\,\,\,\,\,\,\,\,\,\,\,\,\,\,\,\,\,\,\,\,
\downarrow ^{\zeta _{K\cap L}\cap }\\
H_{k}(W;R)\,\,\into \,\,H_{l}(Y;R)\,\,\into
\,\,H_{l}(U;R)\oplus H_{l}(V;R)\,\,\into \,\,H_{l}(W;R) 
\end{array}
\]
is commutative except possibly a $+$ or $-$ sign. Moreover, every 
definably compact in $Y$ has the form $K\cup L$. Passing to the limit gives
a sign 
commutative diagram (where $k=n-q+1$and $l=n-q$)
\[
\begin{array}{clcr}
H_c^{q-1}(W;R)\into H^q_c(Y;R)\into H_c^q(U;R)\oplus H_c^q(V;R)\into 
H_c^q(W;R)\\
\,\,\,\,\,\,\,\,\,\downarrow ^{D}\,\,\,\,\,\,\,\,\,\,\,\,\,\,\,\,\,\,\,\,\,\,\,\,\,\,\,\,
\,\,\,\,\,\downarrow ^{D}\,\,\,\,\,\,\,\,\,\,\,\,\,\,\,\,\,\,\,\,\,\,\,\,\,\,
\,\,\,\,\,\,\,\,\,\,\,\,\,\downarrow ^{D\oplus D}\,\,\,\,\,\,\,\,\,\,\,\,\,\,\,\,\,\,\,\,\,\,\,\,
\,\,\,\,\,\,\,\,\,\,\,\,\,\,\downarrow ^D\,\,\,\,\,\\
H_{k}(W;R)\into H_{l}(Y;R)\into H_{l}(U;R)\oplus
H_{l}(V;R)\into 
H_{l}(W;R)
\end{array}
\]
in which all the rows are exact and all the vertical arrows except
those 
involving $Y$ are isomorphisms. By the $5$-lemma (19.12 \cite{g}), the
result 
holds for $Y$.
\qed

We get from Poincar\'e duality theorem the usual corollaries. In particular,
suppose that $X$ is definably compact orientable. 
If $dim(X)$ is odd then $E(X)=0$ and if $dim(X)$ is even and not
divisible by $4$ then $E(X)$ is even (see 26.10 and 26.11
\cite{g}). Also if $X$ is definably compact with boundary and
$\mbox{\. X}$ is orientable, then $E(\partial X)$ is even (see 28.13
\cite{g}).

\end{subsection}

\begin{subsection}{Alexander duality}
\label{subsection alexander duality}

Through this subsection, $A\subseteq X$ is a closed definable subset and
$U:=X-A$. Let $$\mbox{{\it \u H}}^q(A;R)=\lim _{V}H^q(V;R)$$ be the direct
limit where $V$ are open definable neighbourhoods of $A$ directed by
reverse inclusion.

\begin{lem}\label{lem exact sequence for alexander duality}
Assume that $X$ is definably compact. Then there is an exact
sequence
$$\cdots \rightarrow H^q_c(U;R)\stackrel{i}{\rightarrow}H^q(X;R)
\stackrel{j}{\rightarrow}\mbox{{\it \u H}}^q(A;R)\stackrel{\delta }{\rightarrow}
H^{q+1}_c(U;R)\rightarrow \cdots $$
\end{lem}

\pf
The homomorphism $H^q_c(U;R)\stackrel{i}{\rightarrow}H^q(X;R)$ is the
unique homomorphism making the diagram
\[
\begin{array}{clcr}
H^q_c(U;R)\,\,\,\,\,\,\,\,\stackrel{i}{\rightarrow}\,\,\,\,\,\,\,\, H^q_c(X;R)\\
\downarrow
\,\,\,\,\,\,\,\,\,\,\,\,\,\,\,\,\,\,\,\,\,\,\,\,\,\,\,\,\,\,\,\,\,\,\,\,\,\,\,\,\,\,\,\,\,\,\,\downarrow \\
H^q(U,U-K;R)\,\,\,{\rightarrow}\,\,\, H^q(X, X-K;R)
\end{array}
\]
commutative for all $K\subseteq U$ definably compact, where 
$$H^q(U,U-K;R){\rightarrow}H^q(X, X-K;R)$$ is the inverse to the
excision isomorphism.

$\mbox{{\it \u H}}^q(A;R)\stackrel{\delta }{\rightarrow}H^{q+1}_c(U;R)$ is
induced by the homomorphisms 
$$H^q(V;R)\stackrel{\delta }{\rightarrow}H^{q+1}_c(X,V;R)\simeq
H^{q+1}(U,U-K;R)$$
where $V$ is an open definable neighbourhood of $A$ and $K:=X-V$ is a
definably compact contained in $U$.

To prove that this sequence is exact is a simple diagram chasing, for
details see theorem 27.3 \cite{g}.
\qed

Now let $\zeta _{X,A}$ be the fundamental class determined by the
$R$-orientation of $X$ along $A$. 
For any open definable neighbourhood $V$ of $A$, we
have $H_n(V,V-A;R)\simeq H_n(X,X-A;R)$ by excision; the pre-image of
$\zeta _{X,A}$ under this isomorphism is still denoted by $\zeta _
{X,A}$. Taking the cap product with $\zeta _{X,A}$ gives a
homomorphism
$$\zeta _{X,A}\cap :H^q(V;R)\into H_{n-q}(V,V-A;R)\simeq H_{n-q}(X,X-A;R)$$
which induces a homomorphism
$$D_A:\mbox{{\it \u H}}^q(A;R)\into H_{n-q}(X,X-A;R).$$

\begin{thm}\label{thm alexander duality}
(Alexander Duality Theorem).
Assume that $X$ is definably compact. Then $D_A$ is an
isomorphism for all $q$.
\end{thm}

\pf
The diagram
\[
\begin{array}{clcr}
\rightarrow H^q_c(U;R)\,\,\,\,{\rightarrow}\,\,\,\,H^q(X;R)\,\,\,\,\,{\rightarrow}\,\,\,\,\,\mbox{{\it \u H}}^q(A;R)\,\,\,\,\,\,{\rightarrow} \\
\,\,\,\downarrow ^{D_U}\,\,\,\,\,\,\,\,\,\,\,\,\,\,\,\,\,\,\,\,\,\,\downarrow ^{D_X}\,\,\,\,\,\,\,\,\,\,\,\,\,\,\,\,\,\,\,\,\,\,\,\downarrow ^{D_A}\\
\rightarrow H_{n-q}(U;R){\rightarrow}H_{n-q}(X;R){\rightarrow}H_{n-q}(X,X-A;R)
{\rightarrow}
\end{array}
\]
is sign-commutative (where $D_U$ and $D_X$ are the isomorphisms of the
Poincar\'e Duality Theorem). Now apply the $5$-lemma (19.15 \cite{g})
which still works when the diagrams are only sign commutative. 
\qed

\end{subsection}

\begin{subsection}{General separation theorem}
\label{subsection general separation theorem}

The following lemma is an adaptation of results from 26.17 \cite{g} on 
absolute neighbourhood retracts (namely, lemma 26.17.6 \cite{g}).

\begin{lem}\label{lem anr}
Let $Y$ be a definable manifold, $B$ a closed definable 
subset of $Y$ and $f:B\into X$ a continuous definable map. Then $f$
extends 
to a continuous definable map from an open definable neighbourhood of
$B$ in $Y$.
\end{lem}
 
\pf
By lemma \ref{lem good triangulation}, $X$ can be covered by a finite family 
$U_1,\dots , U_l$ closed under intersection of open definable subsets 
definably homeomorphic to open balls. We prove the result by induction on
$l$. The case $l=1$ is corollary 3.10 \cite{vdd}. 

Now let $X_1=U_1\cup \cdots \cup U_{l-1}$ and let $X_2=U_l$ and assume
that 
the result also holds for $X_1$. Let $A_i=B-f^{-1}(X_i)$ for $i=1,2.$
Then 
$A_1\cap A_2$ is empty, and $A_1,$ $A_2$ are closed. Since $Y$ is
affine, by 
lemma 3.5 \cite{vdd} it is normal i.e., we can separate $A_1$ and
$A_2$ by 
open definable sets $Y_1$ and $Y_2$.

Let $Y_0=Y-(Y_1\cup Y_2)$, closed in $Y$, hence also normal. Let 
$B_i=Y_i\cap B$, $i=0,1,2$. Then $f(B_i)\subseteq X_i$ for $i=1,2$ and 
$f(B_0)\subseteq X_1\cap X_2$. Since the result holds for $X_1\cap
X_2$, 
$f_{|B_0}$ can be extended to a continuous definable map $g_0$ on an open 
definable neighbourhood $U_0$ of $B_0$ in $Y_0$. Then $U_0\cap B=B_0$, so
$f$ 
together with $g_0$ define a continuous definable map $g:U_0\cup B\into X$.

Again by lemma 3.5 \cite{vdd}, there are disjoint (relative) open 
definable subsets $V,W\subseteq Y_0$ such that  $B_0\subseteq V$ and 
$Y_0-U_0\subseteq W$. Then $U'_0=Y_0-W$ is closed and $U'_0\subseteq
U_0$. 
Now for $i=1,2$, $g(U'_0\cup B_i)\subseteq X_i$ and $U'_0\cup B_i$ is
closed 
in $Y$. Since the result holds for $X_i$, $g_{|U'_0\cup B_i}$ extends
to a 
continuous definable map $G_i:U_i\into X_i$ on an open definable neighbourhood
$U_i$, 
$i=1,2.$ Now $U'_i=U_i\cap (U'_0\cup Y_i)$ is closed in $U:=U'_1\cup
U'_2$, 
$i=1,2$ and $U'_0=U'_1\cap U'_2$. Hence we can define a well defined 
continuous definable map $F:U\into X$ by  $F_{|U'_i}:=G_i$, $i=1,2.$
Moreover, $U$ contains the open definable neighbourhood $(U_1\cap (V\cup
Y_1))\cup 
(U_2\cap (Y\cup Y_2))$ of $B$.
\qed

\begin{prop}\label{prop same cohomology}
There is an isomorphism $\kappa :\mbox{{\it \u H}}^q(A;R)\into H^q(A;R)$.
\end{prop}

\pf
Passing the inclusion homomorphisms $H^q(V;R)\into H^q(A;R)$ to the
limit 
we get a canonical homomorphism $\kappa :\mbox{{\it \u H}}^q(A;R)\into H^q(A;R).$
By proposition 3.3 \cite{vdd} there is a definable retraction $r:V\into A$
of an open definable neighbourhood $V$ of $A$. Let $i:A\into V$ be the 
inclusion. Then 
$H^q(ri)=$ identity, and so $H^q(i)$ is an epimorphism, hence $\kappa
$ is also an epimorphism.

\medskip
{\it Claim:} There is an open definable neighbourhood $W$ of $A$ contained in
$V$ such 
that if $j:W\into V$ is the inclusion, then there is a definable homotopy
between 
$i\circ r_{|W}$ and $j$.

\medskip
{\it Proof of Claim:} On the closed subset $(V\times 0)\cup (A\times
[0,1])
\cup (V\times 1)$ of $V\times [0,1]$, set $F(x,t)=x$ if $(x,t)\in
V\times 0$, 
$F(x,t)=r(x)$ if $(x,t)\in V\times 1$ and $F(x,t)=x$ if $(x,t)\in A
\times [0,1].$ Since $V\times [0,1]$ is normal, by lemma \ref{lem anr}
$F$ 
extends to a definable map of an open definable neighbourhood of this set into
$V$. That 
open definable neighbourhood contains a definable set of the form 
$W\times [0,1]$,
where $W$ is an open definable neighbourhood of $A$. This gives the
desired definable homotopy.
\qed

By the Claim, $H^q(j)=H^q(r_{|W})H^q(i)$ and we have a factorisation
\[
\begin{array}{clcr}
H^q(U';R)\,\,\,\,\rightarrow \,\,\,\,\,
H^q(A;R)\,\,\,\,\,\,\,\,\,\,\,\,\,\,\,\,\,\,\,\,\,\,\,\,\,\,\,\,\,\,\,\,\,\,
\,\,\,\,\,\,\,\,\\
\downarrow
\,\,\,\,\,\,\,\,\,\,\,\,\,\,\,\,\,\,\,\,\,\,\,\,\,\,\,\,\,\,
\uparrow\,\,\,\,\,\,\,\,\,\,\,\,\,\,\,\,\,\,\,\,\,\,\,\,\,\,\,\,\,\,\,\,\,\,
\,\,\,\,\,\,\,\, \\
H^q(V;R)\into H^q(A;R)\into H^q(W;R)
\end{array}
\]
where $U'$ is any open definable neighbourhood of $A$ containing
$V$. It follows that 
any class in $H^q(U';R)$ going to zero in $H^q(A;R)$ goes to zero in 
$H^q(W;R)$ and thus $\kappa $ is a monomorphism.
\qed

\begin{cor}
There is an isomorphism $H^q_c(X-A;R)\into H^q(X,A;R)$.
\end{cor}

\pf 
The same as the proof of corollary 27.4 \cite{g}.
\qed

\begin{nrmk}\label{nrmk relative alexander duality}
{\em (For details see 27.6 \cite{g} and 27.7 \cite{g}).
Let $(K,L)$ be a pair of definably compact subsets of $X$. Then there is a 
relative Alexander duality $$\mbox{{\it \u H}}^q(K,L;R)\into H_{n-q}(X-L,X-K;R)$$
which is 
an isomorphism in the case $K=X$. Here, $\mbox{{\it \u H}}^q(K,L;R)$ is the direct
limit 
$\lim _{(U,V)}H^q(U,V;R)$ where $(U,V)$ runs through the directed set
of open definable sets containing $(K,L)$.

We can similarly define a canonical homomorphism $$H^q(X-L,X-K;R)\into
\mbox{{\it \u H}}_
{n-q}(K,L;R)$$
which is an isomorphism if $R$ is a field. Here, $\mbox{{\it \u H}}_{n-q}(K,L;R)$
is the 
inverse limit
$\lim _{(U,V)}H_{n-q}(U,V;R).$ 
}
\end{nrmk}

\begin{cor}\label{cor general separation thm}
(Separation Theorem).
If $A$ is a definably compact submanifold of $N^m$ of dimension
$m-1$ and having $k$ definably connected components, then the complement of
$A$ has $k+1$ definably connected components.
\end{cor}

\pf
By proposition \ref{prop same cohomology} we have
$\mbox{{\it \u H}}^q(A;R)=H^q(A;R)$. Regarding 
$N^m$ as $S^m-$point, we have the isomorphisms $H^q(A;R)\simeq H_{m-q}
(S^m,S^m-A;R)\leftarrow H_{m-q}(N^m,N^m-A;R)\simeq \tilde{H}_{m-q-1}(N^m-A;R).$
\qed

The standard proof of the fact that $\tilde{H}_q(S^n)=\ZZ$ iff $q=n$
and zero otherwise, remains valid in our case, but the usual proof of
the fact that if $e$ is a definable subset of $S^n$ that is definably 
homeomorphic to $[0,1]^r$ for some $r\leq n$, then
$\tilde{H}_q(S^n-e)=0$ for all $q$ depends on the compactness of
closed and bounded subsets of the reals. Woerhiede circumvents this
difficulty by using the definable trivialization theorem. As corollaries
of this fact and the Mayer-Vietoris sequence we get (see \cite{Wo}): $(1)$ 
{\it (The Jordan-Brouwer separation theorem)}
for $n>0$, $\tilde{H}_q(S^n-s_r)=\ZZ$ iff $q=n-r-1$ and zero otherwise
for every definable subspace $s_r$ of $S^n$ which is definably
homeomorphic to $S^r$, in particular $S^n-s_{n-1}$ has exactly two
definable connected components, and $s_{n-1}$ is their common
boundary; $(2)$ {\it (Invariance of domain)} if $U\subseteq N^n$ is a
definable open subset and $h:U\into N^n$ is an injective definable
continuous map, then $h(U)$ is open in $N^n$ and $h$ is a definable
homeomorphism onto $h(U)$.

\end{subsection}

\begin{subsection}{Lefschetz duality theorem}
\label{subsection lefschetz duality theorem}

\begin{lem}\label{lem fundamental class on the boundary}
Let $X$ be a definably compact definable manifold with boundary
and let $s\in \Gamma (\mbox{\. X};R)$ an $R$-orientation of
$\mbox{\. X}$. Then there is a unique homology class $\zeta \in 
H_n(X,\partial X;R)$ such that for any
$x\in \mbox{\. X}$, $s(x)=j_x^{\mbox{\. X}}(\zeta )$. Moreover, 
$\partial \zeta \in H_{n-1}(\partial X;R)$
is the fundamental class for the induced $R$-orientation of $\partial X$.
\end{lem}

\pf
The same as proposition 28.15\cite{g} and corollary 28.16\cite{g}.
\qed

\begin{thm}\label{thm lefschetz duality}
(Lefschetz Duality Theorem).
Suppose that $X$ is a definably compact definable manifold with
boundary, such that
let $\mbox{\. X}$ be $R$-orientable.
Let $\partial \zeta \in H_{n-1}(\partial X;R)$ be the fundamental
class. The 
the diagram
\[
\begin{array}{clcr}
\rightarrow H^{q-1}(X;R)\,\,\,\, {\rightarrow}\,\,\,\, H^{q-1}(\partial X;R)
\,\,\,\,\, \stackrel{\delta }{\rightarrow}\,\,\,\,\,\, H^q(X,\partial
X;R)\,\,\,\,\,\,{\rightarrow}
H^q(X;R)\,\,\,\,\,  \\
\,\,\,\downarrow ^{\zeta \cap}\,\,\,\,\,\,\,\,\,\,\,\,\,\,\,\,\,\,\,\,\,\,\,\,\,\,\,\,\,\,
\downarrow ^{\partial \zeta \cap}\,\,\,\,\,\,\,\,\,\,\,\,\,\,\,\,\,\,\,\,\,\,\,\,\,\,\,\,\,\,\,\,\,\,\,\,
\downarrow ^{\zeta \cap }\,\,\,\,\,\,\,\,\,\,\,\,\,\,\,\,\,\,\,\,\,\,\,\,\,\,\,\,\,\,\downarrow \\
\rightarrow H_{n-q+1}(X,\partial X;R)\stackrel{\partial }{\rightarrow}
H_{n-q}(\partial X;R){\rightarrow}H_{n-q}(X;R)\rightarrow H_{n-q}
(X,\partial X;R)
\end{array}
\]
is sign-commutative and the vertical arrows are isomorphisms.
\end{thm}

\pf
Similar to the proof of 28.18\cite{g}.
\qed

\end{subsection}

\end{section}

\begin{section}{Lefschetz fixed point theorem}
\label{section lefschetz fixed point theorem}

Through this section, ${\mathbf X}=(X,(X_i,\phi _i)_{i\in I})$ is an
$R$-orientable manifold of dimension $n$ with $R$-orientation $s\in 
\Gamma (X,R)$.

\begin{subsection}{Thom class}
\label{subsection thom class}

Consider the dual sheaf $U\into \Gamma ^*(U,R)$ of
the $R$-orientation sheaf $\tilde{O}^X$. $s$ determines a global
section $s^*\in \Gamma ^*(X,R)$ such that $(s(x),s^*(x))=1$
for all $x\in X$.

Let $\Delta _X$ be the diagonal of $X\times X$ and 
$i^U_x:(X,X-x)\into (X\times U, X\times U-\Delta _X)$ be the map given
by $i^U_x(z):=(z,x)$ for $z\in X$.  

\begin{thm}\label{thm thom isomorphism}
(Thom isomorphism theorem).
Let $U\subseteq X$ be an open definable subset. Then for all $q<n$, 
$H^q(X\times U,X\times U-\Delta _X;R)=0$ and there is a unique
isomorphism $\phi :H^n(X\times U,X\times U-\Delta _X;R)\into \Gamma
^*(U,R)$ such that $\phi (\beta )(x)=H^n(i^U_x)(\beta )$
for all $\beta \in H^n(X\times U,X\times U-\Delta _X;R)$, $x\in U$. 
\end{thm}

\pf
As we saw before, by lemma \ref{lem good triangulation}, 
$X$ is the union of a finite family $U_1,\dots ,U_l$ closed under
intersection of open definable subsets
definably homeomorphic to open balls. 

For any $\beta \in H^n(X\times U,X\times U-\Delta _X;R)$ define  
a set-theoretic section $\Phi (\beta ):U\into (\tilde{O}^X)^*$ by
$\Phi (\beta )(x)=H^n(i^U_x)(\beta )$ for all $x\in U,$ where 
$(\tilde{O}^X)^*$ is the \'etal\'e space dual to $\tilde{O}^X$ i.e.,
whose 
fibre at $x$ is the local cohomology $R$-module $H^n(X,X-x;R).$

If $\Gamma '(U;R)$ denotes the $R$-module of set-theoretic sections
$U\into 
(\tilde{O}^X)^*$, then for $V\subseteq U$ we have the commutative diagram
\[
\begin{array}{clcr}
H^n(X\times U, X\times U-\Delta _X;R)\into H^n(X\times V,X\times V-
\Delta _X;R)\\
\,\,\,\,\downarrow ^{\Phi }\,\,\,\,\,\,\,\,\,\,\,\,\,\,\,\,\,\,\,\,\,\,\,\,\,
\,\,\,\,\,\,\,\,\,\,\,\,\,\,\,\,\,\,\,\,\,\,\,\downarrow ^{\Phi }\\
\,\,\,\,\,\Gamma '(U;R)\,\,\,\,\,\,\,\into \,\,\,\,\,\,\,\,\,\,\Gamma '(V;R)
\end{array}
\]
 
Thus to verify that the homomorphism $\Phi $ takes its values in
$\Gamma ^*
(U;R)$, it suffices to consider the following cases:

\medskip
{\it Case (1):} $U$ is definably homeomorphic to an open ball and is
contained in 
an open definable set $V$ which is definably homeomorphic to an open ball.

\medskip
{\it Proof of Case (1):} For each $x\in U$ we have a commutative diagram
\[
\begin{array}{clcr}
H^q(X\times U,X\times U-\Delta _X;R)\into H^q(V\times U, V\times 
U-\Delta _V;R)\\ 
\,\,\,\,\downarrow ^{i^U_x}\,\,\,\,\,\,\,\,\,\,\,\,\,\,\,\,\,\,\,\,\,\,\,\,\,
\,\,\,\,\,\,\,\,\,\,\,\,\,\,\,\,\,\,\,\,\,\,\downarrow ^{i^U_x}\\
H^q(X,X-x;R)\,\,\,\,\into \,\,\,\,\,H^q(V,V-x;R)
\end{array}
\]
where the horizontal isomorphisms are excisions. Thus we may assume
that 
$X=V=N^n.$ In this case we have a homeomorphism
$$ f:(N^n\times U, (N^n-0)\times U)\into (N^n\times U, N^n\times U-
\Delta _{N^n})$$
given by $f(y,x):=(y+x,x)$ and for each $x\in U$, a commutative diagram
\[
\begin{array}{clcr}
(N^n\times 0, N^n_0\times 0)\stackrel{1_{N^n}\times j_x}{\rightarrow
} (N^n\times U, N^n_0\times U)\stackrel{f}{\rightarrow }(N^n\times
U, N^n\times U-\Delta _{N^n})\\
\,\,\,\,\,\,\,\,\,\,\,\,\,\,\,\,\,\,\,\,\,\,\,\,\,\,\,\,\,\,\,\,\,\,\,\,\,\,
\,\,\,\,\,\,\,\,\,\,\,\,\,\,\,\,\,\,\,\,\,\,\,\,\,\,\,\,\,\,\,\,\,\,\,\,\,\,
\downarrow \,\,\,\,\,\,\,\,\,\,\,\,\,\,\,\,\,\,\,\,\,\,\,\,\,\,\,\,\, 
\,\,\,\,\,\,\,\,\,\,\,\,\,\,\,\,\,\uparrow ^{i^U_x}\\
(N^n\times 0, N^n_0\times 0)\,\,\,\,\,\,\,\,\,\,\,\,
{\leftarrow }\,\,\,\,\,\,\,\,\,\,\,\,
(N^n,N^n_0)\,\,\,\,\,\,\,\,\,\stackrel{f_x}{\rightarrow
}\,\,\,\,\,\,\,\,\,\,\,(N^n,N^n_x)\,\,\,\,\,\,\,\,
\end{array}
\]
where $N^n_z=N^n-z$, $f_x(y):=x+y$, and $j_x$ is the map of the point $0$ onto
$x$. We may 
assume that $0\in U$.

\medskip
{\it Claim (1.1):} The map $s'\into s'(0)$ is an isomorphism of 
$\Gamma ^*(U;R)$ onto $H_n(N^n,N^n-0;R)$.

\medskip
{\it Proof of Claim (1.1):} This follows from the fact that
$j^U_0$ is 
an isomorphism (by lemma \ref{lem locally constant}).
\qed

\medskip
{\it Claim (1.2):} If $s'\in \Gamma ^*(U;R)$ and $x\in U$ then $s'(x)=
H^n(f_x)(s'(0)).$

\medskip
{\it Proof of Claim (1.2):} Let $\alpha \in H^n(N^n,N^n-U;R)$ be the
unique 
class such that $\alpha =j^U_x(s'(x))$ for all $x\in U$. Now the maps 
$l^U_x,f_xl^U_0:(N^n,N^n-U)\into (N^n,N^n-x)$ (where $l^U_x$ is the
natural inclusion, and $H^n(l^U_x)=j^U_x$) are definably homotopic (at
time $t$ 
the map is $f_{tx}l^U_0$), whence
$s'(x)=H^n(l^U_x)^{-1}(\alpha )=H^n(f_x)^{-1}H^n(l^U_0)^{-1}H^n(l^U_0)(s'(0))$$
=H^n(f_{-x})(s'(0)).$
\qed

\medskip
{\it Claim (1.3):} $H^q(i^U_0)$ is an isomorphism for all $q$.

\medskip
{\it Proof of Claim (1.3):} Using the commutative diagram above, we
must show 
that
$$H^q(1_{N^n}\times j_0):H^q(N^n\times U,(N^n-0)\times U;R)\into 
H^q(N^n\times 0,(N^n-0)\times 0);R)$$  is an isomorphism, which follows
from 
the fact that $U$ is definably contractible.
\qed

\medskip
Now the theorem for Case (1) follow from Claims (1.1) and (1.3),
since by 
Claim (1.2) and the fact that $i^U_x$ is definably homotopic to $i^U_0\circ
f_{-x}$
we have $\Phi (\beta )(x):=H^n(i^U_x)(\beta )=H^n(f_{-x})H^n(i^U_0)(\beta ).$  
\qed

\medskip
{\it Case (2):} If the theorem holds for open definable subsets $U, V$
and $W=U\cap V$, then it holds for the open $Y=U\cup V$.

\medskip
{\it Proof of Case (2):} Let $U'=X\times U-\Delta _X$,
$V'=X\times V-\Delta _X$ , $W'=X\times W-\Delta _X$ and $Y'=X\times Y-
\Delta _X$.

\medskip
{\it Claim (2.1):} There is an exact sequence
\[
\begin{array}{clcr}
\into H^q(X\times Y, Y';R)\stackrel{i}{\rightarrow }H^q(X\times U, U';R)\oplus H^q(X\times V,V';R)\\
\,\,\,\,\,\,\,\,\,\,\,\,\,\,\stackrel{j}{\rightarrow }H^q(X\times
W,W';R)\stackrel{k}{\rightarrow }H^{q+1}(X\times Y,Y';R)\into
\end{array}
\]
where $i$ is induced by the chain homomorphism $z\into (z,z)$, $j$ by
the 
chain homomorphism  $(z,w)\into z-w$, and $k$ is the connecting homomorphism.

\medskip
{\it Proof of Claim (2.1):} By the Universal Coefficient Theorem, its
enough 
to prove Claim (2.1) for $R=\ZZ$. Consider the monomorphism of chain complexes
$$ i:S_*(X\times W)/S_*(W')\into S_*
(X\times U)/S_*(U')\oplus S_*(X\times V)/S_*(V')$$
given by $i(\bar{z})=(\bar{z},\bar{z})$, and the chain epimorphism
\[
\begin{array}{clcr}
j:S_*(X\times U)/S_*(U')\oplus S_*(X\times V)/
S_*(V')\into \\
\into (S_*(X\times U)+S_*(X\times
V)/(S_*(U')+S_*(V'))
\end{array}
\]
given by $j(\bar{z},\bar{w})=\bar{z}-\bar{w}$. Clearly $ji=0$. Suppose
that 
$j(\bar{z},\bar{w})=0$. If $z=\sum \mu_i\sigma _i$, let $v$ be the sum
of 
those $\mu _i\sigma _i$ such that $|\sigma _i|$ meets $\Delta _X$;
then $v$ 
is equal to the chain defined in the same way using $w$ instead of
$z$, since 
$z-w\in S_*(U')+S_*(V')$, so that $i(\bar{v})
=(\bar{z},\bar{w})$. Therefore, we have an exact sequence of chain
complexes. 
Since these complexes are free, dualizing gives an exact sequence of
cochain 
complexes, hence an infinite exact cohomology sequence
\[
\begin{array}{clcr}
\into H^q(C/C')\stackrel{i}{\rightarrow }H^q(X\times U, U')\oplus 
H^q(X\times V,V')\\
\,\,\,\,\,\,\,\,\,\,\,\,\,\,\stackrel{j}{\rightarrow }H^q(X\times
W,W')
\stackrel{k}{\rightarrow }H^{q+1}(C/C')\into
\end{array}
\]
where $C=S_*(X\times U)+S_*(X\times V)$ and $C'=S_*(U')+
S_*(V').$ Since $X\times Y=X\times U\cup X\times V$, the
inclusion 
$C\into S_*(X\times Y)$ is a chain homotopy (see the proof of
the 
excision theorem in \cite{Wo}). This holds similarly for the inclusion 
$C'\into S_*(Y')$, and by passage to the quotient, for the
inclusion 
$C/C'\into S_*(X\times Y)/S_*(Y')$. Thus we can replace 
$H^q(C/C')$ by $H^q(X\times Y,Y')$ in the above exact sequence.
\qed

We now prove Case (2). By the assumption we see that $$H^q(X\times
Y,Y';R)=0$$ 
for $q<n$ and the commutative diagram
\[
\begin{array}{clcr}
H^n(X\times Y,Y';R)\stackrel{i}{\rightarrow }H^n(X\times U, U')\oplus H^n(X\times V,V')\stackrel{j}{\rightarrow }H^n(X\times W, W')\\
\,\,\,\,\downarrow ^{\Phi }\,\,\,\,\,\,\,\,\,\,\,\,\,\,\,\,\,\,\,\,\,\,\,\,\,
\,\,\,\,\,\,\,\,\,\,\,\,\,\,\,\,\,\,\,\,\,\,\,\,\,\,\,\,\,\,\,\,\,\,\,\,\,
\downarrow ^{\Phi \oplus \Phi }\,\,\,\,\,\,\,\,\,\,\,\,\,\,\,\,\,\,\,\,\,\,\,\,\,\,\,\,\,\,\,\,\,\,\,\,\,\,\,\,\,\,\,\,\,\,\,\,\,\,\downarrow ^{\Phi }\\
0\,\,\,\into \,\,\Gamma ^*(Y;R)\,\,\,\,\,\,\,\,\,\,\,\,
{\stackrel{i}\rightarrow }\,\,\,\,\,\,\,\,\,\,\,\,\,\,\Gamma
^*(U;R)\oplus \Gamma ^*(V;R)\,\,\,\,\stackrel{j}{\rightarrow }
\,\,\,\,\,\Gamma ^*(W;R)\,\,\,\,\,\,\,\,\,\,\,
\end{array}
\]
together with the $5$-lemma imply that $\Phi $ is an isomorphism for $Y$.
\qed
\qed

\begin{defn}\label{defn thom class}
{\em
This means that there is a unique cohomology class $\mu _X$ in 
$H^n(X\times X,X\times X-\Delta _X;R)$ such that $s^*(x)=H^n(i^X_x)(\mu
_X )$. $\mu _X$ is called the {\it Thom class} of the given
$R$-orientation. The {\it Lefschetz class of $X$} 
is the image $\Lambda _X=H^n(j)(\mu _X)$ 
of the Thom class $\mu _X$ under the homomorphism
$H^n(j):H^n(X\times X, X\times X-\Delta _X;R)\into H^n(X\times X;R)$ induced
by the inclusion $j:X\times X\into (X\times X, X\times X-\Delta _X)$. 
}
\end{defn}

\end{subsection}
\begin{subsection}{The Lefschtez isomorphism}
\label{subsection lefschtez isomorphism}

Before we proceed we need the following lemma

\begin{lem}\label{lem commutativity}
Suppose that $X$ is definably compact. 
If $\tau \in H^p(X\times X, X\times X-\Delta _X;R)$ and $\sigma \in
H^q(X;R)$,
then $$H^p(j)(\tau )\cup (\sigma \times 1)=H^p(j)(\tau )(1\times
\sigma ).$$ 
\end{lem}

\pf
\medskip
{\it Claim:} There is an open definable neighbourhood $V$ of $\Delta _X$ in
$X\times X$ and a definable retraction $r:V\into \Delta _X$ such that
$i\circ r$ is definably
homotopic to $k$, where $i:\Delta _X\into X\times X$ and $k:V\into
X\times X$ 
are the inclusions.

\medskip
{\it Proof of Claim:}
Suppose that $X\subseteq N^m$. Then by proposition 3.3 \cite{vdd}
there is an 
open definable neighbourhood $U$ of $X$ having a definable retraction $s:U\into
X$. Let 
$\epsilon =$ distance from $X$ to $N^m-U$, and let $V$ be the $\epsilon
$-neighbourhood of $\Delta _X$ in $X\times X$. Define $F:X\times
X\times [0,1]
\into N^m$ by $F(x,y,t):=(1-t)x+ty.$ Then $F$ maps $V\times [0,1]$
into $U$. 
Let $G:=r\circ F_{|V\times [0,1]}:V\times [0,1]\into X$ so that
$G(x,y,0)=s(x)
=x,$ $G(x,y,1)=s(y)=y.$ The required definable homotopy $H:V\times
[0,1]\into X
\times X$ is defined by $H(x,y,t):=(x,G(x,y,t)).$
\qed

\medskip
Let $k':(V,V-\Delta _X)\into (X\times X, X\times X-\Delta _X)$ denote
the 
inclusion. Note that $k'$ is an excision. We have
\[
\begin{array}{clcr}
H^q(\Delta _X;R)\,\,\,\,\,\,\,\stackrel{i}{\leftarrow }\,\,\,\,\,\,\,
H^q(X\times X;R)\stackrel{\lambda }{\rightarrow}
H^q(X\times X,X\times X-\Delta _X;R)\\
\downarrow \,\,\,\,\,\,\,\,\,\,\,\,\,
\,\,\,\,\,\,\,\,\,\,\,\,\,\,\,\,\,\,\,\,\,\,\,\,\,\,\,\,\,\,\,\,\,\,\,\,\,
\downarrow ^{H^q(k)}\,\,\,\,\,\,\,\,
\,\,\,\,\,\,\,\,\,\,\,\,\,\,\,\,\,\,\,\,\,\,\,\,\,\,\,\,\,\,\,\,\,\,\,\,\,
\,\,\,\,\,\,\,\,\downarrow ^{H^{p+q}(k')}\\
H^q(\Delta _X;R)\,\,\,\,\,\,\,\,\,\,\,\,\stackrel{r}{\rightarrow }\,\,\,\,\,\,\,\,\,\,\,\,H^q(V;R)\,\,\,\,\,\,\,\,\,\,\,\,\stackrel{\rho }
{\rightarrow }\,\,\,\,\,\,\,\,\,\,\,\,H^q(V,V-\Delta _X;R)\,\,\,\,\,
\end{array}
\]
where we are using mixed cup products of absolute and relative
cohomology, 
$\lambda (\beta ):=\tau \cup \beta $ and $\rho (\alpha ):=H^p(k')(\tau
)\cup 
\alpha .$

Let $p_i:X\times X\into X$, $i=1,2$ be the projections. Then $1\times
\sigma 
=H^0(p_1)(1)\cup H^p(p_2)(\sigma )$$=H^q(p_2)(\sigma ),$ $\sigma
\times 1=
H^q(p_1)(\sigma )$. Let $p:\Delta _X\into X$ be the common restriction
of 
$p_1$ and $p_2$ to the diagonal. From the diagram we get 
$$\tau \cup H^q(p_i)(\sigma )=H^{p+q}(k')^{-1}(H^p(k')(\tau )\cup
H^q(p_iir)
(\sigma ))$$
$$=H^{p+q}(k')^{-1}(H^p(k')(\tau )\cup H^q(pr)(\sigma ))$$
for both $i=1$ and $2$. By definition of the mixed cup product,
$$H^p(j)(\tau )\cup H^q(p_i)(\sigma )=H^{p+q}(k')(\tau \cup
H^q(p_i)(\sigma ))$$
which proves the lemma.
\qed

\begin{lem}\label{lem explicit poincare duality}
Suppose that $X$ is  definably connected and definably compact. 
Then for any $p\leq n$, the inverse to the 
Poincar\'{e} duality isomorphism 
$$ D_X:H^p(X;R)\into H_{n-p}(X;R),\,\,\, \sigma \into \zeta _X\cap \sigma $$ 
is given by
$$D_X^{-1}:
H_{n-p}(X;R)\into H^p(X;R),\,\,\, \alpha \into (-1)^p\Lambda _X/\alpha .$$
\end{lem}

\pf
We first show that $\Lambda _X/\zeta _X=1$. For $x\in X$, consider the
commutative diagram
\[
\begin{array}{clcr}
(X,X-x)\,\,\,\stackrel{i_x}{\rightarrow}\,\,\, (X\times X, X\times
X-\Delta _X)\\
\uparrow
^{j_x}\,\,\,\,\,\,\,\,\,\,\,\,\,\,\,\,\,\,\,\,\,\,\,\,\,\,\,\,\,\,\,\,\,\,\,\,\,\,\,\,\,\,\,\,\,\,\,\uparrow ^j\\
X\,\,\,\,\,\,\,\,\,\,\stackrel{i_x}{\rightarrow} \,\,\,\,\,\,\,\,\,\,\,\,\,\,\,\,\,\,\,\,\, X\times X
\end{array}
\]
where $i_x=i^X_x$ and the vertical arrows are inclusions. Note that,
if 
$\bar{x}$ is the homology class of $x$
($x$ is a $0$-cycle), then $H_n(i_x)\zeta _X=\zeta _X\times \bar{x}$ 
(since $X\simeq X\times x$).
We have, $1=(s(x),H^n(i_x)(\mu _X))=
(H_n(j_x)\zeta _X,H^n(i_x)\mu _X)$$=(\zeta _X,H^n(i_xj_x)\mu _X)=
(\zeta _X,H^n(ji_x)\mu _X)$$
=(H_n(i_x)\zeta _X,\Lambda _X)=(\zeta _X\times \bar{x},\Lambda _X)
=(\bar{x},\Lambda _X/\zeta _X)$ (by theorem \ref{thm slant product}).

Consider $\sigma \in H^p(X;R)$, then we have 
$\Lambda _X/\zeta _X\cap \sigma =$$
1\cup (\Lambda _X/\zeta _X\cap \sigma )$
$=(-1)^{p(n+p+0-n)}[(\sigma \times 1)\cup \Lambda _X]/\zeta _X$ 
(by theorem \ref{thm slant product}) $=(-1)^{p^2}
[(1\times \sigma )\cup \Lambda _X]/ \zeta _X$ 
(by lemma \ref{lem commutativity}) 
$=(-1)^{p^2}(-1)^0\sigma \cup [\Lambda _X/\zeta _X\cap 1]$
(by theorem \ref{thm slant product}) 
$=(-1)^p\sigma \cup (\Lambda _X/\zeta _X)$ 
$=(-1)^p\sigma \cup 1=(-1)^p\sigma $.
\qed

Let $f:X\into Y$ be a continuous  definable map, where $Y$ is another  
definably compact, $R$-oriented definable manifold of
dimension $m$. 
We define the cohomology class 
$\mu _f$ of the graph of $f$ by 
$$\mu _f=H^m(f\times 1_Y)(\Lambda _Y)\in H^m(X\times Y;R).$$
If $X=Y$, the {\it Lefschetz class $L_f$ of $f$} is defined by
$L_f:=H^n(\Delta _X)(\mu _f)\in H^n(X;R)$. 

Note that, for $\sigma \in H^p(Y)$ we have $\mu _f/\zeta _Y\cap
\sigma $$=H^p(f)(\mu '_Y/\zeta _Y\cap \sigma )$ (by theorem \ref{thm
slant product} $(3)$) $=(-1)^pH^p(f)(\sigma )$ (by lemma \ref
{lem explicit poincare duality}).

\begin{defn}\label{defn lefschetz isomorphism}
{\em
Let $L^p(X;R):=Hom_R(H^p(X;R),H^p(X;R))$ and let
$$L^*(X;R):=\sum _{p=0}^nL^p(X;R).$$ For each $p$ we have a canonical 
isomorphism 
$$k^p:H^p(X;R)\otimes _RH_p(X;R)\into L^p(X;R)$$ which
induces a 
canonical isomorphism
$$k:\sum _{p=0}^nH^p(X;R)\otimes _RH_p(X;R)\into L^*(X;R)$$ given by 
$k:=\sum _{p=0}^n(-1)^pk^p$.
The {\it Lefschetz isomorphism for $X$} is the isomorphism of $R$-modules 
$$ \lambda _X:L^*(X;R)\into H^n(X\times X;R)$$ given by $\lambda _X:=
\alpha '\circ (1^*_X\otimes _RD^{-1}_X)\circ
k^{-1}$ where, $\alpha '$ is the K\"unneth isomorphism and $D^{-1}_X$
is the 
inverse of the Poincar\'e duality isomorphism.
Note that, $\Lambda _X=\lambda _X(1^*_X)$.
}
\end{defn}

\begin{lem}\label{lem trace and lefschetz isomorphism}
Let $Tr:L^*(X;R)\into R$ be the linear map given by $Tr \sigma :=\sum
_{p=0}^n
(-1)^ptr \sigma ^p$ where $\sigma =
\sum _{p=0}^n \sigma ^p$, $\sigma ^p\in L^p(X;R)$. Then 
$$Tr\sigma =(\zeta _X, \Delta _X^*\lambda _X(\sigma )).$$
\end{lem}

\pf
Its enough to consider $\sigma =k(\beta \otimes _RD_X\gamma )$ with
$\beta 
\in H^p(X;R)$ and $\gamma \in H^{n-p}
(X;R)$. Then, by ordinary linear algebra $Tr\sigma
=(-1)^{np+p}(D_X\gamma ,\beta)$
$=(-1)^{p(n-p)}(D_X\gamma ,\beta )
$$=(-1)^{p(n-p)}(\zeta _X\cap \gamma , \beta )=(-1)^{p(n-p)}(\zeta _X, 
\gamma \cup \beta )$
$=(\zeta _X,\beta \cup \gamma )$$=(\zeta _X, \Delta _X^*\alpha '(\beta 
\otimes _R\gamma ))$$=
(\zeta _X, \Delta _X^*\lambda _X(\sigma )).$
\qed

\begin{defn}\label{defn poincare adjoint}
{\em
The {\it Poincar\'e adjoint of $f^*$} where $f:X\into Y$ is a
continuous definable map,
is the 
unique linear map
$\tilde{f}:=\sum _{p=0}^n\tilde{f}^p$ where
$\tilde{f}^{n-p}:H^{n-p}(X;R)
\into H^{m-p}(Y;R)$ is determined by
$(D_Y\tilde{f}^{n-p}\alpha ,\beta )=(D_X\alpha ,f^{*p}\beta )$ for
all 
$\alpha \in H^{n-p}(X;R)$ and $\beta \in
H^p(Y;R)$.
}
\end{defn}

Its easy to see that the Poincar\'e adjoint of a composition of
continuous definable maps  is the composition of the Poincar\'e
adjoints, and if 
$dim X=dimY$ then $\tilde{f}\circ f^*=
(degf)1_X^*$. In particular, if $X=Y$ then $f^*$ is a linear
isomorphism iff 
$degf\neq 0$, in which case $(f^*)^{-1}=
\frac{1}{degf}\tilde{f}$.

\begin{lem}\label{lem lambda and functions}
Let $\sigma \in L^*(Y;R)$, $f,g:X\into Y$ be continuous definable maps
and suppose
that 
$dimX=dimY$. Then
$$(f\times g)^*(\lambda _Y(\sigma ))=\lambda _X(f^*\circ \sigma \circ 
\tilde{g}).$$
In particular, if $X=Y$ then $\mu _f=\lambda _X(f^*)$.
\end{lem}

\pf
Its enough to take $\sigma =k(\alpha \otimes _RD_Y\beta )$ with
$\alpha 
\in H^p(Y;R)$ and $\beta \in H^{n-p}
(Y;R)$. We have $(f^*\circ \sigma \circ 
\tilde{g})(\gamma )=$$(-1)^{np}
(D_Y\beta ,\tilde{g}(\gamma ))f^*(\alpha )$$=(-1)^{np}(D_Xg^*(\beta)
,\gamma )
f^*(\alpha)$$=
[k(f^*\alpha \otimes _RD_Xg^*\beta)](\gamma )$
for all $\gamma \in H^p(X;R)$ . Therefore, $\lambda _X
(f^*\circ \sigma \circ \tilde{g})=
(f\times g)^*\circ \alpha '(\alpha \otimes _R\beta)=(f\times
g)^*(\lambda _Y
(\sigma )).$
\qed

\end{subsection}
\begin{subsection}{The Lefschetz fixed point theorem}
\label{subsection the lefschetz fixed point theorem}

\begin{defn}\label{defn coincidence number}
{\em
Let $f,g:X\into Y$ be continuous definable maps and suppose that
$dimX=dimY$. The 
{\it coincidence number of $f$ and $g$} is defined by
$$\lambda (f,g;R):=\sum _{p=0}^n(-1)^ptr(f^{*p}\circ \tilde{g}^p).$$
Note that if $X=Y$ then $\lambda (f,1_X;R)=\lambda (f;R)$.
}
\end{defn}

We have the following (see volume I chapter X \cite{ghv}): 
$\lambda (f,g;R)=(-1)^n\lambda (g,f;R)$ and if
$h:Z\into X$ is a third continuous definable map from a definably connected, 
definably compact, $R$-orientable definable manifold
then $\lambda (f\circ h, g\circ h;R)=(degh)\lambda (f,g;R)$.

\begin{thm}\label{thm lefschetz}
Let $X$ and $Y$ be a $R$-orientable, definably compact 
definable manifolds of dimension $n$, where $R$ is a field. 
If $f,g:X\into Y$ are continuous definable maps then 
$$\lambda (f,g;R)=(\zeta _X,\Delta _X^*\circ (f\times g)^*(\Lambda _Y)).$$
If $\lambda (f,g;R)\neq 0$, then there is $x\in X$ such that $f(x)=g(x)$.
\end{thm}

\pf
We have $\lambda (f,g;R)=Tr(f^*\circ \tilde{g})$$=(\zeta _X,\Delta^*
_X(\lambda _X(f^*\circ \tilde{g})))$$=
(\zeta _X,\Delta _X^*\circ (f\times g)^*(\Lambda _Y)).$

If there is no $x\in X$ such that $f(x)=g(x)$, then we have a factorisation
\[
\begin{array}{clcr}
X\,\,\,\,\,\,\,\,\stackrel{f\times g}{\rightarrow}\,\,\,\,\,\,\,\,\,
Y\times Y\\
\downarrow
^{\Delta _{X}}\,\,\,\,\,\,\,\,\,\,\,\,\,\,\,\,\,\,\,\,\,\,\,\,\,\,\,\,\,\,\uparrow ^i\\
X\times X\,\,\,\,\,\,\stackrel{f\times g}{\rightarrow}
\,\,\,\,\,\,\,\,\, Y\times Y-\Delta _Y
\end{array}
\]
where $i$ is the inclusion. 
Since $H^n(i)H^n(j)=0$ and $\Lambda _Y=H^n(j)(\mu _Y)$, we have 
$0=
\Delta _X^*\circ (f\times g)^* \circ i^*(\Lambda _Y)$
$=(f\times g)^*(\Lambda _Y)$ and therefore $\lambda (f,g;R)=0$.
\qed

\begin{cor}
Let $X$ be a definably connected, definably compact definable
manifold. If $X$ admits a continuous definable map $f:X\into X$ definably 
homotopic to the identity and without fixed points, then $E(X)=0$. 
\end{cor}

Note that, like in the classical case, all the results of this
section generalise to definable manifolds 
with boundary (see remark in 30.14 \cite{g}). 
\end{subsection}
\end{section}

\begin{section}{Cohomology rings of definable groups}
\label{section cohomology rings of definable groups}

Below $G$ will be  an definably connected, definably compact definable group of
dimension $n$. We will use $^{-1}$ and $\iota $ for inverse in $G$ and 
$\cdot $ and $m$ for multiplication in $G$. $X$ is like in the
previous section and $R$ will be a field of characteristic zero.

\begin{lem}\label{lem coincidence in G}
If $f,g:X\into G$ are two continuous definable maps then  
$$\lambda (f,g;R)=deg(f^{-1}\cdot g).$$
\end{lem}

\pf
Consider the definable 
map $q:G\times G\into G$ given by $q:=m\circ (\iota \times
1_G)$. A simple calculation using
co-multiplication in $H^*(G;R)$, the fact that $\iota ^*(\alpha )=
(-1)^p
\alpha $ for $\alpha \in H^p(G;R)$
and the 
definition of the Lefschtez isomorphism $\lambda _G$
shows that $q^*\omega _G=(-1)^n\omega _G\times 1+1\times \omega _G$$=
\lambda _G(1_G)=\Lambda _G$. The result follows from
the fact that $f^{-1}\cdot g=q\circ (f\times g)\circ \Delta _X$.
\qed

Before we prove our main theorem, computing the cohomology rings of 
definably compact, definably connected definable groups we recall some
examples of $\RR$-semialgebraic groups and the notion of Lie algebra
cohomology.

\begin{nrmk}\label{nrmk main}
{\em 
(see \cite{mt} and \cite{ghv}).
For $F=\RR$ , $\CC$ or $\mathbf{H}$ the division ring of quaternions over
$\RR$, let 
$$GL(n,F):=\{A\in M(n,F):\exists A^{-1}\in M(n,F),\, AA^{-1}=I_n\}$$ 
(the general linear group over $F$), 
$U(n,F):=\{A\in GL(n,F):AA^*=I_n\}$ (the unitary group over $F$).
For $F=\RR$ or $\CC$ let $SL(n,F):=\{A\in GL(n,F): detA=1\}$
(the special linear group over $F$),
$O(n,F)
:=\{A\in GL(n,F):A^tA=I_n\}$,
(the $F$ orthogonal groups)
$Sp(n,F):=\{A\in M(2n,F):^tAJ_nA=J_n\}$ 
(the $F$ symplectic groups) where $J_n:=(J^n_{i,j})$ with 
$J^n_{1,1}=J^n_{2,2}=0_n$, 
$J^n_{1,2}=-I_n$ and $J^n_{2,1}=I_n$.
   
The orthogonal groups are $O(n):=O(n,\RR$$)=U(n,\RR$$)$, the 
special orthogonal groups are
$SO(n):=U(n,\RR$$)\cap SL(n,\RR$$)$, the unitary groups are $U(n)
:=U(n,\CC$$)$, the special unitary groups 
are defined as
$SU(n):=U(n,\CC$$)\cap SL(n,\CC$$)$ and the 
symplectic groups are $Sp(n):=U(n,\mathbf{H}$$).$
We have a definable extension 
$1\into SO(n)\into O(n)\into \ZZ$$_2\into 1$ and the other groups of
this list are $\RR$-compact and $\RR$-connected
except for $GL(n,F):=U(n,F)\times \RR$$^{dn(n-1)/2+n}$ (for $F=$$\RR$ , 
$\CC$ or $\mathbf{H}$ and $d=1,2$ or $4$ respectively),
$SL(n,F):=SU(n,F)\times \RR$$^{dn(n-1)/2+n}$ (for $F=$$\RR$ or
$\CC$ and $d=1$ or $2$ respectively), $Sp(n,\RR$$)=U(n)\times
\RR$$^{n(n+1)}$, $Sp(n,\CC$$)=Sp(n)\times
\RR$$^{n(2n+1)}$ and $O(n,\CC$$)=O(n)\times \RR$$^{n(n-1)/2}$. 
}
\end{nrmk}

\begin{nrmk}\label{nrmk lie algebra cohomology}
{\em
Let $\mathbf{g}$ be a Lie algebra over $\RR$. For each $k\geq 0$ let
$C^k(\mathbf{g}$$, \RR$$):=L(\wedge ^k\mathbf{g}$$,\RR$$)$ and define
a differential by the Cartan formula
\[
\begin{array}{clcr}
dc(g_1,\dots ,g_{k+1})=\\
=\sum_{1\leq j<l\leq k+1}(-1)^{j+l-1}c([g_j,g_l],g_1,\dots ,
\hat{g}_j,\dots ,\hat{g}_l,\dots ,g_{k+1})\\
+\sum_{j=1}^{k+1}(-1)^jg_jc(g_1,\dots ,\hat{g}_j,\dots ,g_{k+1}).
\end{array}
\]
The cohomology of this complex is denoted by $H^*(\mathbf{g}$$,\RR$$)$
and is called the cohomology of the Lie algebra. If $G$ is a Lie group
with Lie algebra $\mathbf{g}$, and 
$Ad:G\into GL(\mathbf{g}$$)$ the adjoint representation,
$\mathbf{g}$$_I$ denotes the Lie subalgebra of $\mathbf{g}$ invariant
under $Ad$. If $G$ is a connected, compact Lie group, then
$H^*(G,\RR$$)=H^*(\mathbf{g}$$_I,\RR$$)$ (see \cite{ghv}).
}
\end{nrmk}

\begin{lem}\label{lem pi of an h space}
If $(X,e,m )$ is a definable $H$-manifold, then $\pi _1(X,e)$ is an
abelian group. In particular, there are 
non negative integers $s, m_1^{l_1}, \dots , m_u^{l_u}$ 
such that $\pi _1$$(X)=\ZZ$$^s\oplus \ZZ$$_{m_1^{l_1}}\oplus
\cdots \oplus \ZZ$$_{m_u^{l_u}}$.
\end{lem}

\pf
For definable paths $\alpha $ and $\beta $ in $X$ let $\alpha \beta $ 
be defined by $\alpha \beta (t):=m(\alpha (t),\beta (t)).$ It follows 
from the definition of definable $H$-manifold that $[\epsilon
_{e}\alpha]=[\alpha \epsilon _{e}]=[\alpha ]$; if $[\alpha ]=[\alpha
']$ and $[\beta ]=[\beta ']$ then $[\alpha \beta ]=[\alpha '\beta
']$. Further, we have the equality $(\alpha \cdot \beta )(\alpha
'\cdot \beta ')=(\alpha \alpha ')\cdot (\beta \beta ').$ Now the lemma 
follows from the definable homotopies:
$[\alpha \cdot \beta ]=[(\alpha \epsilon _{e})\cdot (\epsilon
_{e}\beta )]=$$[(\alpha \cdot \epsilon _{e})(\epsilon _{e}\cdot \beta
)]$$=[\alpha \beta ],$
$[\beta \cdot \alpha ]=[(\epsilon _{e}\beta )\cdot (\alpha \epsilon _{e})]=$
$[(\epsilon _{e}\cdot \alpha )(\beta \cdot \epsilon _{e})]$$=[\alpha \beta ].$ 
\qed

The next lemma follows from results from \cite{e2} on cover of
definable groups together with lemma \ref{lem pi of an h space}.

\begin{lem}\label{lem torsion points} 
If $G$ is a definably connected, definable abelian group then there are
$s, m_1^{l_1}, \dots , m_u^{l_u}\in \NN$$\cup \{0\}$ 
such that 
$\pi _1$$(G)=\ZZ$$^s\oplus \ZZ$$_{m_1^{l_1}}\oplus
\cdots \oplus \ZZ$$_{m_u^{l_u}}$ and
$card(\{x\in G:kx=0\})= k^s\cdot (k,m_1^{l_1})\cdots (k,m_u^{l_u}).$
\end{lem}

\pf
The definable map $p_k:(G,0)\into (G,0)$ is a surjective (since $G$ is 
divisible) homomorphism with finite kernel. Therefore, by \cite{e2} 
$p_k:(G,0)\into (G,0)$ is a definable covering map and since
$p_{k*}:\pi _1(G)\into
\pi _1(G)$ is given by $p_{k*}([\alpha ])=k[\alpha ]$, we have  that for all
$x,y\in G$, $|p_k^{-1}(x)|=|p_k^{-1}(y)|$$=[\pi _1(G,0):p_{k*}(\pi _1(H,0))]$$=
k^s\cdot (k,m_1^{l_1})\cdots (k,m_u^{l_u})$.
\qed

Finally we prove our main theorem.

\begin{thm}\label{thm main}
Let $G$ be an definably connected, definably compact definable group. 
Then the rank of $G$ equals the 
dimension of a maximal definably connected definably compact definable abelian 
subgroup, $\sum _{i=1}^{rank(G)}g_i=dimG$
and $dimG\equiv rank(G) (mod 2)$. We have,
$$H^*(G;R)=H^*(Z(G)^0;R)\otimes _RH^*(G/Z(G)^0;R),$$
$$\pi _1(G)=\pi _1(Z(G)^0)\oplus\pi _1(G/Z(G)^0).$$
Moreover, we have (1) 
$H^*(Z(G)^0;R)=\bigwedge [y_{1_1},\dots y_{1_{dimZ(G)^0}}]_R$, 
$\pi _1(Z(G)^0)=\ZZ $$^{dimZ(G)^0}$, 
and for each $k>1$, $card(\{x\in Z(G)^0:kx=0\})=k^{dimZ(G)^0}$, (2)
$\pi _1(G/Z(G)^0)$ is finite and
$H^*(G/Z(G)^0;R)$ is a obtained by taking tensor products over $R$ of
the following types of free, skew-commutative graded Hopf
$R$-algebras:
\[
\begin{array}{clcr}
Type \,\, A_l\,(l\geq 1) :\,\,\bigwedge [y_{1_3},\dots ,y_{1_{2l+1}}]_R;\,\,\,\,\,\,\,\,\,\,\,\,\,\,\,\,\,\,\,\,\,\,\,\,\,\,\\
Type \,\, B_l\,(l\geq 2) :\,\,\bigwedge [y_{1_3},\dots ,y_{1_{4l-1}}]_R;\,\,\,\,\,\,\,\,\,\,\,\,\,\,\,\,\,\,\,\,\,\,\,\,\,\,\\
Type \,\, C_l\,(l\geq 3) :\,\,\bigwedge [y_{1_3},\dots ,y_{1_{4l-1}}]_R;\,\,\,\,\,\,\,\,\,\,\,\,\,\,\,\,\,\,\,\,\,\,\,\,\,\,\\
Type \,\, D_l\,(l\geq 4) :\,\,\bigwedge [y_{1_3},\dots
,y_{1_{4l-5}}]_R;\,\,\,\,\,\,\,\,\,\,\,\,\,\,\,\,\,\,\,\,\,\,\,\,\,\,\\
Type \,\, E_6 :\,\,
\bigwedge [y_{1_3},y_{1_9},y_{1_{11}},y_{1_{15}},y_{1_{17}},y_{1_{23}}]_R;\,\,\,\,\,\,\,\,\,\,\,\,\,\\
Type \,\, E_7 :\,\,
\bigwedge [y_{1_3},y_{1_{11}},y_{1_{15}},y_{1_{19}},y_{1_{23}},y_{1_{27}},y_{1_{35}}]_R;\\
Type \,\, E_8 :\,\,
\bigwedge [y_{1_3},y_{1_{15}},y_{1_{23}},y_{1_{35}},y_{1_{39}},y_{1_{47}},y_{1_{59}}]_R;\\
Type \,\, F_4 :\,\,
\bigwedge
[y_{1_3},y_{1_{11}},y_{1_{15}},y_{1_{23}}]_R;\,\,\,\,\,\,\,\,\,\,\,
\,\,\,\,\,\,\,\,\,\,\,\,\,\,\,\,\,\,\,\,\,\,\,\,\\
Type \,\, G_2 :\,\,\bigwedge [y_{1_3},y_{1_{11}}]_R.\,\,\,\,\,\,\,\,\,\,\,\,\,\,\,\,\,\,\,\,\,\,\,\,\,\,\,\,\,\,\,\,\,\,\,\,\,\,\,\,\,\,\,\,\,\,\,\,\,\,\,\,\,\,\,\,\,
\end{array}
\]
\end{thm}

\pf
We have a definable extension 
$$1\into Z(G)^0\into G\into G/Z(G)^0\into 1$$
and therefore,
$H^*(G;R)=H^*(Z(G)^0;R)\otimes _RH^*(G/Z(G)^0;R)$ and its enough to show
the theorem separately for $Z(G)^0$ and $G/Z(G)^0$.

We first prove the result for $Z(G)^0$. By
lemma  \ref{lem coincidence in G} we have that for
each $k\in \ZZ$, $(1-k)^{r}=
\lambda (p_k)=degp_{1-k}$ where $r=rank(Z(G)^0)$. 
Therefore for each $k>1$, $card(\{x\in Z(G)^0:kx=0\})=k^r$ (since
$p_k:Z(G)^0\into Z(G)^0$ does not change local $R$-orientation, recall that
$(p_k)^*(x)=
k^{len(x)}x$).
But by lemma \ref{lem torsion points} 
there are $s, m_1^{l_1}, \dots , m_u^{l_u}\in \NN $$
\cup \{0\}$ such that
$\pi _1(Z(G)^0)=\ZZ $$^s\oplus \ZZ $$_{m_1^{l_1}}\oplus \cdots \oplus \ZZ
$$_{m_u^{l_u}}$ 
and $card(\{x\in Z(G)^0:kx=0\})=
k^s\cdot (k, m_1^{l_1})\cdots (k, m_u^{l_u})$. It follows that $s=r$, 
$\pi _1(Z(G)^0)=\ZZ $$^r$ and
$H^*(Z(G)^0;R)=\bigwedge [y_{1_1},\dots ,y_{1_r}]_R$. 
But since $H^n(Z(G)^0;R)$ is non trivial ($Z(G)^0$ is $R$-orientable), 
we have $dimZ(G)^0=r$.

Let $d:=dim(G/Z(G)^0)$ and $t:=rank(G/Z(G)^0)$.
By \cite{e1} $G/Z(G)^0$ is definably semisimple definable group
and by \cite{pps1} and 
\cite{pps3}, $G/Z(G)^0$ is definably isomorphic to
an $N$-semialgebraic subgroup (definable with parameters from $\QQ$)
of some $GL(m,N)$ and therefore since
the cohomology groups depend only on the triangulation, 
$G/Z(G)^0$ and $(G/Z(G)^0)(\RR$$)$ have the same cohomology rings. But 
$(G/Z(G)^0)(\RR$$)$ is a Lie group and it follows from theorem
IV, volume II chapter IV in \cite{ghv} that $rank(G/Z(G)^0)$ is the 
dimension of a maximal definably connected definably compact definable abelian 
subgroup, $\sum _{i=1}^{t}g_i=d$ and $d\equiv t(mod 2).$ 

Also by results from
\cite{ghv} it follows together with what we have proved so far that
$H^*(G;R)=H^*(\mathbf{g}$$_I;R)$ where $\mathbf{g}$ is the
Lie algebra of $G$. 
On the other hand, by \cite{pps1} we have a definable extension
$$1\into Z(G/Z(G)^0)\into G/Z(G)^0\into G_1\times \cdots \times
G_k\into 1$$
where each $G_i$ is a definably simple centerless definably connected
definable group
(definable with parameters from $\QQ$) and each one of $G_i(\RR$$)$ is
a simple, compact, connected Lie group of one of the 
types $A_l,B_l,C_l,D_l,E_6,$ $E_7,E_8,$$F_4$ or $G_2$. Now the result
follows from a similar result for Lie groups (see \cite{mt} or \cite{ghv}):
For example, the classical groups $SU(l+1)$, $SO(2l+1)$,
$Sp(l)$ and $SO(2l)$ are of type $A_l, B_l, C_l$ and $D_l$ respectively,
$H^*(SU(n);R)$$=\bigwedge [y_{1_3},\dots ,y_{1_{2n-1}}]_R$,
$H^*(SO(2m+1);R)=\bigwedge [y_{1_3},\dots ,y_{1_{4m-1}}]_R$, 
we also have 
$H^*(Sp(n);R)$$=\bigwedge [y_{1_3},\dots ,y_{1_{4n-1}}]_R$ and 
$H^*(SO(2m)$$;R)$$=\bigwedge [y_{1_3},\dots ,$$y_{1_{4m-5}}]_R$. 

Since $H^1(G/Z(G)^0;R)=0$ we have that $\pi _1(G/Z(G)^0)$ is finite. 
\qed

Note that $H^*(G;R)=\otimes _R^{|G/G^0|}H^*(G^0;R)$. In particular,
$H^*(O(n);R)=\otimes _R^2H^*(SO(n);R)$.
The cohomology ring of $U(n)$ is given by
$H^*(U(n);R)=\bigwedge [y_{1_1},\dots ,y_{1_{2n-1}}]_R$,
since $Z(U(n))=U(1)=SO(2)$ and the projective unitary group 
$PU(n):=U(n)/Z(U(n))=SU(n)/\ZZ$$_n$.

\begin{cor}
Let $G$ be a definably compact, definably connected definable
group. Then $G$ is definably semisimple iff $Z(G)$ is finite iff
$Ad:G\into Ad(G)$ is a definable covering map iff $\pi _1(G)$ is
finite iff $H^1(G;\ZZ$$)=0$ iff the universal covering group
$\tilde{G}$ of $G$ is a definably compact (definably semisimple)
definable group.   
\end{cor}

\pf
This follows from theorem \ref{thm main} and results from \cite{e2}.
\qed

\begin{nrmk}\label{nrmk mod p cohomology}
{\em
The cohomology ring $H^*(G;\ZZ$$_p)$ of
a definably compact definably, definably simply connected, definably
simple  definable group $G$ is equal the cohomology
ring of the corresponding simply connected, compact simple Lie group
$G(\RR$$)$ of the same type as $G$. Explicit computation of these cohomology
rings is given
in \cite{mt}. We will not include here the full description of these
rings, but we note that
the (co)homology ring of definably simply connected, definably
compact, definably simple definable group $G$ are $p$-torsion free in
the following cases:
\[
\begin{array}{clcr}
\,\,\,\,\,\,\,\,\,\,p\geq 2,\,\,\,\,G\,\,of\,type\,\,A_l,\,\,C_l;\\
\,\,\,\,\,\,\,\,\,\,\,\,\,\,\,\,\,\,\,\,\,\,\,p\geq 3,\,\,\,\,G\,\,of\,type\,\,B_l,\,\,D_l,\,\,G_2;\\
\,\,\,\,\,\,\,\,\,\,\,\,\,\,\,\,\,\,\,\,\,\,\,p\geq 5,\,\,\,\,G\,\,of\,type\,\,F_4,\,\,E_6,\,\,E_7;\\
p\geq 7,\,\,\,\,G\,\,of\,type\,\,E_8.
\end{array}
\]

Finally note that if $G$ is a definably connected, definably compact 
definable abelian group of dimension $n$ then, since $
H_1(G)=\pi _1(G)=\ZZ $$^n$, $H^*(G)$ is $p$-torsion free for all $p$
and so $H^*(G;K_p)=
\bigwedge [x_{1_1},\dots ,x_{1_n}]_{K_p}$. 
}
\end{nrmk}

\begin{nrmk}
{\em
The following results are obtained by transfering similar results for
simple Lie groups (see \cite{mt}) and from result from \cite{e2}. 
Let $G$ be a definably connected,
definably compact, definably semisimple definable group with universal
covering group $\tilde{G}$. Then 
$$\pi _1(G)\leq Z(\tilde{G})=Z(\tilde{G_1})\times \cdots \times
Z(\tilde{G_k})$$
where $\tilde{G_i}$'s are definably simple, simply connected definably
compact definable groups such that
$\tilde{G}/Z(\tilde{G})=
\tilde{G_1}\times \cdots \times \tilde{G_k}.$

Moreover, 
$Z(\tilde{G_i})$ is $\ZZ$$_{l+1}$, $\ZZ$$_2$, $\ZZ$$_2$, $\ZZ$$_2\oplus
\ZZ$$_2$, $\ZZ$$_4$, $\ZZ$$_3$, $\ZZ$$_2$, $0$, $0$ or  $0$ if
$\tilde{G_i}(\RR$$)$
is of type $A_l$, $B_l$, $C_l$, $D_{2l}$, $D_{2l+1}$, $E_6$, $E_7$,
$E_8$, $F_4$ or $G_2$ respectively.

For example, $\pi _1(SU(n))=0$ and $Z(SU(n+1))=\ZZ$$_{n+1}$;
for $n>1$ we have $\pi _1(SO(n))=\ZZ$$_2$,
$Z(SO(2n+1))=0$, $Z(SO(2n))=\ZZ$$_2$
the universal covering group of $SO(n)$ is
$Spin(n)$ (the spinor groups) and $Z(Spin(2n+1))=\ZZ$$_2$;
$\pi _1(Sp(n))=0$ and $Z(Sp(n))=\ZZ$$_2$;
$Z(Spin(4n))=\ZZ$$_2\oplus \ZZ$$_2$ and $Z(Spin(4n+2))=\ZZ$$_4.$
}
\end{nrmk}

\end{section}

\end{document}